\documentclass[12pt, a4paper]{amsart}
\usepackage{pdfsync}
\usepackage{amsmath,amsthm,amsfonts,amscd,amssymb,eucal,latexsym,mathrsfs}
\usepackage[all]{xy}
\usepackage{latexsym, amssymb, amscd, mathrsfs, graphics, fullpage}
\usepackage[backref=page,colorlinks,bookmarksopen=true]{hyperref}

\def\gG{\mathfrak{G}}
\def\gH{\mathfrak{H}}

\def\gM{\mathfrak{M}}
\def\gN{\mathfrak{N}}
\def\gT{\mathfrak{T}}
\def\ga{\mathfrak{a}}
\def\gb{\mathfrak {b}}

\begin{document}

%%%%%%%%%%%%%%%%%%%%%TITLE
\title[Measured Quantum Groupoid action] 
{The Unitary Implementation of a Measured Quantum Groupoid action}
\author{Michel Enock}
\address{Institut de Math\'ematiques de Jussieu, Unit\'{e} Mixte Paris 6 / Paris 7 /
CNRS de Recherche 7586 \\175, rue du Chevaleret, Plateau 7E, F-75013 Paris}
 \email{enock@math.jussieu.fr}

\begin{abstract}
Mimicking the von Neumann version of Kustermans and Vaes' locally compact quantum groups, Franck Lesieur had introduced a notion of measured quantum groupoid, in the setting of von Neumann algebras. In a former article, the author had introduced the notions of actions, crossed-product, dual actions of a measured quantum groupoid; a biduality theorem for actions has been proved.  This article continues that program : we prove the existence of a standard implementation for an action, and a biduality theorem for weights.  We generalize this way results which were proved, for locally compact quantum groups by S. Vaes, and for measured groupoids by T. Yamanouchi. \end{abstract}

\maketitle
\newpage
 %%%%%%%%%%intro
\section{Introduction}
\label{intro}
\subsection{}
 In two articles (\cite{Val1}, \cite{Val2}), J.-M. Vallin has introduced two notions (pseudo-multiplicative
unitary, Hopf-bimodule), in order to generalize, up to the groupoid
case, the classical notions of multiplicative unitary \cite{BS} and of Hopf-von Neumann algebras \cite{ES}
which were introduced to describe and explain duality of groups, and leaded to appropriate notions
of quantum groups (\cite{ES}, \cite{W1}, \cite{W2}, \cite{BS}, \cite{MN}, \cite{W3}, \cite{KV1}, \cite{KV2}, \cite{MNW}). 
\\ In another article \cite{EV}, J.-M. Vallin and the author have constructed, from a depth 2 inclusion of
von Neumann algebras $M_0\subset M_1$, with an operator-valued weight $T_1$ verifying a regularity
condition, a pseudo-multiplicative unitary, which leaded to two structures of Hopf bimodules, dual
to each other. Moreover, we have then
constructed an action of one of these structures on the algebra $M_1$ such that $M_0$
is the fixed point subalgebra, the algebra $M_2$ given by the basic construction being then
isomorphic to the crossed-product. We construct on $M_2$ an action of the other structure, which
can be considered as the dual action.
\\  If the inclusion
$M_0\subset M_1$ is irreducible, we recovered quantum groups, as proved and studied in former papers
(\cite{EN}, \cite{E2}).
\\ Therefore, this construction leads to a notion of "quantum groupoid", and a construction of a
duality within "quantum groupoids". 
\subsection{}
In a finite-dimensional setting, this construction can be
mostly simplified, and is studied in \cite{NV1}, \cite{BSz1},
\cite{BSz2}, \cite{Sz},\cite{Val3}, \cite{Val4}, and examples are described. In \cite{NV2}, the link between these "finite quantum
groupoids" and depth 2 inclusions of
$II_1$ factors is given. 
\subsection{}
F. Lesieur, in \cite{L}, starting from a Hopf-bimodule, as introduced in \cite{Val1}, when there exist a left-invariant operator-valued weight, and a right-invariant operator-valued weight, mimicking in that wider setting the technics of Kustermans and Vaes (\cite{KV1}, \cite{KV2}), obtained a pseudo-multiplicative unitary, which, as in quantum group theory, "contains" all the information about the object (the von Neumann algebra, the coproduct) and allows to construct important data (an antipod, a co-inverse, etc.) Lesieur gave the name of "measured quantum groupoids" to these objects. A new set of axioms for these had been given in an appendix of \cite{E5}. In \cite{E4} had been shown that, with suitable conditions, the objects constructed from \cite{EV} are "measured quantum groupoids" in the sense of Lesieur. 
\subsection{}
In \cite{E5} have been developped the notions of action (already introduced in \cite{EV}), crossed-product, etc, following what had been done for locally compact quantum groups in (\cite{E1}, \cite{ES1}, \cite{V1}); a biduality theorem for actions had been obtained in (\cite{E5}, 11.6). Several points were left apart in \cite{E5}, namely the generalization of Vaes' theorem (\cite{V1}, 4.4) on the standard implementation of an action of a locally compact quantum group, which was the head light of \cite{V1}, and a biduality theorem for weights, as obtained in \cite{Y3}, \cite{Y4} (in fact, we were much more inspired by the shorter proof given in an appendix of \cite{BV}). 
\newline
We solve here these two problems when there exists a normal semi-finite faithful operator-valued weight from the von Neumann algebra on which the measured quantum groupoid is acting, onto the copy of the basis of this measured quantum groupoid which is put inside this algebra. In fact, these results appear much more as a biduality theorem of operator-valued weights rather than a biduality theorem on weights, which seems quite natural in the spirit of measured quantum groupoids, where, for instance, left-invariant weight on a locally compact quantum group is replaced by a left-invariant operator-valued weight. The strategy for the proofs had been mostly inspired by \cite{V1} and \cite{BV}. 
\subsection{}
This article is organized as follows : 
\newline
In chapter \ref{not}, we recall very quickly all the notations and results needed in that article; we have tried to make these preliminaries as short as possible, and we emphazise that this article should be understood as the continuation of \cite{E5}. 
\newline
In chapter \ref{standarddual}, we follow (\cite{V1}, 4.1 to 4.4), and prove, for any dual action, the result on the standard implementation of an action. 
\newline  
Chapter \ref{auxilliary} is rather technical; let $\gG=(N, M, \alpha, \beta, \Gamma, T, T', \nu)$ be a measured quantum groupoid, and let $b$ be an injective $*$-anti-homomorphism from $N$ into a von Neumann algebra $A$; let us suppose that there exists a normal semi-finite faithful operator-valued weight $\gT$ from $A$ onto $b(N)$, and let us write $\psi=\nu^o\circ b^{-1}\circ\gT$. Then, we can define on $A\underset{N}{_b*_\alpha}\mathcal L(H)$ a weight $\underline{\psi}$, which will generalize the tensor product of $\psi$ and $Tr_{\widehat{\Delta}^{-1}}$ (when $\gG$ is a locally compact quantum group, and therefore $N=\mathbb{C}$).  
\newline
In chapter \ref{using}, using this auxilliary weight introduced in chapter \ref{auxilliary}, and the particular case of the dual actions studied in chapter \ref{standarddual}, we calculate the standard implementation of an action, whenever there exists a normal semi-finite faithful operator-valued weight from $A$ onto $b(N)$. This condition is fulfilled trivially when the measured quantum groupoid is a locally compact quantum group, or is a measured groupoid; therefore, we recover in both cases the results already obtained. 
\newline
Chapter \ref{gamma} is another technical chapter; we define conditions on a weight $\psi$ defined on $A$ which allow us to construct on $A\underset{N}{_b*_\alpha}\mathcal L(H)$ a weight $\underline{\psi_\delta}$ which generalize the tensor product of $\psi$ and $Tr_{(\delta\widehat{\Delta})^{-1}}$(when $\gG$ is a locally compact quantum group, and therefore $N=\mathbb{C}$). 
\newline
In chapter \ref{bidualw} we use both auxilliary weights constructions made in chapters \ref{auxilliary} and \ref{gamma}; then, when there exists a normal semi-finite faithful operator-valued weight $\gT$ from $A$ onto $b(N)$ such that $\psi=\nu^o\circ b^{-1}\circ\gT$, we can define a Radon-Nikodym derivative of the weight $\psi$ with respect to the action, which will be a cocycle for this action. This condition is fulfilled trivially when the measured quantum groupoid is a locally compact quantum group, or is a measured groupoid, and, therefore, we recover in both cases the results already obtained. 
%%%%%%%notations
 \section{Definitions and notations}
 \label{not}
This article is the continuation of \cite{E5}; preliminaries are to be found in \cite{E5}, and we just recall herafter the following definitions and notations :

%%%%spatial
\subsection{Spatial theory; relative tensor products of Hilbert spaces and fiber products of von Neumann algebras (\cite{C1}, \cite{S}, \cite{T}, \cite{EV})}
\label{spatial}
 Let $N$ a von Neumann algebra, $\psi$ a normal semi-finite faithful weight on $N$; we shall denote by $H_\psi$, $\gN_\psi$, $S_\psi$, $J_\psi$, $\Delta_\psi$... the canonical objects of the Tomita-Takesaki theory associated to the weight $\psi$; let $\alpha$ be a non degenerate faithful representation of $N$ on a Hilbert space $\mathcal H$; the set of $\psi$-bounded elements of the left-module $_\alpha\mathcal H$ is :
\[D(_\alpha\mathcal{H}, \psi)= \lbrace \xi \in \mathcal{H};\exists C < \infty ,\| \alpha (y) \xi\|
\leq C \| \Lambda_{\psi}(y)\|,\forall y\in \gN_{\psi}\rbrace\]
Then, for any $\xi$ in $D(_\alpha\mathcal{H}, \psi)$, there exists a bounded operator
$R^{\alpha,\psi}(\xi)$ from $H_\psi$ to $\mathcal{H}$,  defined, for all $y$ in $\gN_\psi$ by :
\[R^{\alpha,\psi}(\xi)\Lambda_\psi (y) = \alpha (y)\xi\]
which intertwines the actions of $N$. 
\newline
If $\xi$, $\eta$ are bounded vectors, we define the operator product 
\[<\xi,\eta>_{\alpha,\psi} = R^{\alpha,\psi}(\eta)^* R^{\alpha,\psi}(\xi)\]
belongs to $\pi_{\psi}(N)'$, which, thanks to Tomita-Takesaki theory, will be identified to the opposite von Neumann algebra $N^o$. 
\newline
If now $\beta$ is a non degenerate faithful antirepresentation of $N$ on a Hilbert space $\mathcal K$, the relative tensor product $\mathcal K\underset{\psi}{_\beta\otimes_\alpha}\mathcal H$ is the completion of the algebraic tensor product $K\odot D(_\alpha\mathcal{H}, \psi)$ by the scalar product defined,  if $\xi_1$, $\xi_2$ are in $\mathcal{K}$, $\eta_1$, $\eta_2$ are in $D(_\alpha\mathcal{H},\psi)$, by the following formula :
\[(\xi_1\odot\eta_1 |\xi_2\odot\eta_2 )= (\beta(<\eta_1, \eta_2>_{\alpha,\psi})\xi_1 |\xi_2)\]
If $\xi\in \mathcal{K}$, $\eta\in D(_\alpha\mathcal{H},\psi)$, we shall denote $\xi\underset{\psi}{_\beta\otimes_\alpha}\eta$ the image of $\xi\odot\eta$ into $\mathcal K\underset{\psi}{_\beta\otimes_\alpha}\mathcal H$, and, writing $\rho^{\beta, \alpha}_\eta(\xi)=\xi\underset{\psi}{_\beta\otimes_\alpha}\eta$, we get a bounded linear operator from $\mathcal K$ into $\mathcal K\underset{\nu}{_\beta\otimes_\alpha}\mathcal H$, which is equal to $1_\mathcal K\otimes_\psi R^{\alpha, \psi}(\eta)$. 
\newline
Changing the weight $\psi$ will give a canonical isomorphic Hilbert space, but the isomorphism will not exchange elementary tensors !
\newline

We shall denote $\sigma_\psi$ the relative flip, which is a unitary sending $\mathcal{K}\underset{\psi}{_\beta\otimes_\alpha}\mathcal{H}$ onto $\mathcal{H}\underset{\psi^o}{_\alpha\otimes _\beta}\mathcal{K}$, defined, for any $\xi$ in $D(\mathcal {K}_\beta ,\psi^o )$, $\eta$ in $D(_\alpha \mathcal {H},\psi)$, by :
\[\sigma_\psi (\xi\underset{\psi}{_\beta\otimes_\alpha}\eta)=\eta\underset{\psi^o}{_\alpha\otimes_\beta}\xi\]
In $x\in \beta(N)'$, $y\in \alpha(N)'$, it is possible to define an operator $x\underset{\psi}{_\beta\otimes_\alpha}y$ on $\mathcal K\underset{\psi}{_\beta\otimes_\alpha}\mathcal H$, with natural values on the elementary tensors. As this operator does not depend upon the weight $\psi$, it will be denoted $x\underset{N}{_\beta\otimes_\alpha}y$. We can define a relative flip $\varsigma_N$ at the level of operators such that $\varsigma_N(x\underset{N}{_\beta\otimes_\alpha}y)=y\underset{N^o}{_\alpha\otimes_\beta}x$. If $P$ is a von Neumann algebra on $\mathcal H$, with $\alpha(N)\subset P$, and $Q$ a von Neumann algebra on $\mathcal K$, with $\beta(N)\subset Q$, then we define the fiber product $Q\underset{N}{_\beta*_\alpha}P$ as $\{x\underset{N}{_\beta\otimes_\alpha}y, x\in Q', y\in P'\}'$, and we get that $\varsigma_N(Q\underset{N}{_\beta*_\alpha}P)=P\underset{N^o}{_\alpha*_\beta}Q$. 
\newline
Moreover, this von Neumann algebra can be defined independantly of the Hilbert spaces on which $P$ and $Q$ are represented; if $(i=1,2)$, $\alpha_i$ is a faithful non degenerate homomorphism from $N$ into $P_i$, $\beta_i$ is a faithful non degenerate antihomomorphism from $N$ into $Q_i$, and $\Phi$ (resp. $\Psi$) an homomorphism from $P_1$ to $P_2$ (resp. from $Q_1$ to $Q_2$) such that $\Phi\circ\alpha_1=\alpha_2$ (resp. $\Psi\circ\beta_1=\beta_2$), then, it is possible to define an homomorphism $\Psi\underset{N}{_{\beta_1}*_{\alpha_1}}\Phi$ from $Q_1\underset{N}{_{\beta_1}*_{\alpha_1}}P_1$ into $Q_2\underset{N}{_{\beta_2}*_{\alpha_2}}P_2$. 
\newline
The operators $\theta^{\alpha, \psi}(\xi, \eta)=R^{\alpha, \psi}(\xi)R^{\alpha, \psi}(\eta)^*$, for all $\xi$, $\eta$ in $D(_\alpha\mathcal H, \psi)$, generates a weakly dense ideal in $\alpha(N)'$. Moreover, there exists a family $(e_i)_{i\in I}$ of vectors in $D(_\alpha\mathcal H, \psi)$ such that the operators $\theta^{\alpha, \psi}(e_i, e_i)$ are 2 by 2 orthogonal projections ($\theta^{\alpha, \psi}(e_i, e_i)$ being then the projection on the closure of $\alpha(N)e_i$). Such a family is called an orthogonal $(\alpha, \psi)$-basis of $\mathcal H$.

%%%%%%%MQG
\subsection{Measured quantum groupoids (\cite{L}, \cite{E5})}
\label{MQG}
 A measured quantum groupoid is an octuplet $\mathfrak {G}=(N, M, \alpha, \beta, \Gamma, T, T', \nu)$ such that (\cite{E5}, 3.8) :
\newline
(i) $(N, M, \alpha, \beta, \Gamma)$ is a Hopf-bimodule (as defined in \cite{E5}, 3.1), 
\newline
(ii) $T$ is a left-invariant normal, semi-finite, faithful operator valued weight $T$ from $M$ to $\alpha (N)$, 
\newline
(iii) $T'$ is a right-invariant normal, semi-finite, faithful operator-valued weight $T'$ from $M$ to $\beta (N)$, 
\newline
(iv) $\nu$ is normal semi-finite faitfull weight on $N$, which is relatively invariant with respect to $T$ and $T'$. 
\newline
We shall write $\Phi=\nu\circ\alpha^{-1}\circ T$, and $H=H_\Phi$, $J=J_\Phi$, and, for all $n\in N$, $\hat{\beta}(n)=J\alpha(n^*)J$, $\hat{\alpha}(n)=J\beta(n^*)J$.  The weight $\Phi$ will be called the left-invariant weight on $M$. 
\newline
Then, $\mathfrak {G}$ can be equipped with a pseudo-multiplicative unitary $W$ from $H\underset{\nu}{_\beta\otimes_\alpha}H$ onto $H\underset{\nu^o}{_\alpha\otimes_{\hat{\beta}}}H$ (\cite{E5}, 3.6), a co-inverse $R$, a scaling group $\tau_t$, an antipod $S$, a modulus $\delta$, a scaling operator $\lambda$, a managing operator $P$, and a canonical one-parameter group $\gamma_t$ of automorphisms on the basis $N$ (\cite{E5}, 3.8). Instead of $\mathfrak {G}$, we shall mostly use $(N, M, \alpha, \beta, \Gamma, T, RTR, \nu)$ which is another measured quantum groupoid, denoted $\underline{\mathfrak {G}}$, which is equipped with the same data ($W$, $R$, ...) as $\gG$. 
\newline
A dual measured quantum group $\widehat{\mathfrak{G}}$, which is denoted $(N, \widehat{M}, \alpha, \hat{\beta}, \widehat{\Gamma}, \widehat{T}, \widehat{R}\widehat{T}\widehat{R}, \nu)$, can be constructed, and we have $\widehat{\widehat{\mathfrak {G}}}=\underline{\mathfrak {G}}$. 
\newline
Canonically associated to $\mathfrak {G}$, can be defined also the opposite measured quantum groupoid is $\mathfrak{G}^o=(N^o, M, \beta, \alpha, \varsigma_N\Gamma, RTR, T, \nu^o)$ and the commutant measured quantum groupoid $\mathfrak{G}^c=(N^o, M', \hat{\beta}, \hat{\alpha}, \Gamma^c, T^c, R^cT^cR^c, \nu^o)$; we have $(\mathfrak{G}^o)^o=(\mathfrak{G}^c)^c=\underline{\mathfrak{G}}$, $\widehat{\mathfrak{G}^o}=(\widehat{\mathfrak {G}})^c$, $\widehat{\mathfrak {G}^c}=(\widehat{\mathfrak {G}})^o$, and $\mathfrak{G}^{oc}=\mathfrak {G}^{co}$ is canonically isomorphic to $\underline{\mathfrak {G}}$ (\cite{E5}, 3.12). 
\newline
The pseudo-multiplicative unitary of $\widehat{\mathfrak{G}}$ (resp. $\mathfrak{G}^o$, $\mathfrak{G}^c$) will be denoted $\widehat{W}$ (resp. $W^o$, $W^c$). The left-invariant weight on $\widehat{\mathfrak{G}}$ (resp. $\mathfrak{G}^o$, $\mathfrak{G}^c$) will be denoted $\widehat{\Phi}$ (resp. $\Phi^o$, $\Phi^c$). 
\newline
Let $_a\gH_b$ be a $N-N$-bimodule, i.e. an Hilbert space $\gH$ equipped with a normal faithful non degenerate representation $a$ of $N$ on $\gH$ and a normal faithful non degenerate anti-representation $b$ on $\gH$, such that $b(N)\subset a(N)'$. A corepresentation of $\gG$ on $_a\gH_b$ is a unitary $V$ from $\gH\underset{\nu^o}{_a\otimes_\beta}H_\Phi$ onto 
 $\gH\underset{\nu}{_b\otimes_\alpha}H_\Phi$, satisfying, for all $n\in N$ :
 \[V(b(n)\underset{N^o}{_a\otimes_\beta}1)=(1\underset{N}{_b\otimes_\alpha}\beta(n))V\]
 \[V(1\underset{N^o}{_a\otimes_\beta}\alpha(x))=(a(n)\underset{N}{_b\otimes_\alpha}1)V\]
such that, for any $\xi\in D(_a\gH, \nu)$ and $\eta\in D(\gH_b, \nu^o)$, the operator $(\omega_{\xi, \eta}*id)(V)$ belongs to $M$ (then, it is possible to define $(id*\theta)(V)$, for any $\theta$ in $M_*^{\alpha, \beta}$ which is the linear set generated by the $\omega_\xi$, with $\xi\in D(_\alpha H, \nu)\cap D(H_\beta, \nu^o)$), and such that the application $\theta\rightarrow (id*\theta)(V)$ from $M_*^{\alpha, \beta}$ into $\mathcal L(\gH)$ is multiplicative (\cite{E5} 5.1, 5.5).

%%%%%%%%%%action
\subsection{Action of a measured quantum groupoid (\cite{E5})}
\label{action}

An action (\cite{E5}, 6.1) of $\mathfrak{G}$ on a von Neumann algebra $A$ is a couple $(b, \mathfrak a)$, where :
\newline
(i) $b$ is an injective $*$-antihomomorphism from $N$ into $A$; 
\newline
(ii) $\mathfrak a$ is an injective $*$-homomorphism from $A$ into $A\underset{N}{_b*_\alpha}M$; 
\newline
(iii) $b$ and $\mathfrak a$ are such that, for all $n$ in $N$:
\[\mathfrak a (b(n))=1\underset{N}{_b\otimes_\alpha}\beta(n)\]
(which allow us to define $\mathfrak a\underset{N}{_b*_\alpha}id$ from $A\underset{N}{_b*_\alpha}M$ into $A\underset{N}{_b*_\alpha}M\underset{N}{_\beta*_\alpha}M$)
and such that :
\[(\mathfrak a\underset{N}{_b*_\alpha}id)\mathfrak a=(id\underset{N}{_b*_\alpha}\Gamma)\mathfrak a\]
The set of invariants is defined as the sub von Neumann algebra :
\[A^\ga=\{x\in A\cap b(N)', \ga(x)=x\underset{N}{_b\otimes_\alpha}1\}\]
If the von Neumann algebra acts on a Hilbert space $\gH$, and if there exists a representation $a$ of $N$ on $\gH$ such that $b(N)\subset A\subset a(N)'$, a corepresentation $V$ of $\gG$ on the bimodule $_a\gH_b$ will be called an implementation of $\ga$ if we have $\ga(x)=V(x\underset{N^o}{_a\otimes_b}1)V^*$ , for all $x\in A$ (\cite{E5}, 6.6); we shall look at the following more precise situation : let $\psi$ is a normal semi-finite faithful weight on $A$, and $V$ an implementation of $\ga$ on $_a(H_\psi)_b$ (with $a(n)=J_\psi b(n^*)J_\psi$), such that :
\[V^*=(J_\psi\underset{\nu^o}{_\alpha\otimes_\beta}J_{\widehat{\Phi}})V(J_\psi\underset{\nu}{_b\otimes_\alpha} J_{\widehat{\Phi}})\]
\newline
Such an implementation had been constructed (\cite{E5} 8.8) in the particular case when the weight $\psi$ is called $\delta$-invariant, which means that,  for all $\eta\in D(_\alpha H_\Phi, \nu)\cap\mathcal D(\delta^{1/2})$, such that $\delta^{1/2}\eta$ belongs to $D((H_\Phi)_\beta, \nu^o)$, and for all $x\in\gN_\psi$, we have:
\[\psi((id\underset{N}{_b*_\alpha}\omega_\eta)\mathfrak a(x^*x))=\|\Lambda_\psi(x)\underset{\nu^o}{_a\otimes_\beta}\delta^{1/2}\eta\|^2\]
and bears the density property, which means that the subset $D((H_\psi)_b, \nu^o)\cap D(_aH_\psi, \nu)$ is dense in $H_\psi$. This standard implementation is then given by the formula (\cite{E5}, 8.4) :
\[V_\psi(\Lambda_\psi (x)\underset{\nu^o}{_a\otimes_\beta}\delta^{1/2}\eta)=\sum_i\Lambda_\psi((id\underset{N}{_b*_\alpha}\omega_{\eta, e_i})\mathfrak a(x))\underset{\nu}{_b\otimes_\alpha}e_i\]
for all $x\in\gN_\psi$, $\eta\in D(_\alpha H, \nu)\cap\mathcal D(\delta^{1/2})$ such that $\delta^{1/2}\eta$ belongs to $D(H_\beta, \nu^o)$, $(e_i)_{i\in I}$ any orthonormal $(\alpha, \nu)$-basis of $H$. Moreover (\cite{E5}, 8.9), it is possible to define one parameter groups of unitaries $\Delta_\psi^{it}\underset{N^o}{_a\otimes_\beta}\delta^{-it}\Delta_{\widehat{\Phi}}^{-it}$ and $\Delta_\psi^{it}\underset{N}{_b\otimes_\alpha}\delta^{-it}\Delta_{\widehat{\Phi}}^{-it}$, with natural values on elementary tensor, and we have :
\[V_\psi(\Delta_\psi^{it}\underset{N^o}{_a\otimes_\beta}\delta^{-it}\Delta_{\widehat{\Phi}}^{-it})=(\Delta_\psi^{it}\underset{N}{_b\otimes_\alpha}\delta^{-it}\Delta_{\widehat{\Phi}}^{-it})V_\psi\]
and, therefore, for any $x$ in $A$, $t$ in $\mathbb{R}$, we have :
\[\mathfrak a(\sigma_t^\psi(x))=(\Delta_\psi^{it}\underset{N}{_b\otimes_\alpha}\delta^{-it}\Delta_{\widehat{\Phi}}^{-it})\mathfrak a(x)(\Delta_\psi^{-it}\underset{N}{_b\otimes_\alpha}\delta^{it}\Delta_{\widehat{\Phi}}^{it})\]

%%%%%%crossed
\subsection{Crossed-product (\cite{E5})}
\label{crossed}
The crossed-product of $A$ by $\mathfrak {G}$ via the action $\mathfrak a$ is the von Neumann algebra generated by $\mathfrak a(A)$ and $1\underset{N}{_b\otimes_\alpha}\widehat{M}'$ (\cite{E5}, 9.1) and is denoted $A\rtimes_\mathfrak a\mathfrak {G}$; then there exists (\cite{E5}, 9.3) an action $(1\underset{N}{_b\otimes_\alpha}\hat{\alpha}, \tilde{\mathfrak a})$ of $(\widehat{\mathfrak {G}})^c$ on $A\rtimes_\mathfrak a\mathfrak {G}$. 
\newline
The biduality theorem (\cite{E5}, 11.6) says that the bicrossed-product $(A\rtimes_\mathfrak a\mathfrak {G})\rtimes_{\tilde{\mathfrak a}}\widehat{\mathfrak {G}}^c$ is canonically isomorphic to $A\underset{N}{_b*_\alpha}\mathcal L(H)$; more precisely, this isomorphism is given by :
\[\Theta (\ga\underset{N}{_b*_\alpha}id)(A\underset{N}{_b*_\alpha}\mathcal L(H))=(A\rtimes_\mathfrak a\mathfrak {G})\rtimes_{\tilde{\mathfrak a}}\widehat{\mathfrak {G}}^c\]
where $\Theta$ is the spatial isomorphism between $\mathcal L(\gH\underset{\nu}{_b\otimes_\alpha}H\underset{\nu}{_\beta\otimes_\alpha}H)$ and $\mathcal L(\gH\underset{\nu}{_b\otimes_\alpha}H\underset{\nu^o}{_{\hat{\alpha}}\otimes_\beta}H)$ implemented by $1_\gH\underset{\nu}{_b\otimes_\alpha}\sigma_\nu W^o\sigma_\nu$; the biduality theorem says also that this isomorphism sends  the action $(1\underset{N}{_b\otimes_\alpha}\hat{\beta}, \underline{\mathfrak a})$ of $\gG$ on $A\underset{N}{_b*_\alpha}\mathcal L(H)$, defined, for any $X\in A\underset{N}{_b*_\alpha}\mathcal L(H)$, by :
\[\underline{\mathfrak a}(X)=(1\underset{N}{_b\otimes_\alpha}\sigma_{\nu^o}W\sigma_{\nu^o})(id\underset{N}{_b*_\alpha}\varsigma_N)(\mathfrak a\underset{N}{_b*_\alpha}id)(X)(1\underset{N}{_b\otimes_\alpha}\sigma_{\nu^o}W\sigma_{\nu^o})^*\]
on the bidual action (of $\mathfrak{G}^{co}$) on $(A\rtimes_\mathfrak a\mathfrak {G})\rtimes_{\tilde{\mathfrak a}}\widehat{\mathfrak {G}}^o$. 
\newline
There exists a normal faithful semi-finite operator-valued weight $T_{\tilde{\ga}}$ from $A\rtimes_\mathfrak a\mathfrak {G}$ onto $\ga(A)$; therefore, starting with a  normal semi-finite weight $\psi$ on $A$, we can construct a dual weight $\tilde{\psi}$ on $A\rtimes_\mathfrak a\mathfrak {G}$ by the formula $\tilde{\psi}=\psi\circ\ga^{-1}\circ T_{\tilde{\ga}}$ (\cite{E5} 13.2). These dual weights are exactly the $\hat{\delta}^{-1}$-invariant weights on $A\rtimes_\mathfrak a\mathfrak {G}$ bearing the density property (\cite{E5} 13.3). 
\newline
Moreover (\cite{E5} 13.3), the linear set generated by all the elements $(1\underset{N}{_b\otimes_\alpha}a)\mathfrak a(x)$, for all $x\in\gN_\psi$, $a\in\gN_{\widehat{\Phi}^c}\cap\gN_{\hat{T}^c}$, is a core for $\Lambda_{\tilde{\psi}}$, and it is possible to identify the GNS representation of $A\rtimes_\mathfrak a\gG$ associated to the weight $\tilde{\psi}$ with the natural representation on $H_\psi\underset{\nu}{_b\otimes_\alpha}H_\Phi$ by writing :
\[\Lambda_\psi(x)\underset{\nu}{_b\otimes_\alpha}\Lambda_{\widehat{\Phi}^c}(a)=\Lambda_{\tilde{\psi}}[(1\underset{N}{_b\otimes_\alpha}a)\mathfrak a(x)]\]
which leads to the identification of $H_{\tilde{\psi}}$ with $H_\psi\underset{\nu}{_b\otimes_\alpha}H$. Moreover, using that identification, the linear set generated by the elements of the form $\mathfrak a(y^*)(\Lambda_\psi(x)\underset{\nu}{_b\otimes_\alpha}\Lambda_{\widehat{\Phi}^c}(a))$, for $x, y$ in $\gN_\psi$, and $a$ in $\gN_{\widehat{\Phi}^c}\cap\gN_{\hat{T}^c}\cap\gN_{\widehat{\Phi}^c}^*\cap\gN_{\hat{T}^c}^*$ is a core for $S_{\tilde{\psi}}$, and we have :
\[S_{\tilde{\psi}}\mathfrak a(y^*)(\Lambda_\psi(x)\underset{\nu}{_b\otimes_\alpha}\Lambda_{\widehat{\Phi}^c}(a))=\mathfrak a(x^*)(\Lambda_\psi(y)\underset{\nu}{_b\otimes_\alpha}\Lambda_{\widehat{\Phi}^c}(a^*))\]
Then, the unitary $U_\psi^\ga=J_{\tilde{\psi}}(J_\psi\underset{N^o}{_a\otimes_\beta}J_{\widehat{\Phi}})$ from $H_\psi\underset{\nu^o}{_a\otimes_\beta}H_\Phi$ onto $H_\psi\underset{\nu}{_b\otimes_\alpha}H_\Phi$ satisfies :
\[U^\mathfrak a_\psi(J_\psi\underset{N}{_b\otimes_\alpha}J_{\widehat{\Phi}})=(J_\psi\underset{N}{_b\otimes_\alpha}J_{\widehat{\Phi}})(U^\mathfrak a_\psi)^*\]
and we have (\cite{E5} 13.4) :
\newline
(i) for all $y\in A$ :
\[\mathfrak a (y)=U^\mathfrak a_\psi(y\underset{N^o}{_a\otimes_\beta}1)(U^\mathfrak a_\psi)^*\]
(ii) for all $b\in M$ :
\[(1\underset{N}{_b\otimes_\alpha}J_\Phi bJ_\Phi)U^\mathfrak a_\psi=U^\mathfrak a_\psi(1\underset{N^o}{_a\otimes_\beta}J_\Phi bJ_\Phi)\]
(iii) for all $n\in N$ :
\[U_\psi^\mathfrak a(b(n)\underset{N^o}{_a\otimes_\beta}1)=(1\underset{N}{_b\otimes_\alpha}\beta(n))U_\psi^\mathfrak a\]
\[U_\psi^\mathfrak a(1\underset{N^o}{_a\otimes_\beta}\alpha(n))=(a(n)\underset{N}{_b\otimes_\alpha}1)U_\psi^\mathfrak a\]
Therefore, we see that this unitary $U^\ga_\psi$ "implements" $\ga$, but we do not know whether it is a corepresentation. If it is, we shall say that it is a standard implemantation of $\ga$. 
\newline
We can define the bidual weight $\tilde{\tilde{\psi}}$ on $(A\rtimes_\mathfrak a\mathfrak {G})\rtimes_{\tilde{\mathfrak a}}\widehat{\mathfrak {G}}^o$, and the weight $\tilde{\tilde{\psi}}\circ\Theta\circ (\ga\underset{N}{_b*_\alpha}id)$ on $A\underset{N}{_b*_\alpha}\mathcal L(H)$, that we shall denote $\overline{\psi_\ga}$ for simplification (or $\overline{\psi}$ if there is no ambiguity about the action). Then we get (\cite{E5}, 13.6) that the spatial derivative $\frac{d\overline{\psi}}{d\psi^o}$ is equal to the modulus operator $\Delta_{\tilde{\psi}}$. There exists a normal semi-finite faithful operator-valued weight  $T_{\underline{\ga}}$ from $A\underset{N}{_b*_\alpha}\mathcal L(H)$ onto $A\rtimes_\ga\gG$ such that $\overline{\psi_\ga}=\tilde{\psi}\circ T_{\underline{\ga}}$
\newline
Using twice (\cite{T} 4.22(ii)), we obtain, for any $x\in A$ and $t\in\mathbb{R}$, that $\sigma_t^{\overline{\psi_\ga}}(\ga(x))=\ga(\sigma_t^\psi(x))$; and if $\psi_1$ and $\psi_2$ are two normal semi-finite faithful weights on $A$, , we get :
\[(D\overline{\psi_{1\ga}}:D\overline{\psi_{2\ga}})_t=(D\tilde{\psi_1} : D\tilde{\psi_2})_t=\ga((D\psi_1:D\psi_2)_t)\]

%%%ex
\subsection{Examples of measured quantum groupoids}
\label{ex}
Examples of measured quantum groupoids are the following :
\newline
(i) locally compact quantum groups, as defined and studied by J. Kustermans and S. Vaes (\cite{KV1}, \cite {KV2}, \cite{V1}); these are, trivially, the measured quantum groupoids with the basis $N=\mathbb{C}$. 
\newline
(ii) measured groupoids, equipped with a left Haar system and a quasi-invariant measure on the set of units, as studied mostly by T. Yamanouchi (\cite{Y1}, \cite{Y2}, \cite{Y3}, \cite{Y4}); it was proved in \cite{E6} that these measured quantum groupoids are exactly those whose underlying von Neumann algebra is abelian. 
\newline
(iii) the finite dimensional case had been studied by D. Nikshych and L. Vainermann (\cite{NV1}, \cite{NV2}, \cite{NV3}), J.-M. Vallin (\cite{Val3}, \cite{Val4}) and M.-C. David (\cite{D}); in that case, non trivial examples are given, for instance Temperley-Lieb algebras (\cite{NV3}, \cite{D}), which had appeared in subfactor theory (\cite{J}).  . 
\newline
(iv) continuous fields of ($\bf{C}^*$-version of) locally compact quantum groups, as studied by E. Blanchard in (\cite{Bl1}, \cite{Bl2}); it was proved in \cite{E6} that these measured quantum groupoids are exactly those whose basis is central in the underlying von Neumann algebras of both the measured quantum groupoid and its dual. 
\newline
(v) in \cite{DC}, K. De Commer proved that, in the case of a monoidal equivalence between two locally compact quantum groups (which means that these two locally compact quantum group have commuting ergodic and integrable actions on the same von Neumann algebra), it is possible to construct a measurable quantum groupoid of basis $\mathbb{C}^2$ which contains all the data. Moreover, this construction was usefull to prove new results on locally compact quantum groups, namely on the deformation of a locally compact quantum group by a unitary $2$-cocycle; he proved that these measured quantum groupoids are exactly those whose basis $\mathbb{C}^2$ is central in the underlying von Neumann algebra of the measured quatum groupoid, but not in the underlying von Neumann algebra of the dual measured quantum groupoid. 
\newline
(vi) starting from a depth 2 inclusion $M_0\subset M_1$ of von Neumann algebras, equipped with an operator-valued weight $T_1$ from $M_1$ onto $M_0$, satisfying appropriate conditions, such that there exists a normal semi-finite faithful weight $\chi$ on the first relative commutant $M'_0\cap M_1$, invariant under the modular automorphism group $\sigma_t^{T_1}$, it has been proved (\cite{EV}, \cite{E4}) that it is possible to put on the second relative commutant $M'_0\cap M_2$ (where $M_0\subset M_1\subset M_2\subset M_3 ...$ is Jones' tower associated to the inclusion $M_0\subset M_1$) a canonical structure of a measured quantum groupoid; moreover, its dual is given then by the same construction associated to the inclusion $M_1\subset M_2$, and this dual measured quantum groupoid acts canonically on the von Neumann algebra $M_1$, in such a way that $M_0$ is equal to the subalgebra of invariants, and the inclusion $M_1\subset M_2$ is isomorphic to the inclusion of $M_1$ into its crossed-product. This gives a "geometrical" construction of measured quantum groupoids; in another article in preparation (\cite{E7}), in which is used the biduality theorem for weights proved in  \ref{cocycle}, had been proved that any measured quantum groupoid has an outer action on some von Neumann algebra, and can be, therefore, obtained by this "geometrical construction". The same result for locally compact quantum groups relies upon \cite{V2} and the corresponding result for measured quantum groupoids had been pointed out in \cite{E5}. 
\newline
(vii) in \cite{VV} and \cite{BSV} was given a specific procedure for constructing locally compact quantum groups, starting from a locally compact group $G$, whose almost all elements belong to the product $G_1G_2$ (where $G_1$ and $G_2$ are closed subgroups of $G$ whose intersection is reduced to the unit element of $G$); such $(G_1, G_2)$ is called a "matched pair" of locally compact groups (more precisely, in \cite{VV}, the set $G_1G_2$ is required to be open, and it is not the case in \cite{BSV}).Then, $G_1$ acts naturally on $L^\infty(G_2)$ (and vice versa), and the two crossed-products obtained bear the structure of two locally compact quantum groups in duality. In \cite{Val5}, J.-M. Vallin generalizes this constructions up to groupoids, and, then, obtains examples of measured quantum groupoids; more specific examples are then given by the action of a matched pair of groups on a locally compact space, and also more exotic examples.

%%%%%standarddual
\section{The standard implementation of an action : the case of a dual action}
\label{standarddual}
In this chapter, following \cite{V1}, we prove that the unitary $U_\psi^\ga$ introduced in \ref{crossed} is a standard implementation of $\ga$, for all normal semi-finite faithful weight $\psi$ on $A$, whenever $\ga$ is a dual action (\ref{cordual}). For this purpose, we prove first that, if for some weight $\psi_1$, the unitary $U^\ga_{\psi_1}$ is a standard implementation, then, for any weight $\psi$, $U^\ga_\psi$ is a standard implementation (\ref{propu}). Second (\ref{Uinv}), we prove, for a $\delta$-invariant weight $\psi$, that $U_\psi^\ga$ is equal to the implementation $V_\psi$ constructed in (\cite{E5} 8.8) and recalled in \ref{action}. Thanks to (\cite{E5} 13.3), recalled in \ref{crossed}, we then get the result.

%%%%%propu
\subsection{Proposition}
\label{propu}
{\it Let $\gG$ be a measured quantum groupoid, and $(b,\mathfrak a)$ an action of $\gG$ on a von Neumann algebra $A$; let $\psi_1$ and $\psi_2$ be two normal faithful semi-finite weights on $A$ and $U^\mathfrak a_{\psi_1}$ and $U^\mathfrak a_{\psi_2}$ the two unitaries constructed in \ref{crossed}; let $u$ be the unitary from $H_{\psi_1}$ onto $H_{\psi_2}$ intertwining the representations $\pi_{\psi_1}$ and $\pi_{\psi_2}$; then :
\newline 
(i) the unitary $u\underset{N}{_b\otimes_\alpha}1$ intertwines the representations of $A\rtimes_\mathfrak a\gG$ on $H_{\psi_1}\underset{\nu}{_b\otimes_\alpha}H_\Phi$ and on $
H_{\psi_2}\underset{\nu}{_b\otimes_\alpha}H_\Phi$; moreover, we have :
\[(u\underset{N}{_b\otimes_\alpha}1)U^\mathfrak a_{\psi_1}=U^\mathfrak a_{\psi_2}(u\underset{N^o}{_{a_1}\otimes_\beta}1)\]
where $a_1(n)=J_{\psi_1}\pi_{\psi_1}(b(n^*))J_{\psi_1}$, for all $n\in N$.
\newline
(ii) if $U^\ga_{\psi_1}$ is a corepresentation of $\gG$ on $H_{\psi_1}$, then $U^\ga_{\psi_2}$ is a corepresentation of $\gG$ on $H_{\psi_2}$. 
\newline
(iii) if $U^\ga_{\psi_1}$ is a standard implementation of $\ga$, then $U^\ga_{\psi_2}$ is a 
standard implementation of $\ga$. 
 }

\begin{proof}
Let us write $J_{2,1}$ the relative modular conjugation, which is an antilinear surjective isometry from $H_{\psi_1}$ onto $H_{\psi_2}$. Then we have $u=J_{2,1}J_{\psi_1}=J_{\psi_2}J_{2,1}$, by (\cite{St} 3.16). 
Moreover, let us define, for $x\in A$, and $t\in\mathbb{R}$ $\sigma_t^{2,1}(x)=[D\psi_2:D\psi_1]_t\sigma_t^{\psi_1}(x)$; then, by (\cite{St}, 3.15), for $x\in\gN_{\psi_1}$, $y\in D(\sigma^{2,1}_{-i/2})$, $xy^*$ belongs to $\gN_{\psi_2}$ and :
\[\Lambda_{\psi_2}(xy^*)=J_{2,1}\pi_{\psi_1}(\sigma^{2,1}_{-i/2}(y))J_{\psi_1}\Lambda_{\psi_1}(x)\]
Therefore, if $a\in\gN_{\widehat{\Phi}^c}$, $(1\underset{N}{_b\otimes_\alpha}a)\mathfrak a(xy^*)$ belongs to $\gN_{\tilde{\psi_2}}$, and, we have, where $V_i$ ($i=(1,2)$) denotes the unitary from $H_{\psi_i}\underset{\nu}{_b\otimes_\alpha}H_\Phi$ onto $H_{\tilde{\psi_i}}$ defined in \ref{crossed} :
\begin{eqnarray*}
\Lambda_{\tilde{\psi_2}}[(1\underset{N}{_b\otimes_\alpha}a)\mathfrak a(xy^*)]
&=&V_2(\Lambda_{\psi_2}(xy^*)\underset{\nu}{_b\otimes_\alpha}\Lambda_{\widehat{\Phi}^c}(a))\\
&=&V_2J_{2,1}\pi_{\psi_1}(\sigma^{2,1}_{-i/2}(y))J_{\psi_1}\Lambda_{\psi_1}(x)\underset{\nu}{_b\otimes_\alpha}\Lambda_{\widehat{\Phi}^c}(a))
\end{eqnarray*}
which is equal to :
\[V_2(J_{2,1}\pi_{\psi_1}(\sigma^{2,1}_{-i/2}(y))J_{\psi_1}\underset{N}{_b\otimes_\alpha}1)V_1^*\Lambda_{\tilde{\psi_1}}[(1\underset{N}{_b\otimes_\alpha}a)\mathfrak a(x)]\]
and, as the linear set generated by the elements of the form $(1\underset{N}{_b\otimes_\alpha}a)\mathfrak a(x)$ is a core for $\Lambda_{\tilde{\psi_1}}$, we get, for any $z\in\gN_{\tilde{\psi_1}}$, that $z\mathfrak a(y^*)$ belongs to $\gN_{\tilde{\psi_2}}$, and that :
\[\Lambda_{\tilde{\psi_2}}(z\mathfrak a (y^*))=V_2(J_{2,1}\pi_{\psi_1}(\sigma^{2,1}_{-i/2}(y))J_{\psi_1}\underset{N}{_b\otimes_\alpha}1)V_1^*\Lambda_{\tilde{\psi_1}}(z)\]
Let us denote by $\tilde{J_{2,1}}$ the relative modular conjugation constructed from the weights $\tilde{\psi_1}$ and $\tilde{\psi_2}$, and $\tilde{\sigma}^{2,1}_t$ the one-parameter group of isometries of $A\rtimes_\mathfrak a \gG$ constructed from these two weights by the formula, for any $X\in A\rtimes_\mathfrak a \gG$ :
\[\tilde{\sigma}^{2,1}_t(X)=[D\tilde{\psi_2}:D\tilde{\psi_1}]_t\sigma_t^{\tilde{\psi_1}}(X)\]
Using (\cite{St}, 3.15) applied to these two weights, we get that $\mathfrak a(y)$ belongs to $D(\tilde{\sigma}^{2,1}_{-i/2})$ and that :
\[\tilde{J_{2,1}}\pi_{\tilde{\psi_1}}(\tilde{\sigma}^{2,1}_{-i/2}(\mathfrak a (y)))J_{\tilde{\psi_1}}=
V_2(J_{2,1}\pi_{\psi_1}(\sigma^{2,1}_{-i/2}(y))J_{\psi_1}\underset{N}{_b\otimes_\alpha}1)V_1^*\]
We easily get that $\tilde{\sigma}^{2,1}_t(\mathfrak a(y))=\mathfrak a (\sigma^{2,1}_t(y))$ and, therefore, we have :
\[\pi_{\tilde{\psi_1}}(\mathfrak a(\sigma^{2,1}_{-i/2}(y))=\tilde{J_{2,1}}^*V_2(J_{2,1}\pi_{\psi_1}(\sigma^{2,1}_{-i/2}(y))J_{\psi_1}\underset{N}{_b\otimes_\alpha}1)V_1^*J_{\tilde{\psi_1}}\]
As we have, using \ref{crossed} :
\[(J_{\psi_1}\underset{N}{_b\otimes_\alpha}J_{\widehat{\Phi}})V_1^*J_{\tilde{\psi_1}}=U_{\psi_1}^\mathfrak a V_1^*\]
we get :
\[\pi_{\tilde{\psi_1}}(\mathfrak a(\sigma^{2,1}_{-i/2}(y))=
\tilde{J_{2,1}}^*V_2(J_{2,1}\underset{N^o}{_{a_1}\otimes_\beta}J_{\widehat{\Phi}})(\pi_{\psi_1}(\sigma^{2,1}_{-i/2}(y))\underset{N^o}{_{a_1}\otimes_\beta}1)U_{\psi_1}^\mathfrak a V_1^*\]
and, therefore, using \ref{crossed} :
\begin{eqnarray*}
\tilde{J_{2,1}}^*V_2(J_{2,1}\underset{N^o}{_{a_1}\otimes_\beta}J_{\widehat{\Phi}})(\pi_{\psi_1}(\sigma^{2,1}_{-i/2}(y))\underset{N^o}{_{a_1}\otimes_\beta}1)
&=&
\pi_{\tilde{\psi_1}}(\mathfrak a(\sigma^{2,1}_{-i/2}(y))V_1(U_{\psi_1}^\mathfrak a)^*\\
&=&V_1\mathfrak a(\sigma^{2,1}_{-i/2}(y))(U_{\psi_1}^\mathfrak a)^*
\end{eqnarray*}
which, using \ref{crossed}, is equal to :
\[V_1U_{\psi_1}^\mathfrak a(\pi_{\psi_1}(\sigma^{2,1}_{-i/2}(y))\underset{N^o}{_{a_1}\otimes_\beta}1)\]
By density, we get :
\[U_{\psi_1}^\mathfrak a=V_1^*\tilde{J_{2,1}}^*V_2(J_{2,1}\underset{N^o}{_{a_1}\otimes_\beta}J_{\widehat{\Phi}})\]
and, therefore, using \ref{crossed} again :
\begin{eqnarray*}
1_{H_{\psi_1}}\underset{N}{_b\otimes_\alpha}1_{H_\Phi}&=&
V_1^*\tilde{J_{2,1}}^*V_2(J_{2,1}\underset{N^o}{_{a_1}\otimes_\beta}J_{\widehat{\Phi}})
(J_{\psi_1}\underset{N}{_b\otimes_\alpha}J_{\widehat{\Phi}})V_1^*J_{\tilde{\psi_1}}V_1\\
&=&
V_1^*\tilde{J_{2,1}}^*V_2(u\underset{N}{_b\otimes_\alpha}1)V_1^*J_{\tilde{\psi_1}}V_1
\end{eqnarray*}
which implies that :
\[1_{H_{\tilde{\psi_1}}}\underset{N}{_b\otimes_\alpha}1_{H_\Phi}=\tilde{J_{2,1}}^*V_2(u\underset{N}{_b\otimes_\alpha}1)V_1^*J_{\tilde{\psi_1}}\]
and :
\[V_2(u\underset{N}{_b\otimes_\alpha}1)V_1^*=\tilde{J_{2,1}}J_{\tilde{\psi_1}}\]
But $\tilde{J_{2,1}}J_{\tilde{\psi_1}}=J_{\tilde{\psi_2}}\tilde{J_{2,1}}$ is the unitary from $H_{\tilde{\psi_1}}$ onto $H_{\tilde{\psi_2}}$ which intertwines $\pi_{\tilde{\psi_1}}$ and $\pi_{\tilde{\psi_2}}$; from which we get the first result. 
\newline
This formula gives also, where $a_2(n)=J_{\psi_2}\pi_{\psi_2}(b(n^*))J_{\psi_2}$, for all $n\in N$ :
\begin{eqnarray*}
U_{\psi_2}^\mathfrak a
&=&
V_2^*J_{\tilde{\psi_2}}V_2(J_{\psi_2}\underset{N^o}{_{a_2}\otimes_\beta}J_{\widehat{\Phi}})\\
&=&(u\underset{N}{_b\otimes_\alpha}1)V_1^*J_{\tilde{\psi_1}}\tilde{J_{2,1}}^*J_{\tilde{\psi_2}}V_2(J_{\psi_2}\underset{N^o}{_{a_2}\otimes_\beta}J_{\widehat{\Phi}})\\
&=&(u\underset{N}{_b\otimes_\alpha}1)V_1^*\tilde{J_{2,1}}^*V_2(J_{\psi_2}\underset{N^o}{_{a_2}\otimes_\beta}J_{\widehat{\Phi}})\\
&=&(u\underset{N}{_b\otimes_\alpha}1)U_{\psi_1}^\mathfrak a(J_{\psi_1}\underset{N}{_b\otimes_\alpha}J_{\widehat{\Phi}})V_1^*J_{\tilde{\psi_1}}\tilde{J_{2,1}}^*V_2(J_{\psi_2}\underset{N^o}{_{a_2}\otimes_\beta}J_{\widehat{\Phi}})\\
&=&(u\underset{N}{_b\otimes_\alpha}1)U_{\psi_1}^\mathfrak a(J_{\psi_1}\underset{N}{_b\otimes_\alpha}J_{\widehat{\Phi}})V_1^*V_1(u^*\underset{N}{_b\otimes_\alpha}1)V_2^*V_2(J_{\psi_2}\underset{N^o}{_{a_2}\otimes_\beta}J_{\widehat{\Phi}})\\
&=&(u\underset{N}{_b\otimes_\alpha}1)U_{\psi_1}^\mathfrak a(J_{\psi_1}\underset{N}{_b\otimes_\alpha}J_{\widehat{\Phi}})(u^*\underset{N}{_b\otimes_\alpha}1)(J_{\psi_2}\underset{N^o}{_{a_2}\otimes_\beta}J_{\widehat{\Phi}})\\
&=&(u\underset{N}{_b\otimes_\alpha}1)U_{\psi_1}^\mathfrak a(u^*\underset{N^o}{_{a_2}\otimes_\beta}1)
\end{eqnarray*}
from which we finish the proof of (i). Using the intertwining properties of $u$, (i) and (\cite{E5} 5.2), we then get (ii). Using then (ii) and the properties of $U^\ga_\psi$ (\cite{E5} 13.4) recalled in \ref{crossed}, we get (iii). \end{proof}

%%%%Uinv
\subsection{Proposition}
\label{Uinv}
{\it Let $\gG$ be a measured quantum groupoid, and $(b,\mathfrak a)$ an action of $\gG$ on a von Neumann algebra $A$; let $\psi$ be a $\delta$-invariant weight on $A$, bearing the density condition, as defined in \ref{action}; then :
\newline
(i) the unitary $U^\mathfrak a_\psi$ constructed in \ref{crossed} is equal to the implementation $V_\psi$ of $\mathfrak a$ constructed in \ref{action}. 
\newline
(ii) the dual weight satisfies $\Delta_{\tilde{\psi}}^{it}=\Delta_\psi^{it}\underset{N}{_b\otimes_\alpha}(\delta\Delta_{\widehat{\Phi}})^{-it}$, where this last one-parameter group of unitaries had been defined in \ref{action}. 

\begin{proof}
Let $\xi\in D(_\alpha H_\Phi, \nu)$, $x$, $y$ in $\gN_\psi\cap\gN_\psi^*$, $a\in\gN_{\hat{T}^c}\cap\gN_{\hat{T}^c}^*\cap\gN_{\widehat{\Phi}^c}\cap\gN_{\widehat{\Phi}^c}^*$, such that $\Lambda_{\widehat{\Phi}^c}(a^*)$ belongs to the set $\widehat{\mathcal E_{\hat{\tau}}}$ introduced in (\cite{E5}4.4). We have, using \ref{crossed}:
\[(\rho_\xi^{b, \alpha})^*S_{\tilde{\psi}}\mathfrak a(x^*)(\Lambda_\psi(y)\underset{\nu}{_b\otimes_\alpha}\Lambda_{\widehat{\Phi}^c}(a))
=
(\rho_\xi^{b, \alpha})^*\mathfrak a(y^*)(\Lambda_\psi(x)\underset{\nu}{_b\otimes_\alpha}\Lambda_{\widehat{\Phi}^c}(a^*))\]
and, as $\Lambda_{\widehat{\Phi}^c}(a^*)$ belongs to $D(_\alpha H_\Phi, \nu)$, thanks to (\cite{E5}4.4) it is equal to :
\[(id\underset{N}{_b*_\alpha}\omega_{\Lambda_{\widehat{\Phi}^c}(a^*), \xi})\mathfrak a(y^*)\Lambda_\psi(x)
=\Lambda_\psi((id\underset{N}{_b*_\alpha}\omega_{\Lambda_{\widehat{\Phi}^c}(a^*), \xi})\mathfrak a(y^*)x)\]
Let us suppose now that $x$ is analytic with respect to $\psi$; as $\delta^{1/2}\Lambda_{\widehat{\Phi}^c}(a^*)$ belongs to $D((H_\Phi)_\beta, \nu^o)$, thanks again to (\cite{E5}4.4), we get, using (\cite{E5} 8.4.(iii)), that it is equal to :
\begin{multline*}
J_\Psi\sigma_{-i/2}^\psi(x^*)J_\psi\Lambda_\psi[(id\underset{N}{_b*_\alpha}\omega_{\Lambda_{\widehat{\Phi}^c}(a^*), \xi})\mathfrak a(y^*)]=\\
=
J_\Psi\sigma_{-i/2}^\psi(x^*)J_\psi(id*\omega_{\delta^{1/2}\Lambda_{\widehat{\Phi}^c}(a^*), \xi})(V_\psi)\Lambda_\psi(y^*)\\
=
(\rho_\xi^{b, \alpha})^* (J_\Psi\sigma_{-i/2}^\psi(x^*)J_\psi\underset{\nu}{_b\otimes_\alpha}1)V_\psi(\Lambda_\psi(y^*)\underset{\nu^o}{_a\otimes_\beta}\delta^{1/2}\Lambda_{\widehat{\Phi}^c}(a^*))
\end{multline*}
from which we get that :
\[S_{\tilde{\psi}}\mathfrak a(x^*)(\Lambda_\psi(y)\underset{\nu}{_b\otimes_\alpha}\Lambda_{\widehat{\Phi}^c}(a))=
(J_\Psi\sigma_{-i/2}^\psi(x^*)J_\psi\underset{\nu}{_b\otimes_\alpha}1)V_\psi(\Lambda_\psi(y^*)\underset{\nu^o}{_a\otimes_\beta}\delta^{1/2}\Lambda_{\widehat{\Phi}^c}(a^*))\]
and, taking a bounded net $x_i$ strongly converging to $1$, such that $\sigma_{-i/2}^\psi(x_i^*)$ is also converging to $1$, and using the fact that $S_{\tilde{\psi}}$ is closed, we get :
\[S_{\tilde{\psi}}(\Lambda_\psi(y)\underset{\nu}{_b\otimes_\alpha}\Lambda_{\widehat{\Phi}^c}(a))=
V_\psi[J_\psi\Delta_\psi^{1/2}\Lambda_\psi(y)\underset{\nu^o}{_a\otimes_\beta}J_{\widehat{\Phi}}(\overline{\delta\Delta_{\widehat{\Phi}}})^{-1/2}\Lambda_{\widehat{\Phi}^c}(a)]\]
from which we deduce that :
\[V_\psi(J_\psi\underset{N}{_b\otimes_\alpha}J_{\widehat{\Phi}})(\Delta_\psi^{1/2}\underset{N}{_b\otimes_\alpha}\overline{\delta\Delta_{\widehat{\Phi}}}^{-1/2})\subset S_{\tilde{\psi}}\]
where $\Delta_\psi^{1/2}\underset{N}{_b\otimes_\alpha}(\overline{\delta\Delta_{\widehat{\Phi}}})^{-1/2}$ is the infinitesimal generator of the one-parameter group of unitaries $\Delta^{it}_\psi\underset{N}{_b\otimes_\alpha}\delta^{-it}\Delta_{\widehat{\Phi}}^{-it}$ introduced in \ref{action}.  
But, on the other hand, for all $t\in\mathbb{R}$, we have, using \ref{crossed} :
\begin{multline*}
(\Delta^{it}_\psi\underset{N}{_b\otimes_\alpha}\delta^{-it}\Delta_{\widehat{\Phi}}^{-it})S_{\tilde{\psi}}\mathfrak a(x^*)(\Lambda_\psi(y)\underset{\nu}{_b\otimes_\alpha}\Lambda_{\widehat{\Phi}^c}(a))=\\
=(\Delta^{it}_\psi\underset{N}{_b\otimes_\alpha}\delta^{-it}\Delta_{\widehat{\Phi}}^{-it})\mathfrak a(y^*)(\Lambda_\psi(x)\underset{\nu}{_b\otimes_\alpha}\Lambda_{\widehat{\Phi}^c}(a^*))
\end{multline*}
which, using \ref{action}, is equal to :
\begin{multline*}
\mathfrak a(\sigma^\psi_t(y^*))(\Lambda_\psi(\sigma_t^\psi(x))\underset{\nu}{_b\otimes_\alpha}S_{\widehat{\Phi}^c}\delta^{-it}\Delta_{\widehat{\Phi}}^{-it}\Lambda_{\widehat{\Phi}^c}(a))=\\
\mathfrak a(\sigma^\psi_t(y^*))(\Lambda_\psi(\sigma_t^\psi(x))\underset{\nu}{_b\otimes_\alpha}S_{\widehat{\Phi}^c}\delta^{-it}\Delta_{\widehat{\Phi}}^{-it}\Lambda_{\widehat{\Phi}^c}(a))
\end{multline*}
which is equal, using again \ref{crossed}, to :
\[S_{\tilde{\psi}}\mathfrak a(\sigma^\psi_t(x^*))(\Lambda_\psi(\sigma_t^\psi(y))\underset{\nu}{_b\otimes_\alpha}\delta^{-it}\Delta_{\widehat{\Phi}}^{-it}\Lambda_{\widehat{\Phi}^c}(a))\]
Taking again a family $x_i$ converging to $1$, and using the closedness of $S_{\tilde{\psi}}$, we get that :
\begin{multline*}
(\Delta^{it}_\psi\underset{N}{_b\otimes_\alpha}\delta^{-it}\Delta_{\widehat{\Phi}}^{-it})S_{\tilde{\psi}}(\Lambda_\psi(y)\underset{\nu}{_b\otimes_\alpha}\Lambda_{\widehat{\Phi}^c}(a))=\\
S_{\tilde{\psi}}(\Lambda_\psi(\sigma_t^\psi(y))\underset{\nu}{_b\otimes_\alpha}\delta^{-it}\Delta_{\widehat{\Phi}}^{-it}\Lambda_{\widehat{\Phi}^c}(a))=\\
S_{\tilde{\psi}}(\Delta^{it}_\psi\underset{N}{_b\otimes_\alpha}\delta^{-it}\Delta_{\widehat{\Phi}}^{-it})
(\Lambda_\psi(y)\underset{\nu}{_b\otimes_\alpha}\Lambda_{\widehat{\Phi}^c}(a))
\end{multline*}
from which, using \ref{crossed}, we deduce that 
\[(\Delta^{it}_\psi\underset{N}{_b\otimes_\alpha}\delta^{-it}\Delta_{\widehat{\Phi}}^{-it})S_{\tilde{\psi}}=S_{\tilde{\psi}}(\Delta^{it}_\psi\underset{N}{_b\otimes_\alpha}\delta^{-it}\Delta_{\widehat{\Phi}}^{-it})\]
and, therefore, we have :
\[V_\psi(J_\psi\underset{N}{_b\otimes_\alpha}J_{\widehat{\Phi}})(\Delta_\psi^{1/2}\underset{N}{_b\otimes_\alpha}\overline{\delta\Delta_{\widehat{\Phi}}}^{-1/2})= S_{\tilde{\psi}}\]
and, by polar decomposition, we have :
\[J_{\tilde{\psi}}=V_\psi(J_\psi\underset{N}{_b\otimes_\alpha}J_{\widehat{\Phi}})\]
which, by definition of $U_\psi^\ga$, leads to (i).
\newline
We also get :
\[\Delta_{\tilde{\psi}}^{1/2}=\Delta_\psi^{1/2}\underset{N}{_b\otimes_\alpha}\overline{\delta\Delta_{\widehat{\Phi}}}^{-1/2}\]
which leads to (ii). \end{proof}

%%%%corstandard
\subsection{Corollary}
\label{corstandard}
{\it Let $\gG$ be a measured quantum groupoid, and $(b,\mathfrak a)$ an action of $\gG$ on a von Neumann algebra $A$; let us suppose that there exists on $A$ a $\delta$-invariant weight on $A$, bearing the density condition, as defined in \ref{action}; then, for any normal semi-finite faithful weight $\psi$ on $A$, the unitary $U^\mathfrak a_\psi$ constructed in \ref{crossed} is a standard implementation of $\mathfrak a$ as defined in \ref{crossed}. }

\begin{proof}
If $\psi$ is a $\delta$-invariant weight on $A$, bearing the density condition, as defined in \ref{action}, we have the result using \ref{Uinv}; for another weight, using \ref{propu}(iii), we get the result. \end{proof}

%%%%%cordual
\subsection{Corollary}
\label{cordual}
{\it Let $\gG$ be a measured quantum groupoid, and $(b,\mathfrak a)$ an action of $\gG$ on a von Neumann algebra $A$; let us suppose that $A$ is isomorphic to a crossed-product $B\rtimes_\gb\widehat{\gG}^o$ where $\gb$ is an action of $\widehat{\gG}^o$ on a von Neumann algebra $B$, and that this isomorphism sends $\ga$ on $\tilde{\gb}$. Then, for any normal semi-finite faithful weight $\psi$ on $A$, the unitary $U^\mathfrak a_\psi$ constructed in \ref{crossed} is a standard implementation of $\mathfrak a$ as defined in \ref{crossed}. }

\begin{proof}
We have recalled in \ref{crossed} that any dual weight on $B\rtimes_\gb\widehat{\gG}^o$ is a $\delta$-invariant weight on $B\rtimes_\gb\widehat{\gG}^o$, bearing the density condition; therefore, using \ref{corstandard}, we get the result. 
\end{proof}

%%%%cora
\subsection{Corollary}
\label{cora}
{\it Let $\gG$ be a measured quantum groupoid, and $(b,\mathfrak a)$ an action of $\gG$ on a von Neumann algebra $A$; let us consider the action $(1\underset{N}{_b\otimes_\alpha}\hat{\beta}, \underline{\ga})$ of $\gG$ on $A\underset{N}{_b*_\alpha}\mathcal L(H)$, introduced in \ref{crossed}; then, for any normal semi-finite faithful weight $\psi$ on $A\underset{N}{_b*_\alpha}\mathcal L(H)$, the unitary $U^{\underline{\ga}}_\psi$ is a standard implementation of the action $\underline{\ga}$}. }
\begin{proof}
This is just a corollary of \ref{cordual} and of the biduality theorem, recalled in \ref{crossed}. \end{proof}

%%%%corsigma
\subsection{Corollary}
\label{corsigma}
{\it Let $\gG$ be a measured quantum groupoid, and $(b, \ga)$ an action of $\gG$ on a von Neumann algebra $A$; let $\psi$ be a $\delta$-invariant weight on $A$, bearing the density condition, as defined in \ref{action}; then, for any $x\in\widehat{M}'$, $t\in\mathbb{R}$, we have :}
\[\sigma_t^{\tilde{\psi}}(1\underset{N}{_b\otimes_\alpha}x)=1\underset{N}{_b\otimes_\alpha}\Delta_{\Phi}^{it}x\Delta_{\Phi}^{-it}\] 

\begin{proof}
Using \ref{Uinv}(ii), we get that :
\[\sigma_t^{\tilde{\psi}}(1\underset{N}{_b\otimes_\alpha}x)=1\underset{N}{_b\otimes_\alpha}(\delta\Delta_{\widehat{\Phi}})^{-it}x(\delta\Delta_{\widehat{\Phi}})^{it}\]
But, using (\cite{E5}3.11(ii)), we know that $(\delta\Delta_{\widehat{\Phi}})^{it}=(\hat{\delta}\Delta_\Phi)^{-it}$; as $\hat{\delta}$ is affiliated to $\widehat{M}$, we get the result. 
\end{proof}

%%%%%corsigma2
\subsection{Corollary}
\label{corsigma2}
{\it Let $\gG$ be a measured quantum groupoid, and $(b, \ga)$ an action of $\gG$ on a von Neumann algebra $A$; let $\psi$ be a normal semi-finite faithful weight on $A$; then, for any $x$ in $M'$, $t\in\mathbb{R}$, we have :}
\[\sigma_t^{\tilde{\tilde{\psi}}}(1\underset{N}{_b\otimes_\alpha}1\underset{N^o}{_{\hat{\alpha}}\otimes_\beta}x)=1\underset{N}{_b\otimes_\alpha}1\underset{N^o}{_{\hat{\alpha}}\otimes_\beta}\Delta_{\widehat{\Phi}}^{-it}x\Delta_{\widehat{\Phi}}^{it}\]
\begin{proof}
Let's apply \ref{corsigma} to the dual action $(1\underset{N}{_b\otimes_\alpha}\hat{\alpha}, \tilde{\ga})$ of $\gG^c$ on $A\rtimes_\ga\gG$, and the dual weight $\tilde{\psi}$, and we get the result. \end{proof}

%%%%%corsigma3
\subsection{Corollary}
\label{corsigma3}
{\it Let $\gG$ be a measured quantum groupoid, and $(b, \ga)$ an action of $\gG$ on a von Neumann algebra $A$; let $\psi$ be a normal semi-finite faithful weight on $A$; let $(1\underset{N}{_b\otimes_\alpha}\hat{\beta}, \underline{\ga})$ be the action of $\gG$ on $A\underset{N}{_b*_\alpha}\mathcal L(H)$ obtained by transporting on $A\underset{N}{_b*_\alpha}\mathcal L(H)$ the bidual action and $\overline{\psi}_\ga$ be the normal semi-finite faithful weight on $A\underset{N}{_b*_\alpha}\mathcal L(H)$ obtained by transporting the bidual weight. Then,  for any $x$ in $M'$, $t\in\mathbb{R}$, we have :
\[\sigma_t^{\overline{\psi}_\ga}(1\underset{N}{_b\otimes_\alpha}x)=1\underset{N}{_b\otimes_\alpha}\Delta_{\widehat{\Phi}}^{-it}x\Delta_{\widehat{\Phi}}^{it}\]}

\begin{proof} The canonical isomorphism between $A\underset{N}{_b*_\alpha}\mathcal L(H)$ and $(A\rtimes_\ga\gG)\rtimes_{\tilde{\ga}}\widehat{\gG}^c$ sends, for all $x\in M'$, $1\underset{N}{_b\otimes_\alpha}x$ on $1\underset{N}{_b\otimes_\alpha}1\underset{N^o}{_{\hat{\alpha}}\otimes_\beta}x$ (cf. \cite{E5} 11.2). So, the result is a straightforward consequence of \ref{corsigma2}. \end{proof}

%%%%%lift
\section{An auxilliary weight $\underline{\psi}$.}
\label{auxilliary}
If $b$ is a normal faithful non degenerate anti-homomorphism from $N$ into a  von Neumann algebra $A$, such that there exists a normal faithful semi-finite operator-valued weight $\gT$ from $A$ on $b(N)$, we associate to the weight $\psi=\nu^o\circ b^{-1}\circ\gT$ a weight $\underline{\psi}$ on $A\underset{N}{_b*_\alpha}\mathcal L (H)$ (\ref{psibarre}); we calculate its modular automorphism group (\ref{psibarre2}), and the GNS representation of $A\underset{N}{_b*_\alpha}\mathcal L (H)$ given by this weight (\ref{psibarre3}).

%%%defw
\subsection{Definitions}
\label{defw}
Let $b$ be an injective $*$-antihomorphism from a von Neumann algebra $N$ into a von Neumann algebra $A$; we shall then say that $(N, b, A)$ (or simply $A$) is a  faithful right von Neumann $N$-module. If there exists a normal semi-finite faithful operator-valued weight $\gT$ from $A$ onto $b(N)$, we shall say that this faithful right $N$-module is weighted. 
\newline
Let then $\psi$ be a normal faithful semi-finite weight on $A$; if, for all $t$ in $\mathbb{R}$, $n$ in $N$, we have $\sigma_t^\psi(b(n))=b(\sigma_{-t}^\nu(n))$, then there exists a normal semi-finite faithful operator-valued weight $\gT$ from $A$ onto $b(N)$ such that $\psi=\nu^o\circ b^{-1}\circ \gT$; such a weight $\psi$ on $A$ will be said lifted from $\nu^o$ by $\gT$ (or, simply, a lifted weight).
\newline
If $\psi$ is a normal semi-finite faithful weight on $A$, lifted from $\nu^o$ by $\gT$, then the weight $\psi$ bears the density property introduced in (\cite{E5}, 8.1), recalled in \ref{action}. Namely, using (\cite{E5}, 2.2.1), one gets that $D(_aH_\psi, \nu)\cap D((H_\psi)_b, \nu^o)$ contains all the vectors of the form $\Lambda_\psi(x)$, where $x\in \gN_{\gT}\cap\gN_{\gT}^*\cap\gN_\psi\cap\gN_\psi^*$ is analytical with respect to $\psi$, and such that, for any $z\in\mathbb{C}$, $\sigma_z(x)$ belongs to $\gN_{\gT}\cap\gN_{\gT}^*\cap\gN_\psi\cap\gN_\psi^*$;  therefore $D(_aH_\psi, \nu)\cap D((H_\psi)_b, \nu^o)$ is dense in $H_\psi$, which is the density property. 
\newline
If $(b, \ga)$ is an action of a measured quantum goupoid $\gG=(N, M, \alpha, \beta, \Gamma, T, T', \nu)$ on a von Neumann algebra $A$, we shall say that this action is weighted if the faithful right $N$-module $(N, b, A)$ is weighted.

%%%%%propw
\subsection{Lemma}
\label{propw}
{\it Let $(N, b, A)$ be a faithful weighted right von Neumann $N$-module, and let $\gT$ be a normal semi-finite faithful operator-valued weight from $A$ onto $b(N)$; let $\alpha$ be a nomal faithful representation of $N$ on a Hilbert space $H$ and $\nu$ a normal semi-finite faithful weight on $N$; then, it is possible to define a canonical normal semi-finite faithful operator-valued weight $(\gT\underset{N}{_b*_\alpha}id)$ from $A\underset{N}{_b*_\alpha}\mathcal L(H)$ onto $1\underset{N}{_b\otimes_\alpha}\alpha(N)'$ (which is equal to $b(N)\underset{N}{_b*_\alpha}\mathcal L(H)$, by (\cite{E5}, 2.4)), such that, if $\psi$ denotes the weight on $A$ lifted from $\nu^o$ by $\gT$, we get, for any $X\in (A\underset{N}{_b*_\alpha}\mathcal L(H))^+$, that $(\gT\underset{N}{_b*_\alpha}id)(X)=1\underset{N}{_b\otimes_\alpha}(\psi\underset{\nu}{_b*_\alpha}id)(X)$, where $\gT\underset{N}{_b*_\alpha}id$ and $\psi\underset{\nu}{_b*_\alpha}id$ are slice maps introduced in \cite{E4} and recalled in (\cite{E5}, 2.5).} 

\begin{proof}
Let us represent $A$ on a Hilbert space $\mathcal H$; using Haagerup's theorem (\cite{T}, 4.24), there exists a canonical normal semi-finite faithful operator valued weight $\gT^{-1}$ from $b(N)'$ onto $A'$; considering the representation of $b(N)'$ on $\mathcal H\underset{\nu}{_b\otimes_\alpha}H$, and using again Haagerup's theorem, we obtain another normal semi-finite faithful operator-valued weight $(\gT^{-1})^{-1}$ from the commutant of $A'$ on $\mathcal H\underset{\nu}{_b\otimes_\alpha}H$ (which is $A\underset{N}{_b*_\alpha}\mathcal L(H)$) onto the commutant of $b(N)'$ on $\mathcal H\underset{\nu}{_b\otimes_\alpha}H$ (which is $b(N)\underset{N}{_b*_\alpha}\mathcal L(H)$). As both $\gT$ and $(\gT^{-1})^{-1}$ are obtained by taking the commutants, within two different representations, of the same operator-valued weight $\gT^{-1}$, a closer look at this construction leads (\cite{EN}, 10.2) to the fact that $(\gT^{-1})^{-1}=(\gT\underset{N}{_b*_\alpha}id)$. The link between $(\gT\underset{N}{_b*_\alpha}id)$ and $(\psi\underset{\nu}{_b*_\alpha}id)$ is recalled in (\cite{E5} 2.5).
\end{proof}

%%%subalgebra
\subsection{ Proposition}
\label{subalgebra}
{\it Let $(N, b, A)$ be a von Neumann faithful right $N$-module, and let $\alpha$ be a normal faithful non degenerate representation of $N$ on a Hilbert space $\mathcal H$ and $\nu$ a normal semi-finite faithful weight on $N$; then :
\newline
(i) let's represent $A$ on a Hilbert space $\mathcal K$; the linear set generated by all operators on $\mathcal K\underset{\nu}{_\beta\otimes_\alpha}\mathcal H$, of the form $\rho^{\beta, \alpha}_{\xi_1}a(\rho^{\beta, \alpha}_{\xi_2})^*$, with $a$ in $A$ and $\xi_1$, $\xi_2$ in $D(_\alpha \mathcal H, \nu)$,  is a $*$-algebra, which is weakly dense in $A\underset{N}{_b*_\alpha}\mathcal L(\mathcal{H})$. 
\newline
(ii) let $\psi$ be a normal faithful semi-finite weight on $A$, and let's represent $A$ on $H_\psi$; then, for any $a$ in $\gN_\psi$ and $\xi$ in $D(_\alpha \mathcal H, \nu)$, $\Lambda_\psi(a)\underset{\nu}{_b\otimes_\alpha}\xi$ belongs to $D(H_\psi\underset{\nu}{_b\otimes_\alpha}\mathcal H, \psi^o)$ (where we deal with the representation $x\mapsto x\underset{N}{_b\otimes_\alpha}1$ of $A^o=J_\psi AJ_\psi$), and we have $\theta^{\psi^o}( \Lambda_\psi(a)\underset{\nu}{_b\otimes_\alpha}\xi, \Lambda_\psi(a)\underset{\nu}{_b\otimes_\alpha}\xi)=\rho^{b, \alpha}_\xi a a^*(\rho^{b, \alpha}_\xi)^*$. 
\newline
(iii) for all $n\in N$, let us define $a(n)=J_\psi b(n^*)J_\psi$; let $\gG=(N, M, \alpha, \beta, \Gamma, T, T', \nu)$ be a measured quantum groupoid; then, the representation of $A\underset{N}{_b*_\alpha}\mathcal L(H)$ on $H\underset{\nu}{_\beta\otimes_a}H_\psi\underset{\nu}{_b\otimes_\alpha}H$ defined by $x\mapsto 1\underset{N}{_\beta\otimes_a}x$ is standard, when we equip the Hilbert space with the antilinear involutive isometry $J$ defined, for any $\xi$, $\eta$ in $D(_\alpha H, \nu)$, $\zeta$ in $H_\psi$, by :
\[J(J_{\widehat{\Phi}}\eta\underset{\nu}{_\beta\otimes_a}\zeta\underset{\nu}{_b\otimes_\alpha}\xi)=J_{\widehat{\Phi}}\xi\underset{\nu}{_\beta\otimes_a}J_\psi\zeta\underset{\nu}{_b\otimes_\alpha}\eta\]
and with the closed cone $\mathcal P$ generated by all elements of the form $J_{\widehat{\Phi}}\xi\underset{\nu}{_\beta\otimes_a}\zeta'\underset{\nu}{_b\otimes_\alpha}\xi$, when $\xi$ is in $D(_\alpha H, \nu)$, and $\zeta'$ in the cone $\mathcal P_\psi$ given by the Tomita-Takesaki theory associated to the weight $\psi$. 
\newline
(iv) let $\varphi$ be a normal semi-finite faithful weight on $A\underset{N}{_b*_\alpha}\mathcal L(H)$; then $\Lambda_\psi (a)\underset{\nu}{_b\otimes_\alpha}\xi$ belongs to $\mathcal D((\frac{d\varphi}{d\psi^o})^{1/2})$ if and only if $\rho^{b, \alpha}_\xi a a^*(\rho^{b, \alpha}_\xi)^*$ belongs to $\gM_{\varphi}^+$, and then :
\[\varphi(\rho^{b, \alpha}_\xi a a^*(\rho^{b, \alpha}_\xi)^*)=\|(\frac{d\varphi}{d\psi^o})^{1/2}(\Lambda_\psi(a)\underset{\nu}{_b\otimes_\alpha}\xi)\|^2\]
Moreover, if $\Lambda_\psi (a)\underset{N}{_b\otimes_\alpha}\xi$ belongs to $\mathcal D((\frac{d\varphi}{d\psi^o})^{1/2})$, the vector $(\frac{d\varphi}{d\psi^o})^{1/2}(\Lambda_\psi (a)\underset{N}{_b\otimes_\alpha}\xi)$ belongs to $D(H_\psi\underset{\nu}{_b\otimes_\alpha}H, \varphi)$, and the canonical isomorphism between $H_\varphi$ and $H\underset{\nu}{_\beta\otimes_a}H_\psi\underset{\nu}{_b\otimes_\alpha}H$ sends $R^{\varphi}((\frac{d\varphi}{d\psi^o})^{1/2}(\Lambda_\psi (a)\underset{\nu}{_b\otimes_\alpha}\xi))^*(\zeta\underset{\nu}{_b\otimes_\alpha}\eta)$ on $J_{\widehat{\Phi}}\xi\underset{\nu}{_\beta\otimes_a}J_\psi aJ_\psi\zeta\underset{\nu}{_b\otimes_\alpha}\eta$. }

\begin{proof}
Using \ref{spatial}, we get, for $a_1$, $a_2$ in $A$, and $\xi_1$, $\xi_2$, $\xi_3$, $\xi_4$ in $D(_\alpha \mathcal H, \nu)$, that  :
\[\rho^{\beta, \alpha}_{\xi_1}a_1(\rho^{\beta, \alpha}_{\xi_2})^*\rho^{\beta, \alpha}_{\xi_3}a_2(\rho^{\beta, \alpha}_{\xi_4})^*=\rho^{\beta, \alpha}_{\xi_1}a_1b(<\xi_3, \xi_2>_{\alpha, \nu})a_2(\rho^{\beta, \alpha}_{\xi_4})^*\]
from which we see that this linear set is indeed an algebra; moreover, it is clear that it is invariant under taking the adjoint. Let's take $c\in A'$; we have :
\begin{eqnarray*}
\rho^{\beta, \alpha}_{\xi_1}a(\rho^{\beta, \gamma}_{\xi_2})^*(c\underset{N}{_\beta\otimes_\alpha}1)
&=&\rho^{\beta, \alpha}_{\xi_1}ac(\rho^{\beta, \alpha}_{\xi_2})^*\\
&=&\rho^{\beta, \alpha}_{\xi_1}ca(\rho^{\beta, \alpha}_{\xi_2})^*\\
&=&(c\underset{N}{_\beta\otimes_\alpha}1)
\rho^{\beta, \alpha}_{\xi_1}a(\rho^{\beta, \alpha}_{\xi_2})^*
\end{eqnarray*}
from which we get that $\rho^{\beta, \alpha}_{\xi_1}a(\rho^{\beta, \alpha}_{\xi_2})^*$ belongs to $A\underset{N}{_\beta*_\alpha}\mathcal L(\mathcal{H})$. Let now $X\in A\underset{N}{_\beta*_\alpha}\mathcal L(\mathcal{H})$, and let $(e_i)_{i\in I}$ be a $(\alpha, \nu)$-orthogonal basis of $\mathcal H$; we get that $(id\underset{\nu}{_\beta*_\alpha}\omega_{e_i, e_j})(X)$ belongs to $A$, and we have, when we take the weak limits over the finite subsets $J$, $J'$ of $I$ :
\begin{eqnarray*}
X&=&lim_{J, J'}\sum_{i\in J, j\in J'}(1\underset{N}{_\beta\otimes_\alpha}\theta^{\alpha, \nu}(e_i, e_i))X(1\underset{N}{_\beta\otimes_\alpha}\theta^{\alpha, \nu}(e_j, e_j))\\
&=&lim_{J, J'}\sum_{i\in J, j\in J'}\rho^{\beta, \alpha}_{e_i}(id\underset{\nu}{_\beta*_\alpha}\omega_{e_i, e_j})(X)(\rho^{\beta, \alpha}_{e_j})^*
\end{eqnarray*}
which proves (i).   
\newline
Let $a\in \gN_\psi$, $\xi\in D(_\alpha \mathcal H, \nu)$; then, for all $x\in\gN_\psi$, we have :
\begin{eqnarray*}
J_\psi xJ_\psi\Lambda_\psi (a)\underset{\nu}{_b\otimes_\alpha}\xi
&=&
aJ_\psi\Lambda_\psi (x)\underset{\nu}{_b\otimes_\alpha}\xi\\
&=&\rho^{b, \alpha}_\xi aJ_\psi\Lambda_\psi(x)
\end{eqnarray*}
Therefore, $\Lambda_\psi (a)\underset{\nu}{_b\otimes_\alpha}\xi$ belongs to $D(H_\psi\underset{\nu}{_b\otimes_\alpha}\mathcal H, \psi^o)$, and $R^{\psi^o}(\Lambda_\psi (a)\underset{\nu}{_b\otimes_\alpha}\xi)=\rho^{b, \alpha}_\xi a$. So, we get that $\theta^{\psi^o}(\Lambda_\psi (a)\underset{\nu}{_b\otimes_\alpha}\xi, \Lambda_\psi (a)\underset{\nu}{_b\otimes_\alpha}\xi)=\rho_\xi^{b, \alpha}aa^*(\rho_\xi^{b, \alpha})^*$, which is (ii). 
\newline
By \cite{S}(3.1), we know that $A\underset{N}{_b*_\alpha}\mathcal L(H)$ has a standard representation $x\mapsto x\otimes_\psi 1$ on the Hilbert space $(H_\psi\underset{\nu}{_b\otimes_\alpha}H)\otimes_\psi\overline{(H_\psi\underset{\nu}{_b\otimes_\alpha}H)}$. Using then (ii) and \cite{S}(0.3.1), we get that this Hilbert space is isomorphic to $H\underset{\nu}{_\beta\otimes_a}H_\psi\underset{\nu}{_b\otimes_\alpha}H$, and that this isomorphism sends, for $b\in\gN_\psi$, $\eta\in D(_\alpha \mathcal H, \nu)$ :
\newline
a) the vector $(\Lambda_\psi(a)\underset{\nu}{_b\otimes_\alpha}\xi)\otimes_\psi\overline{(\Lambda_\psi(b)\underset{\nu}{_b\otimes_\alpha}\eta)}$ on $J_{\widehat{\Phi}}\eta\underset{\nu}{_\beta\otimes_a}J_\psi bJ_\psi \Lambda_\psi (a)\underset{\nu}{_b\otimes_\alpha}\xi$, 
 \newline
b) the standard representation $x\mapsto x\otimes_\psi 1$ on the representation $x\mapsto 1\underset{N}{_\beta\otimes_a}x$, 
 \newline
c) the antilinear involutive isometry which sends $(\Lambda_\psi(a)\underset{\nu}{_b\otimes_\alpha}\xi)\otimes_\psi\overline{(\Lambda_\psi(b)\underset{\nu}{_b\otimes_\alpha}\eta)}$ to $(\Lambda_\psi(b)\underset{\nu}{_b\otimes_\alpha}\eta)\otimes_\psi\overline{(\Lambda_\psi(a)\underset{\nu}{_b\otimes_\alpha}\xi)}$ on $J$, 
 \newline
d) the cone generated by all elements of the form $(\Lambda_\psi(a)\underset{\nu}{_b\otimes_\alpha}\xi)\otimes_\psi\overline{(\Lambda_\psi(a)\underset{\nu}{_b\otimes_\alpha}\xi)}$ on $\mathcal P$, which gives (iii). 
\newline
 Using (ii), we get that :
\[\varphi(\rho^{b, \alpha}_\xi a a^*(\rho^{b, \alpha}_\xi)^*)=\varphi(\theta^{\psi^o}( \Lambda_\psi(a)\underset{\nu}{_b\otimes_\alpha}\xi, \Lambda_\psi(a)\underset{\nu}{_b\otimes_\alpha}\xi))\]
and, by definition of the spatial derivative, we know that, if $\Lambda_\psi(a)\underset{N}{_b\otimes_\alpha}\xi$ belongs to $\mathcal D((\frac{d\varphi}{d\psi^o})^{1/2})$, 
we have :
\[\varphi(\rho^{b, \alpha}_\xi a a^*(\rho^{b, \alpha}_\xi)^*)=\varphi(\theta^{\psi^o}( \Lambda_\psi(a)\underset{\nu}{_b\otimes_\alpha}\xi, \Lambda_\psi(a)\underset{\nu}{_b\otimes_\alpha}\xi)=\|(\frac{d\varphi}{d\psi^o})^{1/2}(\Lambda_\psi(a)\underset{\nu}{_b\otimes_\alpha}\xi)\|^2\]
and, if $\Lambda_\psi(a)\underset{\nu}{_b\otimes_\alpha}\xi$ does not belong to $\mathcal D((\frac{d\varphi}{d\psi^o})^{1/2})$, we know that $\varphi(\rho^{b, \alpha}_\xi a a^*(\rho^{b, \alpha}_\xi)^*)=+\infty$. So, we have the first part of (iv). Then, the second part of (iv) is given by \cite{S}(3.2) and (iii). 
\end{proof}

%%%%%psibarre
\subsection{Proposition}
\label{psibarre}
{\it Let $\gG=(N, M, \alpha, \beta, \Gamma, T, T', \nu)$ be a measured quantum groupoid; let $(N, b, A)$ be a faithful weighted right von Neumann $N$-module, and let $\psi$ be a normal semi-finite faithful weight on $A$ lifted from $\nu^o$ in the sense of \ref{defw}. Then :
\newline
(i) it possible to define a one-parameter group of unitaries on $\Delta_\psi^{it}\underset{\nu}{_b\otimes_\alpha}\Delta_{\widehat{\Phi}}^{-it}$ on $H_\psi\underset{\nu}{_b\otimes_\alpha}H$, with natural values on elementary tensors. This one-parameter group of unitaries implements on $A\underset{N}{_b*_\alpha}\mathcal L(H)$ the one-parameter group of automorphisms $\sigma_t^\psi\underset{N}{_b*_\alpha}Ad\Delta_{\widehat{\Phi}}^{-it}$. 
\newline
(ii) there exists a normal semi-finite faithful weight $\underline{\psi}$ on $A\underset{N}{_b*_\alpha}\mathcal L(H)$ such that the spatial derivative $\frac{d\underline{\psi}}{d\psi^o}$ is equal to the generator $\Delta_\psi\underset{\nu}{_b\otimes_\alpha}\Delta_{\widehat{\Phi}}^{-1}$ of the one-parameter group of unitaries constructed in (i); the modular automorphism group $\sigma_t^{\underline{\psi}}$ is equal to the automorphism group $\sigma_t^\psi\underset{N}{_b*_\alpha}Ad \Delta_{\widehat{\Phi}}^{-it}$ constructed in (i).
\newline
(iii) for any $a$ in $\gN_\psi\cap\gN_\psi^*$, and $\xi\in D(_\alpha H, \nu)\cap \mathcal D(\Delta_{\widehat{\Phi}}^{-1/2})$, such that $\Delta_{\widehat{\Phi}}^{-1/2}\xi$ belongs to $D(_\alpha H, \nu)$, we have :}
\[\underline{\psi}(\rho_\xi^{b, \alpha}aa^*(\rho_\xi^{b, \alpha})^*)=\|\Delta_\psi^{1/2}\Lambda_\psi (a)\underset{\nu}{_b\otimes_\alpha}\Delta_{\widehat{\Phi}}^{-1/2}\xi\|^2\]

\begin{proof}
If $\xi\in D(_\alpha H, \nu)$, we get, for all $t\in \mathbb{R}$ and $n\in \gN_\nu$ :
\begin{eqnarray*}
\alpha(n)\Delta_{\widehat{\Phi}}^{-it}\xi&=&\Delta_{\widehat{\Phi}}^{-it}\sigma_t^{\widehat{\Phi}}(\alpha(n))\xi\\
&=&\Delta_{\widehat{\Phi}}^{-it}\alpha(\sigma_t^\nu(n))\xi\\
&=&\Delta_{\widehat{\Phi}}^{-it}R^{\alpha, \nu}(\xi)\Delta_\nu^{it}\Lambda_\nu(n)
\end{eqnarray*}
from which we get that $\Delta_{\widehat{\Phi}}^{-it}\xi$ belongs to $D(_\alpha H, \nu)$, and $R^{\alpha, \nu}(\Delta_{\widehat{\Phi}}^{-it}\xi)=\Delta_{\widehat{\Phi}}^{-it}R^{\alpha, \nu}(\xi)\Delta_\nu^{it}$. Therefore, we have $<\Delta_{\widehat{\Phi}}^{-it}\xi, \Delta_{\widehat{\Phi}}^{-it}\xi>_{\alpha, \nu}^o=\sigma_{-t}^\nu(<\xi, \xi>_{\alpha, \nu}^o)$. Taking now $\eta\in H_\psi$, we get :
\begin{eqnarray*}
\|\Delta_\psi^{it}\eta\underset{\nu}{_b\otimes_\alpha}\Delta_{\widehat{\Phi}}^{-it}\xi\|^2
&=&(b(<\Delta_{\widehat{\Phi}}^{-it}\xi, \Delta_{\widehat{\Phi}}^{-it}\xi>_{\alpha, \nu}^o)\Delta_\psi^{it}\eta|\Delta_\psi^{it}\eta)\\
&=&(b(\sigma_{-t}^\nu(<\xi, \xi>_{\alpha, \nu}^o))\Delta_\psi^{it}\eta|\Delta_\psi^{it}\eta)\\
&=&(\sigma_t^\psi(b(<\xi, \xi>_{\alpha, \nu}^o))\Delta_\psi^{it}\eta|\Delta_\psi^{it}\eta)\\
&=&(b(<\xi, \xi>_{\alpha, \nu}^o)\eta|\eta)\\
&=&\|\eta\underset{\nu}{_b\otimes_\alpha}\xi\|^2
\end{eqnarray*}
from which we get the existence of the one-parameter group of unitaries. It is then easy to finish the proof of (i). 
\newline
As $(\Delta_\psi^{it}\underset{\nu}{_b\otimes_\alpha}\Delta_{\widehat{\Phi}}^{-it})(J_\psi xJ_\psi\underset{N}{_b\otimes_\alpha}1)(\Delta_\psi^{-it}\underset{\nu}{_b\otimes_\alpha}\Delta_{\widehat{\Phi}}^{it})=
J_\psi \sigma_t^\psi(x)J_\psi\underset{N}{_b\otimes_\alpha}1$, we obtain (\cite{T}, 3.11) that there exists a normal faithful semi-finite weight $\underline{\psi}$ on $A\underset{N}{_b*_\alpha}\mathcal L(H)$ such that :
\[\frac {d\underline{\psi}}{d\psi^o}=\Delta_\psi\underset{\nu}{_b\otimes_\alpha}\Delta_{\widehat{\Phi}}^{-1}\]
Moreover, the modular automorphism group $\sigma_t^{\underline{\psi}}$ is then equal to the one-parameter automorphism group $\sigma_t^\psi\underset{N}{_b*_\alpha}Ad\Delta_{\widehat{\Phi}}^{-it}$, constructed in (i), which finishes the proof of (ii). 
\newline
 So, using now \ref{subalgebra}(iv) applied to $\underline{\psi}$, we get that $\rho_\xi^{b, \alpha}aa^*(\rho_\xi^{b, \alpha})^*$ belongs to $\gM_{\underline{\psi}}^+$ if and only if $\Lambda_\psi (a)\underset{\nu}{_b\otimes_\alpha}\xi$ belongs to $\mathcal D(\Delta_\psi^{1/2}\underset{\nu}{_b\otimes_\alpha}\Delta_{\widehat{\Phi}}^{-1/2})$, and then, we have :
\[\underline{\psi}(\rho_\xi^{b, \alpha}aa^*(\rho_\xi^{b, \alpha})^*)=\|\Delta_\psi^{1/2}\Lambda_\psi (a)\underset{\nu}{_b\otimes_\alpha}\Delta_{\widehat{\Phi}}^{-1/2}\xi\|^2\]
from which we get (iii). 
\end{proof}

%%%corunderline
\subsection{Corollary}
\label{corunderline}
{\it Let $\gG=(N, M, \alpha, \beta, \Gamma, T, T', \nu)$ be a measured quantum groupoid; let $(N, b, A)$ be a faithful weighted right von Neumann $N$-module, 
$\psi_1$ (resp. $\psi_2$) be a normal faithful semi-finite weight on $A$ lifted from $\nu^o$ and $\underline{\psi_1}$ (resp. $\underline{\psi_2}$) be  the normal semi-finite faithful weight on $A\underset{N}{_b*_\alpha}\mathcal L(H)$ constructed in \ref{psibarre}(ii); then :
\newline
(i) the cocycle $(D\psi_1:D\psi_2)_t$ belongs to $A\cap b(N)'$; 
\newline
(ii) we have : $(D\underline{\psi_1} : D\underline{\psi_2})_t=(D\psi_1:D\psi_2)_t\underset{N}{_b\otimes_\alpha}1$. }

\begin{proof}
As $\psi_1$ and $\psi_2$ are lifted weights, (i) is well known (\cite{T}, 4.22. (iii)). 
\newline
Let $(\gH, \pi, J, \mathcal P)$ be a standard representation of the von Neumann algebra $A$; then $A^o$ is represented on $\gH$ by $JAJ$; for any normal semi-finite faithful weight $\psi$ on $A$, we have $\frac{d\psi}{d\psi^o}=\Delta_\psi^{1/2}$; moreover, we have then :
\begin{align*}
(\frac{d\psi_1}{d\psi_1^o})^{it}(D\psi_1^o:D\psi_2^o)_t(\frac{d\psi_2^o}{d\psi_2})^{it}
&=
(\frac{d\psi_1}{d\psi_1^o})^{it}(\frac{d\psi_1^o}{d\psi_1})^{it}(\frac{d\psi_2^o}{d\psi_1})^{-it}(\frac{d\psi_2^o}{d\psi_2})^{it}\\
&=
(\frac{d\psi_1}{d\psi_2^o})^{it}(\frac{d\psi_2}{d\psi_2^o})^{-it}\\
&=
(D\psi_1: D\psi_2)_{t}
\end{align*}
and, therefore $(D\psi_1^o:D\psi_2^o)_t=\Delta_{\psi_1}^{-it}(D\psi_1: D\psi_2)_{t}\Delta_{\psi_2}^{it}$. 
By similar arguments, we have on $\gH\underset{\nu}{_b\otimes_\alpha}H$ :
\begin{align*}
(D\underline{\psi_1} : D\underline{\psi_2})_t
&=
(\frac{d\underline{\psi_1}}{d\psi_1^o})^{it}(\frac{d\psi_1^o}{d\underline{\psi_2}})^{it}\\
&=(\frac{d\underline{\psi_1}}{d\psi_1^o})^{it}(D\psi_1^o:D\psi_2^o)_t(\frac{d\underline{\psi_2}}{d\psi_2^o})^{-it}
\end{align*}
As $(D\psi_1^o:D\psi_2^o)_t$ belongs to $JAJ\underset{N}{_b\otimes_\alpha}1_H$ and is therefore equal to :
\[\Delta_{\psi_1}^{-it}(D\psi_1:D\psi_2)_{t}\Delta_{\psi_2}^{it}\underset{N}{_b\otimes_\alpha}1_H\]
we obtain, using \ref{psibarre}(ii), that $(D\underline{\psi_1} : D\underline{\psi_2})_t$ is equal to :
\[(\Delta_{\psi_1}^{it}\underset{N}{_b\otimes_\alpha}\Delta_{\widehat{\Phi}}^{-it})(\Delta_{\psi_1}^{-it}(D\psi_1:D\psi_2)_{t}\Delta_{\psi_2}^{it}\underset{N}{_b\otimes_\alpha}1_H)(\Delta_{\psi_2}^{-it}\underset{N}{_b\otimes_\alpha}\Delta_{\widehat{\Phi}}^{it})\]
from which we get the result. \end{proof}

%%%%example
\subsection{Remarks}
\label{example}
Let us consider the trivial action $(id, \beta)$ of $\gG$ on $N^o$ (\cite{E5}, 6.2); it is clearly a weighted action (with the identity of $N^o$ as operator-valued weight); the crossed product is then $\widehat{M}'$, and the dual action is equal to $\widehat{\Gamma}^c$ (\cite{E5}, 9.4); the operator-valued weight from $\widehat{M}'$ onto $\beta(N)$ is then $\hat{T}^c$, and, therefore, the dual weight $\tilde{\nu^o}$ of the weight $\nu^o$ on $N^o$ is the Haar weight $\widehat{\Phi}^c$; by the biduality theorem (\ref{action}), we get that the crossed-product $\widehat{M}'\rtimes_{\widehat{\Gamma}^c}\widehat{\gG}^c$
 is isomorphic to $N^o*\mathcal L(H)=\alpha(N)'$; transporting the bidual weight $\tilde{\tilde{\nu^o}}$ on $\widehat{M}'\rtimes_{\widehat{\Gamma}^c}\widehat{\gG}^c$ by this isomorphism, we obtain the weight $\overline{\nu^o_\beta}$ on $\alpha(N)'$, which verifies, thanks to \ref{crossed}, for all $\xi\in D(_\alpha H, \nu)\cap\mathcal D(\Delta_{\widehat{\Phi}}^{-1/2})$ :
\[\overline{\nu^o_\beta}(\theta^{\alpha, \nu}(\xi, \xi))=\|\Delta_{\widehat{\Phi}}^{-1/2}\xi\|^2\]
and, for all $t\in\mathbb{R}$, $x\in\alpha(N)'$, $\sigma_t^{\overline{\nu^o_\beta}}(x)=\Delta_{\widehat{\Phi}}^{-it}x\Delta_{\widehat{\Phi}}^{it}$. 
\newline
On the other hand, for any $y\in\gN_{\hat{T}^c}\cap\gN_{\widehat{\Phi}^c}$, $z\in \gN_{T^{oc}}\cap\gN_{\Phi^{oc}}$, we have, by construction of $\overline{\nu^o_\beta}$ :
\[\overline{\nu^o_\beta}(y^*z^*zy)=\widehat{\Phi^c}(y^*T^{oc}(z^*z)y)=\|\Lambda_{\widehat{\Phi}^c}(y)\underset{\nu^o}{_{\hat{\alpha}}\otimes_\beta}\Lambda_{\Phi^{oc}}(z)\|^2\]
Let now $(b, \ga)$ be an action of $\gG$ on a von Neumann algebra $A$, and $\psi$ a normal semi-finite faithful weight on $A$; by construction of $\overline{\psi_\ga}$, we have, for any $x\in\gN_\psi$ :
\[\overline{\psi_\ga}(\ga(x^*)(1\underset{N}{_b\otimes_\alpha}y^*z^*zy)\ga(x))=
\|\Lambda_\psi(x)\underset{\nu}{_b\otimes_\alpha}\Lambda_{\widehat{\Phi}^c}(y)\underset{\nu^o}{_{\hat{\alpha}}\otimes_\beta}\Lambda_{\Phi^{oc}}(z)\|^2\]
and, by applying (\cite{E5},13.3) to the weight $\tilde{\tilde{\nu^o}}$, we get, for any $X\in\gN_{\overline{\nu^o_\beta}}$ such that $\Lambda_{\overline{\nu^o_\beta}}(X)$ belongs to $D(_\alpha H_{\overline{\nu^o_\beta}}, \nu)$ :
\[\overline{\psi_\ga}(\ga(x^*)(1\underset{N}{_b\otimes_\alpha}X^*X)\ga(x))=
\|\Lambda_\psi(x)\underset{\nu}{_b\otimes_\alpha}\Lambda_{\overline{\nu^o_\beta}}(X)\|^2\]

%%%%%lem
\subsection{Lemma}
\label{lem}
{\it Let $\gG=(N, M, \alpha, \beta, \Gamma, T, T', \nu)$ be a measured quantum groupoid; let $(N, b, A)$ be a faithful weighted right von Neumann right $N$-module; then :
\newline
(i) if $\xi$, $\eta$ are in $D(_\alpha H, \nu)\cap \mathcal D(\Delta_{\widehat{\Phi}}^{-1/2})$, such that $\Delta_{\widehat{\Phi}}^{-1/2}\xi$ and $\Delta_{\widehat{\Phi}}^{-1/2}\eta$ belong to $D(_\alpha H, \nu)$, $<\Delta_{\widehat{\Phi}}^{-1/2}\xi, \eta>_{\alpha, \nu}^o$ belongs to $\mathcal D(\sigma_{-i/2}^\nu)$ and $\sigma_{-i/2}^\nu(<\Delta_{\widehat{\Phi}}^{-1/2}\xi, \eta>_{\alpha, \nu}^o)=<\xi, \Delta_{\widehat{\Phi}}^{-1/2}\eta>_{\alpha, \nu}^o$. 
\newline
(ii) there exists an $(\alpha, \nu)$-orthogonal basis of $H$ such that, for all $i\in I$, $e_i$ belongs to $D(_\alpha H, \nu)\cap \mathcal D(\Delta_{\widehat{\Phi}}^{-1/2})$ and $\Delta_{\widehat{\Phi}}^{-1/2}e_i$ belongs to $D(_\alpha H, \nu)$; 
\newline
(iii) for any such basis, the weight $\overline{\nu^o_\beta}$ defined in \ref{example} satisfies, for all $x\in \alpha(N)'^+$ :}
\[\overline{\nu^o_\beta} (x)=\sum_i(x\Delta_{\widehat{\Phi}}^{-1/2}e_i|\Delta_{\widehat{\Phi}}^{-1/2}e_i)\]

\begin{proof}
We get, for any $n\in\gN_\nu$, analytic with respect to $\nu$ :
\[R^{\alpha, \nu}(\Delta_{\widehat{\Phi}}^{-1/2}\xi)\Lambda_\nu (n)=\alpha(n)\Delta_{\widehat{\Phi}}^{-1/2}\xi=\Delta_{\widehat{\Phi}}^{-1/2}\alpha(\sigma_{-i/2}^\nu(n))\xi=\Delta_{\widehat{\Phi}}^{-1/2}R^{\alpha, \nu}(\xi)\Delta_\nu^{1/2}\Lambda_\nu (n)\]
and, using (\cite{C2}, 1.5) :
\begin{eqnarray*}
\Lambda_\nu(<\eta, \Delta_{\widehat{\Phi}}^{-1/2}\xi>_{\alpha, \nu}^o)&=&
J_\nu\Delta_\nu^{1/2}\Lambda_\nu(<\Delta_{\widehat{\Phi}}^{-1/2}\xi, \eta>_{\alpha, \nu}^o)\\
&=&
J_\nu\Delta_\nu^{1/2}R^{\alpha, \nu}(\Delta_{\widehat{\Phi}}^{-1/2}\xi)^*\eta\\
&=&
J_\nu R^{\alpha, \nu}(\xi)^*\Delta_{\widehat{\Phi}}^{-1/2}\eta\\
&=&J_\nu\Lambda_\nu(<\xi, \Delta_{\widehat{\Phi}}^{-1/2}\eta>_{\alpha, \nu})
\end{eqnarray*}
from which we get (i).  
\newline
Applying (\cite{E3}2.10) to the inclusion $\alpha(N)\subset \widehat{M}$ and the operator-valued weight $\hat{T}$, we get that it is possible to construct an orthogonal $(\alpha, \nu)$-basis $(e_i)_{i\in I}$ such that $e_i=J_{\widehat{\Phi}}\Lambda_{\widehat{\Phi}}(x_i)$, with $x_i\in\gN_{\widehat{\Phi}}\cap\gN_{\widehat{\Phi}}^*\cap\gN_{\hat{T}}\cap\gN_{\hat{T}}^*$; so, $e_i$ belongs to $\mathcal D(\Delta_{\widehat{\Phi}}^{-1/2})$, and $\Delta_{\widehat{\Phi}}^{-1/2}e_i=J_{\widehat{\Phi}}\Lambda_{\widehat{\Phi}}(x_i^*)$ which belongs to $D(_\alpha H, \nu)$; which is (ii). 
\newline
Using (ii) and (i), we have :
\begin{eqnarray*}
(\theta^{\alpha, \nu}(\xi, \xi)\Delta_{\widehat{\Phi}}^{-1/2}e_i|\Delta_{\widehat{\Phi}}^{-1/2}e_i)
&=&
\|R^{\alpha, \nu}(\xi)^*\Delta_{\widehat{\Phi}}^{-1/2}e_i\|^2\\
&=&\|\Lambda_\nu(<e_i, \Delta_{\widehat{\Phi}}^{-1/2}\xi>_{\alpha, \nu}^o)\|^2\\
&=&\nu((R^{\alpha, \nu}(\Delta_{\widehat{\Phi}}^{-1/2}\xi)^*\theta^{\alpha, \nu}(e_i, e_i)R^{\alpha, \nu}(\Delta_{\widehat{\Phi}}^{-1/2}\xi)^o)
\end{eqnarray*}
and we get, using (i) and \ref{example} :
\[\sum_i(\theta^{\alpha, \nu}(\xi, \xi)\Delta_{\widehat{\Phi}}^{-1/2}e_i|\Delta_{\widehat{\Phi}}^{-1/2}e_i)=
\nu(<\Delta_{\widehat{\Phi}}^{-1/2}\xi, \Delta_{\widehat{\Phi}}^{-1/2}\xi>_{\alpha, \nu}^o)=\|\Delta_{\widehat{\Phi}}^{-1/2}\xi\|^2=\overline{\nu^o_\beta}(\theta^{\alpha, \nu}(\xi, \xi))\]
from which we get that $\sum_i\omega_{\Delta_{\widehat{\Phi}}^{-1/2}e_i}$ is a normal semi-finite weight on $\alpha(N)'$, and, by unicity of the spatial derivative, we get this weight is equal to $\overline{\nu^o_\beta}$. 
\end{proof}

%%%%psibarre2
\subsection{Theorem}
\label{psibarre2}
{\it Let $\gG=(N, M, \alpha, \beta, \Gamma, T, T', \nu)$ be a measured quantum groupoid; let $(e_i)_{i\in I}$ be an $(\alpha, \nu)$-orthogonal basis of $H$ such that, for all $i\in I$, $e_i$ belongs to $D(_\alpha H, \nu)\cap \mathcal D(\Delta_{\widehat{\Phi}}^{-1/2})$ and $\Delta_{\widehat{\Phi}}^{-1/2}e_i$ belongs to $D(_\alpha H, \nu)$; let $(N, b, A)$ be a faithful weighted right von Neumann right $N$-module, and let $\psi$ be a normal semi-finite faithful weight on $A$ lifted from $\nu^o$; then, we have, with the notations of \ref{psibarre},  \ref{propw} and \ref{example} :
\[\underline{\psi}=\sum_i\psi\underset{\nu}{_b*_\alpha}\omega_{\Delta_{\widehat{\Phi}}^{-1/2}e_i}
 =\overline{\nu^o_\beta}\circ (\psi\underset{\nu}{_b*_\alpha}id)\] }

\begin{proof}
Let $X\in (A\underset{N}{_b*_\alpha}\mathcal L(H))^+$; we have :
\begin{eqnarray*}
\sum_i\psi\underset{\nu}{_b*_\alpha}\omega_{\Delta_{\widehat{\Phi}}^{-1/2}e_i}(X)
&=&\sum_i \nu\circ b^{-1}\underset{\nu}{_b*_\alpha}\omega_{\Delta_{\widehat{\Phi}}^{-1/2}e_i}(\gT\underset{N}{_b*_\alpha}id)(X)\\
&=&(\nu\circ b^{-1}\underset{\nu}{_b*_\alpha}\overline{\nu^o_\beta})(\gT\underset{N}{_b*_\alpha}id)(X)\\
&=&\overline{\nu^o_\beta}(\nu\circ b^{-1}\underset{\nu}{_b*_\alpha}id)(\gT\underset{N}{_b*_\alpha}id)(X)\\
&=&\overline{\nu^o_\beta}\circ (\psi\underset{\nu}{_b*_\alpha}id)(X)
\end{eqnarray*}
which is the second equality, and proves therefore that $\sum_i\psi\underset{\nu}{_b*_\alpha}\omega_{\Delta_{\widehat{\Phi}}^{-1/2}e_i}$ defines a normal semi-finite faithful weight on $A\underset{N}{_b*_\alpha}\mathcal L(H)$, which does not depend on the choice of the $(\alpha, \nu)$-orthogonal basis $(e_i)_{i\in I}$. Let us denote $\psi_0$ that weight. 
\newline
We get :
\[\psi_0(\rho_\xi^{b, \alpha}aa^*(\rho_\xi^{b, \alpha})^*)
=\sum_i \psi(b(<\Delta_{\widehat{\Phi}}^{-1/2}e_i, \xi>_{\alpha, \nu}^o)^*aa^*b(<\Delta_{\widehat{\Phi}}^{-1/2}e_i, \xi>_{\alpha, \nu}^o))\]
Applying \ref{lem}(i), if $\xi$ belongs to $D(_\alpha H, \nu)\cap \mathcal D(\Delta_{\widehat{\Phi}}^{-1/2})$, and is such that $\Delta_{\widehat{\Phi}}^{-1/2}\xi$ belongs to $D(_\alpha H, \nu)$,
we get that $b(<\Delta_{\widehat{\Phi}}^{-1/2}e_i, \xi>_{\alpha, \nu}^o)^*)$ belongs to $\mathcal D(\sigma_{-i/2}^\psi)$ and that :
\[\sigma_{-i/2}^\psi(b(<\Delta_{\widehat{\Phi}}^{-1/2}e_i, \xi>_{\alpha, \nu}^o)^*)=b(\sigma_{i/2}^\nu(<\xi, \Delta_{\widehat{\Phi}}^{-1/2}e_i>_{\alpha, \nu})^o)= b(<\Delta_{\widehat{\Phi}}^{-1/2}\xi, e_i>_{\alpha, \nu}^o)\]
So, with such an hypothesis on $\xi$, and if $a$ belongs to $\gN_\psi\cap\gN_\psi^*$, we get that :
\begin{eqnarray*}
\psi_0(\rho_\xi^{b, \alpha}aa^*(\rho_\xi^{b, \alpha})^*)
&=&
\sum_i\|J_\psi b(<e_i, \Delta_{\widehat{\Phi}}^{-1/2}\xi>_{\alpha, \nu}^o)J_\psi\Lambda_\psi (a^*)\|^2\\
&=&\sum_i\|b(< \Delta_{\widehat{\Phi}}^{-1/2}\xi, e_i>_{\alpha, \nu}^o)\Delta_\psi^{1/2}\Lambda_\psi (a)\|^2\\
&=&\sum_i(b(<\Delta_{\widehat{\Phi}}^{-1/2}\xi, e_i>_{\alpha, \nu}^o<e_i, \Delta_{\widehat{\Phi}}^{-1/2}\xi>_{\alpha, \nu}^o)\Delta_\psi^{1/2}\Lambda_\psi (a)|\Delta_\psi^{1/2}\Lambda_\psi (a))\\
&=&(b(<\Delta_{\widehat{\Phi}}^{-1/2}\xi, \Delta_{\widehat{\Phi}}^{-1/2}\xi>_{\alpha, \nu}^o)\Delta_\psi^{1/2}\Lambda_\psi (a)|\Delta_\psi^{1/2}\Lambda_\psi (a))\\
&=&\|\Delta_\psi^{1/2}\Lambda_\psi (a)\underset{\nu}{_b\otimes_\alpha}\Delta_{\widehat{\Phi}}^{-1/2}\xi\|^2
\end{eqnarray*}
Using \ref{psibarre}(iii), we get that $\underline{\psi}(\rho_\xi^{b, \alpha}aa^*(\rho_\xi^{b, \alpha})^*)=\psi_0(\rho_\xi^{b, \alpha}aa^*(\rho_\xi^{b, \alpha})^*)$, for all $a$ in $\gN_\psi\cap\gN_\psi^*$ and $\xi\in D(_\alpha H, \nu)\cap \mathcal D(\Delta_{\widehat{\Phi}}^{-1/2})$, and is such that $\Delta_{\widehat{\Phi}}^{-1/2}\xi$ belongs to $D(_\alpha H, \nu)$. By polarisation, we get $\underline{\psi}(\rho_\xi^{b, \alpha}ab(\rho_\eta^{b, \alpha})^*)=\psi_0(\rho_\xi^{b, \alpha}ab(\rho_\eta^{b, \alpha})^*)$, for all $a$, $b$ in $\gN_\psi\cap\gN_\psi^*$ and $\xi$, $\eta$ in $D(_\alpha H, \nu)\cap \mathcal D(\Delta_{\widehat{\Phi}}^{-1/2})$, such that $\Delta_{\widehat{\Phi}}^{-1/2}\xi$ and $\Delta_{\widehat{\Phi}}^{-1/2}\eta$ belong to $D(_\alpha H, \nu)$. The linear set generated by such elements is an involutive algebra, whose weak closure contains, using (\cite{E5} 2.2.1) and the semi-finiteness of $\psi$, all operators of the form $\rho_{\xi_1}^{b, \alpha}c(\rho_{\xi_2}^{b, \alpha})^*$, for any $\xi_1$, $\xi_2$ in $D(_\alpha H, \nu)$ and $c$ in $A$; therefore, using \ref{subalgebra}, we get that $\underline{\psi}$ and $\psi_0$ are equal on a dense involutive algebra. 
\newline
We easily get that $\psi_0\circ\sigma_t^{\underline{\psi}}=\psi_0\circ (\sigma_t^\psi\underset{N}{_b*_\alpha}Ad \Delta_{\widehat{\Phi}}^{-it})$ is equal to $\sum_i\psi\underset{\nu}{_b*_\alpha}\omega_{\Delta_{\widehat{\Phi}}^{-1/2}\Delta_{\widehat{\Phi}}^{-it}e_i}$; the family $\Delta_{\widehat{\Phi}}^{-it}e_i$ is another $(\alpha, \nu)$-orthogonal basis of $H$, which bears the same properties as $(e_i)_{i\in I}$. As, using (i), we know that the definition of $\psi_0$ does not depend on the choice of the orthogonal $(\alpha, \nu)$-basis, we get that $\psi_0$ is invariant under $\sigma_t^{\underline{\psi}}$, and, therefore $\underline{\psi}=\psi_0$, which finishes the proof. 
\end{proof}

%%%%ex2
\subsection{Example}
\label{ex2}
Looking again at the particular example given in \ref{example}, we get, using \ref{psibarre2}, that $\underline{\nu^o}=\overline{\nu^o_\beta}$. 

%%%psibarre3
\subsection{Theorem}
\label{psibarre3}
{\it Let $\gG=(N, M, \alpha, \beta, \Gamma, T, T', \nu)$ be a measured quantum groupoid; let $(e_i)_{i\in I}$ be an $(\alpha, \nu)$-orthogonal basis of $H$ such that, for all $i\in I$, $e_i$ belongs to $D(_\alpha H, \nu)\cap \mathcal D(\Delta_{\widehat{\Phi}}^{-1/2})$ and $\Delta_{\widehat{\Phi}}^{-1/2}e_i$ belongs to $D(_\alpha H, \nu)$; let $(N, b, A)$ be a faithful weighted right von Neumann right $N$-module, and let $\gT$ be a normal semi-finite faithful operator-valued weight from $A$ onto $b(N)$; let us write $\psi=\nu^o\circ b^{-1}\circ\gT$ and $\underline{\psi}$ the normal semi-finite faithful weight on $A\underset{N}{_b*_\alpha}\mathcal L(H)$ constructed in \ref{psibarre}; for $n\in N$, let us define $a(n)=J_\psi b(n^*)J_\psi$.  Let $\xi$ be in $D(_\alpha H, \nu)\cap \mathcal D(\Delta_{\widehat{\Phi}}^{-1/2})$, such that $\Delta_{\widehat{\Phi}}^{-1/2}\xi$ belongs to $D(_\alpha H, \nu)$; let $\eta$, $\xi_1$ be in $D(_\alpha H, \nu)$, and $\xi_2\in D(H_\beta, \nu^o)$; let $z$ be in $\gN_\psi$, $\zeta$ be in $H_\psi$, $X$ be in $A\underset{N}{_b*_\alpha}\mathcal L(H)$. Then :
\newline
(i) the operator  $\rho^{b, \alpha}_\eta z(\rho^{b, \alpha}_\xi)^*$ belongs to $\gN_{\underline{\psi}}$, and we have :
\[\Lambda_{\underline{\psi}}(\rho^{b, \alpha}_\eta z(\rho^{b, \alpha}_\xi)^*)=J_{\widehat{\Phi}}\Delta_{\widehat{\Phi}}^{-1/2}\xi\underset{\nu}{_\beta\otimes_a}\Lambda_\psi(z)\underset{\nu}{_b\otimes_\alpha}\eta\]
Moreover, the linear set generated by the operators $\rho^{b, \alpha}_\eta z(\rho^{b, \alpha}_\xi)^*$, where $z$ is in $\gN_\psi$, $\eta$ is in $D(_\alpha H, \nu)$, and $\xi$ is in $D(_\alpha H, \nu)\cap \mathcal D(\Delta_{\widehat{\Phi}}^{-1/2})$, such that $\Delta_{\widehat{\Phi}}^{-1/2}\xi$ belongs to $D(_\alpha H, \nu)$, is a core for $\Lambda_{\underline{\psi}}$. 
\newline
(ii) we have : $J_{\underline{\psi}}(\xi_2\underset{\nu}{_\beta\otimes_a}\zeta\underset{\nu}{_b\otimes_\alpha}\xi_1)=J_{\widehat{\Phi}}\xi_1\underset{\nu}{_\beta\otimes_a}J_\psi\zeta\underset{\nu}{_b\otimes_\alpha}J_{\widehat{\Phi}}\xi_2$. 
\newline
(iii) we have : $\pi_{\underline{\psi}}(X)=1\underset{N}{_\beta\otimes_a}X$. 
\newline
(iv) it is possible to define a one parameter group of unitaries $\Delta_{\widehat{\Phi}}^{-it}\underset{\nu}{_\beta\otimes_a}\Delta_\psi^{it}\underset{\nu}{_b\otimes_\alpha}\Delta_{\widehat{\Phi}}^{-it}$ on $H\underset{\nu}{_\beta\otimes_a}H_\psi\underset{\nu}{_b\otimes_\alpha}H$ with natural values on elementary tensors, and $\Delta_{\underline{\psi}}^{1/2}$ is equal to its generator $\Delta_{\widehat{\Phi}}^{-1/2}\underset{\nu}{_\beta\otimes_a}\Delta_\psi^{1/2}\underset{\nu}{_b\otimes_\alpha}\Delta_{\widehat{\Phi}}^{-1/2}$. }

\begin{proof}
We have $(\rho^{b, \alpha}_\eta z(\rho^{b, \alpha}_\xi)^*)^*\rho^{b, \alpha}_\eta z(\rho^{b, \alpha}_\xi)^*=
\rho^{b, \alpha}_\xi z^*b(<\eta, \eta>_{\alpha, \nu}^o)z(\rho^{b, \alpha}_\xi)^*$, which belongs to $\gM_{\underline{\psi}}$, by \ref{psibarre}(iii). 
\newline
Let $a$ in $\gN_\psi\cap\gN_\psi^*$; let us take $\eta_1$ satisfying the same hypothesis as $\xi$. We have, using \ref{subalgebra}(iv) applied to the weight $\underline{\psi}$, and \ref{psibarre}(ii) :
\[J_{\widehat{\Phi}}\eta_1\underset{\nu}{_\beta\otimes_a}J_\psi aJ_\psi\zeta\underset{\nu}{_b\otimes_\alpha}\eta_1=
R^{\underline{\psi}}(\Delta_\psi^{1/2}\Lambda_\psi (a)\underset{\nu}{_b\otimes_\alpha}\Delta_{\widehat{\Phi}}^{-1/2}\eta_1)^*(\zeta\underset{\nu}{_b\otimes_\alpha}\xi_1)\]
and, therefore :
\begin{multline*}
(\Lambda_{\underline{\psi}}(\rho^{b, \alpha}_\eta z(\rho^{b, \alpha}_\xi)^*)|J_{\widehat{\Phi}}\eta_1\underset{\nu}{_\beta\otimes_a}J_\psi aJ_\psi\zeta\underset{\nu}{_b\otimes_\alpha}\xi_1)=\\
(\rho^{b, \alpha}_\eta z(\rho^{b, \alpha}_\xi)^*(\Delta_\psi^{1/2}\Lambda_\psi (a)\underset{\nu}{_b\otimes_\alpha}\Delta_{\widehat{\Phi}}^{-1/2}\eta_1)|\zeta\underset{\nu}{_b\otimes_\alpha}\xi_1)
\end{multline*}
which, using  \ref{lem}(i), and the definition of $\psi$, is equal to :
\begin{multline*}
(\rho^{b, \alpha}_\eta zb(<\Delta_{\widehat{\Phi}}^{-1/2}\eta_1, \xi>_{\alpha, \nu}^o)\Delta_\psi^{1/2}\Lambda_\psi (a)|\zeta\underset{\nu}{_b\otimes_\alpha}\xi_1)=\\
(zb(<\Delta_{\widehat{\Phi}}^{-1/2}\eta_1, \xi>_{\alpha, \nu}^o)\Delta_\psi^{1/2}\Lambda_\psi (a)\underset{\nu}{_b\otimes_\alpha}\eta|\zeta\underset{\nu}{_b\otimes_\alpha}\xi_1)=\\
(zb(\sigma^\nu_{-i/2}(<\eta_1, \Delta_{\widehat{\Phi}}^{-1/2}\xi>_{\alpha, \nu}^o))\Delta_\psi^{1/2}\Lambda_\psi (a)\underset{\nu}{_b\otimes_\alpha}\eta|\zeta\underset{\nu}{_b\otimes_\alpha}\xi_1)=\\
(z\Delta_\psi^{1/2}b(<\eta_1, \Delta_{\widehat{\Phi}}^{-1/2}\xi>_{\alpha, \nu}^o)\Lambda_\psi (a)\underset{\nu}{_b\otimes_\alpha}\eta|\zeta\underset{\nu}{_b\otimes_\alpha}\xi_1)
\end{multline*}
Let us suppose that $z$ belongs to $\mathcal D(\sigma_{i/2}^\psi)$; we get that :
\begin{align*}
z\Delta_\psi^{1/2}b(<\eta_1, \Delta_{\widehat{\Phi}}^{-1/2}\xi>_{\alpha, \nu}^o)\Lambda_\psi (a)
&=
\Delta_\psi^{1/2}\sigma^\psi_{i/2}(z)b(<\eta_1, \Delta_{\widehat{\Phi}}^{-1/2}\xi>_{\alpha, \nu}^o)\Lambda_\psi (a)\\
&=
J_\psi\Lambda_\psi((a^*b(<\Delta_{\widehat{\Phi}}^{-1/2}\xi, \eta_1>_{\alpha, \nu}^o)\sigma_{-i/2}(z^*))\\
&=
J_\psi a^*b(<\Delta_{\widehat{\Phi}}^{-1/2}\xi, \eta_1>_{\alpha, \nu}^o)J_\psi\Lambda_\psi(z)\\
&=J_\psi a^*J_\psi a(<\eta_1, \Delta_{\widehat{\Phi}}^{-1/2}\xi>_{\alpha, \nu}^o)\Lambda_\psi(z)\\
&=J_\psi a^*J_\psi a(<J_{\widehat{\Phi}}\Delta_{\widehat{\Phi}}^{-1/2}\xi, J_{\widehat{\Phi}}\eta_1>_{\beta, \nu^o})\Lambda_\psi(z)
\end{align*}
which remains true for all $z\in\gN_\psi$; therefore, we then get that :
\begin{multline*}
(\Lambda_{\underline{\psi}}(\rho^{b, \alpha}_\eta z(\rho^{b, \alpha}_\xi)^*)|J_{\widehat{\Phi}}\eta_1\underset{\nu}{_\beta\otimes_a}J_\psi aJ_\psi\zeta\underset{\nu}{_b\otimes_\alpha}\xi_1)=\\
(J_\psi a^*J_\psi a(<J_{\widehat{\Phi}}\Delta_{\widehat{\Phi}}^{-1/2}\xi, J_{\widehat{\Phi}}\eta_1>_{\beta, \nu^o})\Lambda_\psi(z)\underset{\nu}{_b\otimes_\alpha}\eta|\zeta\underset{\nu}{_b\otimes_\alpha}\xi_1)=\\(a(<J_{\widehat{\Phi}}\Delta_{\widehat{\Phi}}^{-1/2}\xi, J_{\widehat{\Phi}}\eta_1>_{\beta, \nu^o})\Lambda_\psi(z)\underset{\nu}{_b\otimes_\alpha}\eta|J_\psi aJ_\psi \zeta\underset{\nu}{_b\otimes_\alpha}\xi_1)=\\
(J_{\widehat{\Phi}}\Delta_{\widehat{\Phi}}^{-1/2}\xi\underset{\nu}{_\beta\otimes_a}\Lambda_\psi(z)\underset{\nu}{_b\otimes_\alpha}\eta|J_{\widehat{\Phi}}\eta_1\underset{\nu}{_\beta\otimes_a}J_\psi aJ_\psi \zeta\underset{\nu}{_b\otimes_\alpha}\xi_1)
\end{multline*}
from which, by density, we get the first result of (i). 
\newline
Using \ref{psibarre}(ii), we get that $\sigma_t^{\underline{\psi}}(\rho^{b, \alpha}_\eta z(\rho^{b, \alpha}_\xi)^*)=\rho^{b, \alpha}_{\Delta_{\widehat{\Phi}}^{-it}\eta}\sigma_t^\psi(z)(\rho^{b, \alpha}_{\Delta_{\widehat{\Phi}}^{-it}\xi})^*$;  so, the linear set generated by the operators $\rho^{b, \alpha}_\eta z(\rho^{b, \alpha}_\xi)^*$, where $z$ belongs to $\gN_\psi\cap\gN_\psi^*$, and $\xi$ (resp. $\eta$) is in $D(_\alpha H, \nu)\cap \mathcal D(\Delta_{\widehat{\Phi}}^{-1/2})$, such that $\Delta_{\widehat{\Phi}}^{-1/2}\xi$ (resp. $\Delta_{\widehat{\Phi}}^{-1/2}\eta$) belongs to $D(_\alpha H, \nu)$ is a $*$-subalgebra of $\gN_{\underline{\psi}}\cap \gN_{\underline{\psi}}^*$, dense in $A\underset{N}{_b*_\alpha}\mathcal L(H)$ by \ref{subalgebra}, and globally invariant under $\sigma_t^{\underline{\psi}}$. It is possible to put on the image of this algebra under $\Lambda_{\underline{\psi}}$ a structure of left-Hilbert algebra, which, in turn, leads to a faithful normal semi-finite weight $\psi_0$ on $A\underset{N}{_b*_\alpha}\mathcal L(H)$, equal to $\underline{\psi}$ on this subalgebra, and invariant under $\sigma_t^{\underline{\psi}}$. So, we get $\psi_0=\underline{\psi}$, which finishes the proof of (i). 
\newline
On the other hand, let's apply \ref{subalgebra}(iii) to the standard representation of $A\underset{N}{_b*_\alpha}\mathcal L(H)$ given by the weight $\underline{\psi}$, and we get (ii) and (iii). 
\newline
Let now $\xi\in D(H_\beta, \nu^o)$; we have, for all $t\in\mathbb{R}$, $n\in\gN_\nu$ :
\begin{align*}
\beta(n^*)\Delta_{\widehat{\Phi}}^{-it}\xi
&=\Delta_{\widehat{\Phi}}^{-it}\tau_t(\beta(n^*))\xi\\
&=\Delta_{\widehat{\Phi}}^{-it}\beta(\sigma_t^\nu(n^*))\xi\\
&=\Delta_{\widehat{\Phi}}^{-it}R^{\beta, \nu^o}(\xi)J_\nu\Lambda_\nu(\sigma_t^\nu(n))\\
&=\Delta_{\widehat{\Phi}}^{-it}R^{\beta, \nu^o}(\xi)J_\nu\Delta_\nu^{it}\Lambda_\nu(n)
\end{align*}
and, therefore, $\Delta_{\widehat{\Phi}}^{-it}\xi$ belongs to $D(H_\beta, \nu^o)$, and $R^{\beta, \nu^o}(\Delta_{\widehat{\Phi}}^{-it}\xi)=\Delta_{\widehat{\Phi}}^{-it}R^{\beta, \nu^o}(\xi)\Delta_\nu^{it}$, and $<\Delta_{\widehat{\Phi}}^{-it}\xi, \Delta_{\widehat{\Phi}}^{-it}\xi>_{\beta, \nu^o}=\sigma_{-t}^\nu(<\xi, \xi>_{\beta, \nu^o})$. Therefore, if $\xi'$ belongs to $D(_\alpha H, \nu)$, $\eta\in H_\psi$, we get :
\begin{align*}
\|\Delta_{\widehat{\Phi}}^{-it}\xi\underset{\nu}{_\beta\otimes_a}\Delta_\psi^{it}\eta\underset{\nu}{_b\otimes_\alpha}\Delta_{\widehat{\Phi}}^{-it}\xi'\|^2
&=
(b(<\Delta_{\widehat{\Phi}}^{-it}\xi', \Delta_{\widehat{\Phi}}^{-it}\xi'>_{\alpha, \nu}^o)a(<\Delta_{\widehat{\Phi}}^{-it}\xi, \Delta_{\widehat{\Phi}}^{-it}\xi>_{\beta, \nu^o})\Delta_\psi^{it}\eta|\Delta_\psi^{it}\eta)\\
&=(b(\sigma_{-t}^{\nu^o}(<\xi', \xi'>_{\alpha, \nu}^o))a(\sigma_{-t}^\nu(<\xi, \xi>_{\beta, \nu^o}))\Delta_\psi^{it}\eta|\Delta_\psi^{it}\eta)\\
&=(\sigma_t^\psi(b(<\xi', \xi'>_{\alpha, \nu}^o))J_\psi b(\sigma_{-t}^\nu(<\xi, \xi>_{\beta, \nu^o}))J_\psi\Delta_\psi^{it}\eta|\Delta_\psi^{it}\eta)\\
&=(b(<\xi', \xi'>_{\alpha, \nu}^o)\Delta_\psi^{-it}J_\psi\sigma_t^\psi(b(<\xi, \xi>_{\beta, \nu^o}))J_\psi\Delta_\psi^{it}\eta|\eta)\\
&=(b(<\xi', \xi'>_{\alpha, \nu}^o)J_\psi b(<\xi, \xi>_{\beta, \nu^o})J_\psi\eta|\eta)\\
&=(b(<\xi', \xi'>_{\alpha, \nu}^o)a(<\xi, \xi>_{\beta, \nu^o})\eta|\eta)\\
&=\|\xi\underset{\nu}{_\beta\otimes_a}\eta\underset{\nu}{_b\otimes_\alpha}\xi'\|^2
\end{align*}
Now, from (i) and (ii), we get  that the infinitesimal generator $\Delta_{\widehat{\Phi}}^{-1/2}\underset{\nu}{_\beta\otimes_a}\Delta_\psi^{1/2}\underset{\nu}{_b\otimes_\alpha}\Delta_{\widehat{\Phi}}^{-1/2}$ of this one-parameter of unitaries is included in $\Delta_{\underline{\psi}}^{1/2}$; these operators being self-adjoint, we get (iv). \end{proof}

%%%psia
\subsection{Corollary}
\label{psia}
{\it Let $\gG=(N, M, \alpha, \beta, \Gamma, T, T', \nu)$ be a measured quantum groupoid, $(b,\ga)$ an action of $\gG$ on a von Neumann algebra $A$, $\psi$ a normal semi-finite faithful weight on $A$, $\overline{\psi_\ga}$ the normal semi-finite faithful weight constructed on $A\underset{N}{_b*_\alpha}\mathcal L(H)$ by transporting the bidual weight. Then, for any $x\in \gN_\psi$, $\xi\in D(_\alpha H, \nu)$, $\eta\in D(_\alpha H, \nu)\cap \mathcal D(\Delta_{\widehat{\Phi}}^{-1/2})$ such that  $\Delta_{\widehat{\Phi}}^{-1/2}\eta$ belongs to $D(H_\beta, \nu^o)$, the operator $(1\underset{N}{_b\otimes_\alpha}\theta^{\alpha, \nu}(\xi, \eta))\ga(x)$ belongs to $\gN_{\overline{\psi_\ga}}$, and we have:}
\[\overline{\psi_\ga}(\ga(x^*)((1\underset{N}{_b\otimes_\alpha}\theta^{\alpha, \nu}(\xi, \eta)^*\theta^{\alpha, \nu}(\xi, \eta))\ga(x))=\|\Lambda_\psi(x)\underset{\nu}{_b\otimes_\alpha}J_{\widehat{\Phi}}\Delta_{\widehat{\Phi}}^{-1/2}\eta\underset{\nu}{_\beta\otimes_\alpha}\xi\|^2\]

\begin{proof}
Using \ref{psibarre3} applied to $\underline{\nu^o}$, we get that $\Lambda_{\underline{\nu^o}}(\theta^{\alpha, \nu}(\xi, \eta))=J_{\widehat{\Phi}}\Delta_{\widehat{\Phi}}^{-1/2}\eta\underset{\nu}{_\beta\otimes_\alpha}\xi$, which belongs to $D(_\alpha H_{\underline{\nu^o}}, \nu)$; so, using \ref{example} and \ref{ex2}, we get the result. 
\end{proof}

%%%%%
\section{Standard implementation: using the weight $\underline{\psi}$. }
\label{using}
In that section, we calculate (\ref{propW*sigma}) the dual weight $\widetilde{(\underline{\psi})}$ of $\underline{\psi}$, with respect to the action $\underline{\ga}$ (\ref{psitildetheta}(ii)); this will allow us to calculate $J_{\widetilde{(\underline{\psi})}}$ (\ref{propW*sigma}), and then, to obtain a formula linking $U^\ga_\psi$ and $U^{\underline{\ga}}_{\underline{\psi}}$ (\ref{Upsibarre}). As $U^{\underline{\ga}}_{\underline{\psi}}$ is a corepresentation by \ref{cora}, we obtain then that $U^\ga_\psi$ is a corepresentation (and, therefore, a standard implementation) whenever it is possible to construct $\underline{\psi}$ (\ref{standard}).

%%%%tildeTheta
\subsection{Proposition}
\label{tildeTheta}
{\it Let $\gG=(N, M, \alpha, \beta, \Gamma, T, T', \nu)$ be a measured quantum groupoid; let $A$ be a von Neumann algebra acting on a Hilbert space $\gH$, $(b, \mathfrak a)$ be an action of $\gG$ on $A$, $(1\underset{N}{_b\otimes_\alpha}\hat{\beta}, \underline{\mathfrak a})$ be the action of $\gG$ on $A\underset{N}{_b*_\alpha}\mathcal L(H)$ introduced in \ref{crossed}; then, let us write, for any $Y$ in $\mathcal L(\gH\underset{\nu}{_b\otimes_\alpha}H\underset{\nu}{_{\hat{\beta}}\otimes_\alpha}H)$, 
\[\tilde{\Theta}(Y)=(1\underset{N}{_b\otimes_\alpha}W)^*(id\underset{N}{_b*_\alpha}\varsigma_N)(Y)(1\underset{N}{_b\otimes_\alpha}W)\]
which belongs to $\mathcal L(\gH\underset{\nu}{_b\otimes_\alpha}H\underset{\nu}{_\beta\otimes_\alpha}H)$; then, we have :
\newline
(i) for any $X\in A\underset{N}{_b*_\alpha}\mathcal L(H)$, $\tilde{\Theta}(\underline{\mathfrak a}(X))=(\mathfrak a\underset{N}{_b*_\alpha}id)(X)$ and :
\[\tilde{\Theta}((A\underset{N}{_b*_\alpha}\mathcal L(H))\rtimes_{\underline{\mathfrak a}}\gG)=
(A\rtimes_\mathfrak a\gG)\underset{N}{_\beta*_\alpha}\mathcal L(H)\]
\newline
(ii) $(1\underset{N}{_b\otimes_\alpha}\hat{\alpha}, (id\underset{N}{_b*_\alpha}\varsigma_N)(\tilde{\mathfrak a}\underset{N}{_\beta*_\alpha}id))$ is an action of $\hat{\gG}^c$ on $(A\rtimes_\mathfrak a\gG)\underset{N}{_\beta*_\alpha}\mathcal L(H)$, and :
\[(\tilde{\Theta}\underset{N^o}{_\alpha*_\beta}id)\widetilde{(\underline{\mathfrak a})}=(id\underset{N}{_b*_\alpha}\varsigma_{N^o})(\tilde{\mathfrak a}\underset{N^o}{_{\hat{\alpha}}*_\beta}id)\tilde{\Theta}\]
where $\widetilde{(\underline{\ga})}$ is the dual action  of $\underline{\ga}$ (it is therefore an action of $\widehat{\gG}^c$ on $(A\underset{N}{_b*_\alpha}\mathcal L(H))\rtimes_{\underline{\ga}}\gG$).}

\begin{proof}
By the definition of $\underline{\ga}$, we get the first formula of (i). The second formula of (i) was already proved in (\cite{E5} 11.4). 
Moreover, using (i), we have, for all $X\in A\underset{N}{_b*_\alpha}\mathcal L(H)$ :
\begin{eqnarray*}
(\tilde{\Theta}\underset{N^o}{_{\hat{\alpha}}*_\beta}id)\widetilde{(\underline{\mathfrak a})}(\underline{\mathfrak a}(X))&=&
(\tilde{\Theta}\underset{N^o}{_{\hat{\alpha}}*_\beta}id)(\underline{\mathfrak a}(X)\underset{N^o}{_{\hat{\alpha}}\otimes_\beta}1)\\
&=&
(\mathfrak a\underset{N}{_b*_\alpha}id)(X)\underset{N^o}{_{\hat{\alpha}}\otimes_\beta}1\\
&=&(id\underset{N}{_b*_\alpha}\varsigma_{N^o})(\tilde{\mathfrak a}\underset{N^o}{_{\hat{\alpha}}*_\beta}id)(\mathfrak a\underset{N}{_b*_\alpha}id)(X)\\
&=&(id\underset{N}{_b*_\alpha}\varsigma_{N^o})(\tilde{\mathfrak a}\underset{N^o}{_{\hat{\alpha}}*_\beta}id)\tilde{\Theta}(\underline{\mathfrak a}(X))
\end{eqnarray*}
and, for all $z\in \widehat{M}'$, we have :
\[(\tilde{\Theta}\underset{N^o}{_{\hat{\alpha}}*_\beta}id)\widetilde{(\underline{\mathfrak a})}(1\underset{N}{_b\otimes_\alpha}1\underset{N^o}{_{\hat{\alpha}}\otimes_\beta}z)
=
(\tilde{\Theta}\underset{N^o}{_{\hat{\alpha}}*_\beta}id)(1\underset{N}{_b\otimes_\alpha}\widehat{\Gamma}^c(z))\]
which, thanks again to (\cite{E5} 11.4), is equal to :
\[(1\underset{N}{_b\otimes_\alpha}1\underset{N}{_\beta\otimes_{\hat{\alpha}}}J_\Phi J_{\widehat{\Phi}}
\underset{N}{_\beta\otimes_{\hat{\alpha}}}1)(\widehat{\Gamma}^{oc}\underset{N^o}{_{\hat{\alpha}}*_\beta}id)\widehat{\Gamma}^c(z)(1\underset{N}{_b\otimes_\alpha}1\underset{N}{_\beta\otimes_\alpha}J_{\widehat{\Phi}}J_\Phi 
\underset{N^o}{_{\hat{\alpha}}\otimes_\beta}1)\]
and we have :
\begin{multline*}
(id\underset{N}{_b*_\alpha}\varsigma_{N^o})(\tilde{\mathfrak a}\underset{N^o}{_{\hat{\alpha}}*_\beta}id)\tilde{\Theta}(1\underset{N}{_b\otimes_\alpha}1\underset{N^o}{_{\hat{\alpha}}\otimes_\beta}z)=\\
(id\underset{N}{_b*_\alpha}\varsigma_{N^o})(\widehat{\Gamma}^c\underset{N}{_\beta*_\alpha}id)[(1\underset{N}{_\beta\otimes_{\hat{\alpha}}}J_\Phi J_{\widehat{\Phi}})\widehat{\Gamma}^{oc}(z)(1\underset{N}{_\beta\otimes_\alpha}J_{\widehat{\Phi}}J_\Phi)]
\end{multline*}
from which we deduce that :
\[(\tilde{\Theta}\underset{N^o}{_{\hat{\alpha}}*_\beta}id)\widetilde{(\underline{\mathfrak a})}(1\underset{N}{_b\otimes_\alpha}1\underset{N^o}{_{\hat{\alpha}}\otimes_\beta}z)=
(id\underset{N}{_b*_\alpha}\varsigma_{N^o})(\tilde{\mathfrak a}\underset{N^o}{_{\hat{\alpha}}*_\beta}id)\tilde{\Theta}(1\underset{N}{_b\otimes_\alpha}1\underset{N^o}{_{\hat{\alpha}}\otimes_\beta}z)\]
and we get that :
\[(\tilde{\Theta}\underset{N^o}{_\alpha*_\beta}id)\widetilde{(\underline{\mathfrak a})}=(id\underset{N}{_b*_\alpha}\varsigma_{N^o})(\tilde{\mathfrak a}\underset{N^o}{_{\hat{\alpha}}*_\beta}id)\tilde{\Theta}\]
from which we deduce that $(1\underset{N}{_b\otimes_\alpha}\hat{\alpha}, (id\underset{N}{_b*_\alpha}\varsigma_N)(\tilde{\mathfrak a}\underset{N}{_\beta*_\alpha}id))$ is an action of $\hat{\gG}^o$ on the von Neumann algebra $(A\rtimes_\mathfrak a\gG)\underset{N}{_\beta*_\alpha}\mathcal L(H)$, which finishes the proof. 

\end{proof}

%%%%cortildeTheta
\subsection{Corollary}
\label{cortildeTheta}
{\it Let $\gG=(N, M, \alpha, \beta, \Gamma, T, T', \nu)$ be a measured quantum groupoid; let $A$ be a von Neumann algebra acting on a Hilbert space $\gH$, $(b, \mathfrak a)$ an action of $\gG$ on $A$, $(1\underset{N}{_b\otimes_\alpha}\hat{\beta}, \underline{\mathfrak a})$ be the action of $\gG$ on $A\underset{N}{_b*_\alpha}\mathcal L(H)$ introduced in \ref{crossed} and $\tilde{\Theta}$ the isomorphism introduced in \ref{tildeTheta} which sends 
$(A\underset{N}{_b*_\alpha}\mathcal L(H))\rtimes_{\underline{\mathfrak a}}\gG$ onto 
$(A\rtimes_\mathfrak a\gG)\underset{N}{_\beta*_\alpha}\mathcal L(H)$; then, we have $\tilde{\Theta}\circ T_{\widetilde{(\underline{\ga})}}=(T_{\tilde{\ga}}\underset{N}{_\beta*_\alpha}id)\Tilde{\Theta}$. }

\begin{proof}
Using \ref{tildeTheta}(ii), we get :
\begin{align*}
\tilde{\Theta}\circ T_{\tilde{\underline{\ga}}}
&=
(id\underset{N}{_b*_\alpha}id\underset{N}{_{\hat{\alpha}}*_\beta}\widehat{\Phi}^c)(\tilde{\Theta}\underset{N}{_{\hat{\alpha}}*_\beta}id)\tilde{\underline{\ga}}\\
&=
(id\underset{N}{_b*_\alpha}id\underset{N}{_{\hat{\alpha}}*_\beta}\widehat{\Phi}^c)(id\underset{N}{_b*_\alpha}\varsigma_{N^o})(\tilde{\mathfrak a}\underset{N^o}{_{\hat{\alpha}}*_\beta}id)\tilde{\Theta}\\
&=
((id\underset{N}{_{\hat{\alpha}}*_\beta}\widehat{\Phi}^c)\tilde{\ga}\underset{N}{_\beta*_\alpha}id)\Tilde{\Theta}\\
&=(T_{\tilde{\ga}}\underset{N}{_\beta*_\alpha}id)\Tilde{\Theta}
\end{align*}
which is the result. \end{proof}

%%%%psitildetheta
\label{psitildetheta}
\subsection{Theorem}
{\it Let $\gG=(N, M, \alpha, \beta, \Gamma, T, T', \nu)$ be a measured quantum groupoid; let $A$ be a von Neumann algebra acting on a Hilbert space $\gH$, $(b, \mathfrak a)$ a weighted action of $\gG$ on $A$, $(1\underset{N}{_b\otimes_\alpha}\hat{\beta}, \underline{\mathfrak a})$ be the action of $\gG$ on $A\underset{N}{_b*_\alpha}\mathcal L(H)$ introduced in \ref{crossed} and $\tilde{\Theta}$ the isomorphism introduced in \ref{tildeTheta} which sends 
$(A\underset{N}{_b*_\alpha}\mathcal L(H))\rtimes_{\underline{\mathfrak a}}\gG$ onto 
$(A\rtimes_\mathfrak a\gG)\underset{N}{_\beta*_\alpha}\mathcal L(H)$; then :
\newline
(i) $(N, 1\underset{N}{_b\otimes_\alpha}\beta, A\rtimes_\ga\gG)$ is a von Neumann faithful right $N$-module; let $\psi$ be a lifted weight on $A$, then $\tilde{\psi}$ is a lifted weight on $A\rtimes_\ga\gG$. Let's denote then $\underline{\psi}$ and $\underline{(\tilde{\psi})}$ the weights constructed by \ref{psibarre} applied to $\psi$ and $\tilde{\psi}$. 
\newline
(ii) we have $\underline{(\tilde{\psi})}\circ\tilde{\Theta}=\widetilde{(\underline{\psi})}$ and, for all $t\in\mathbb{R}$, $\sigma_t^{\underline{(\tilde{\psi})}}\circ\tilde{\Theta}=\tilde{\Theta}\circ\sigma_t^{\widetilde{(\underline{\psi})}}$. 
\newline
(iii) moreover, $\overline{\psi_\ga}$ is a lifted weight on $A\underset{N}{_b*_\alpha}\mathcal L(H)$, and we can define a normal semifinite faithful weight $\underline{(\overline{\psi_\ga})}$ on $A\underset{N}{_b*_\alpha}\mathcal L(H)\underset{N}{_\beta*_\alpha}\mathcal L(H)$ . On the other hand, we can define the normal semi-finite faithful weight $\overline{(\underline{\psi})_{\underline{\ga}}}$ on $A\underset{N}{_b*_\alpha}\mathcal L(H)\underset{N}{_{\hat{\beta}}*_\alpha}\mathcal L(H)$. Then, we have $\underline{(\overline{\psi_\ga})}\circ\tilde{\Theta}=\overline{(\underline{\psi})_{\underline{\ga}}}$.

}

\begin{proof}
Let $\gT$ be a normal faithful semi-finite operator valued weight from $A$ into $b(N)$; then $\ga\circ\gT\circ \ga^{-1}$ is a normal faithful semi-finite operator valued weight from $\ga(A)$ into $1\underset{N}{_b\otimes_\alpha}\beta(N)$, and $\ga\circ\gT\circ \ga^{-1}\circ T_{\tilde{\ga}}$ is a normal faithful semi-finite operator-valued weight from $A\rtimes_\ga\gG$ into $1\underset{N}{_b\otimes_\alpha}\beta(N)$; then, if we write $\psi=\nu^o\circ b^{-1}\circ\gT$, the dual weight $\tilde{\psi}$ can be written as $\nu^o\circ (1\underset{N}{_b\otimes_\alpha}\beta)^{-1}\circ (\ga\circ\gT\circ \ga^{-1}\circ T_{\tilde{\ga}})$, which finishes the proof of (i). 
\newline
We have then, using the notations of \ref{psibarre}, and results \ref{cortildeTheta} and \ref{tildeTheta}(i) :
\begin{align*}
\underline{(\tilde{\psi})}\circ\tilde{\Theta}
&=
\sum_i(\tilde{\psi}\underset{\nu}{_\beta*_\alpha}\omega_{\Delta_{\widehat{\Phi}}^{-1/2}e_i})\circ\tilde{\Theta}\\
&=
\sum_i(\psi\circ\ga^{-1}\circ T_{\tilde{a}}\underset{\nu}{_\beta*_\alpha}\omega_{\Delta_{\widehat{\Phi}}^{-1/2}e_i})\circ\tilde{\Theta}\\
&=
\sum_i(\psi\underset{\nu}{_b*_\alpha}\omega_{\Delta_{\widehat{\Phi}}^{-1/2}e_i})\circ(\ga\underset{N}{_b*_\alpha}id)^{-1}\circ (T_{\tilde{\ga}}\underset{N}{_\beta*_\alpha}id)\circ\tilde{\Theta}\\
&=
\underline{\psi}\circ(\ga\underset{N}{_b*_\alpha}id)^{-1}\circ \tilde{\Theta}\circ T_{\tilde{\underline{\ga}}}\\
&=
\underline{\psi}\circ(\underline{\ga})^{-1}\circ T_{\tilde{\underline{\ga}}}\\
&=
\widetilde{(\underline{\psi})}
\end{align*}
which finishes the proof of (ii). 
\newline
We have :
\[\overline{\psi_\ga}=\nu^o\circ (1\underset{N}{_b\otimes_\alpha}\beta)^{-1}\circ (\ga\circ\gT\circ \ga^{-1}\circ T_{\tilde{\ga}})\circ T_{\underline{\ga}}\]
So, by composition of operator-valued weights, we get that $\overline{\psi_\ga}$ is a lifted weight on the faithful right $N$-module $(N, A\underset{N}{_b*_\alpha}\mathcal L(H), 1\underset{N}{_b\otimes_\alpha}\beta)$, and, applying \ref{psibarre}, we can construct the normal semi-finite faithful weight $\underline{(\overline{\psi_\ga})}$ on $A\underset{N}{_b*_\alpha}\mathcal L(H)\underset{N}{_\beta*_\alpha}\mathcal L(H)$.
\newline
On the other hand, as $\underline{\psi}$ is a normal semi-finite faithful weight on $A\underset{N}{_b*_\alpha}\mathcal L(H)$, and as $(1\underset{N}{_b\otimes_\alpha}\hat{\beta}, \underline{\ga})$ (\ref{crossed}) is an action of $\gG$ on $A\underset{N}{_b*_\alpha}\mathcal L(H)$, we can define (\ref{crossed}) a weight $\overline{(\underline{\psi})_{\underline{\ga}}}$ on $A\underset{N}{_b*_\alpha}\mathcal L(H)\underset{N}{_{\hat{\beta}}*_\alpha}\mathcal L(H)$. As $\tilde{\Theta}$ is an isomorphism from $A\underset{N}{_b*_\alpha}\mathcal L(H)\underset{N}{_{\hat{\beta}}*_\alpha}\mathcal L(H)$ onto $A\underset{N}{_b*_\alpha}\mathcal L(H)\underset{N}{_\beta*_\alpha}\mathcal L(H)$, we can define then another normal semi-finite faithful weight $\underline{(\overline{\psi_\ga})}\circ\tilde{\Theta}$ on $A\underset{N}{_b*_\alpha}\mathcal L(H)\underset{N}{_{\hat{\beta}}*_\alpha}\mathcal L(H)$. 
\newline
Let's represent $A$ on $H_\psi$ and consider the isomorphism $\tilde{\Theta}$ from $\mathcal L(H_\psi\underset{\nu}{_b\otimes_\alpha}H\underset{\nu}{_{\hat{\beta}}\otimes_\alpha}H)$ onto $\mathcal L(H_\psi\underset{\nu}{_b\otimes_\alpha}H\underset{\nu}{_\beta\otimes_\alpha}H)$. The commutant of $A\underset{N}{_b*_\alpha}\mathcal L(H)\underset{N}{_\beta*_\alpha}\mathcal L(H)$ on the Hilbert space $H_\psi\underset{\nu}{_b\otimes_\alpha}H\underset{\nu}{_\beta\otimes_\alpha}H$ is $A'\underset{N}{_b\otimes_\alpha}1_H\underset{N}{_\beta\otimes_\alpha}1_H$, which is isomorphic to $A^o$. Let us consider the spatial derivative $\frac{d\overline{(\underline{\psi})_{\underline{\ga}}}\circ\tilde{\Theta}^{-1}}{d\psi^o}$ on $H_\psi\underset{\nu}{_b\otimes_\alpha}H\underset{\nu}{_{\hat{\beta}}\otimes_\alpha}H$. As, for $x\in A'$, $\tilde{\Theta}$ sends $x\underset{N}{_b\otimes_\alpha}1_H\underset{N}{_{\hat{\beta}}\otimes_\alpha}1_H$ on $x\underset{N}{_b\otimes_\alpha}1_H\underset{N}{_\beta\otimes_\alpha}1_H$, we get that :
\[\frac{d\overline{(\underline{\psi})_{\underline{\ga}}}\circ\tilde{\Theta}^{-1}}{d\psi^o}=\tilde{\Theta}(\frac{d\overline{(\underline{\psi})_{\underline{\ga}}}}{d\psi^o})\]
where the spatial derivative $\frac{d\overline{(\underline{\psi})_{\underline{\ga}}}}{d\psi^o}$ is taken on the Hilbert space $H_\psi\underset{\nu}{_b\otimes_\alpha}H\underset{\nu}{_{\hat{\beta}}\otimes_\alpha}H$. 
But, (using \cite{St} 12.11), we get that :
\[\frac{d\overline{(\underline{\psi})_{\underline{\ga}}}}{d\psi^o}=\frac{d\widetilde{(\underline{\psi})}\circ T_{\underline{\ga}}}{d\psi^o}=\frac{d\widetilde{(\underline{\psi})}}{d\tilde{\psi}^o}\]
where we write, for simplification, $\tilde{\psi}^o$ for the weight taken on $(A\underset{N}{_b*_\alpha}\mathcal L(H)\rtimes_{\underline{\ga}}\gG)'$, whose image by $\tilde{\Theta}$ is, thanks to (i), equal to $(A\rtimes_a\gG)'\underset{N}{_\beta\otimes_\alpha}1_H$. 
\newline
Therefore, using (ii), we get that :
\[\tilde{\Theta}(\frac{d\overline{(\underline{\psi})_{\underline{\ga}}}}{d\psi^o})=\frac{d\widetilde{(\underline{\psi})}\circ\Theta^{-1}}{d\tilde{\psi}^o}=\frac{d\underline{(\tilde{\psi})}}{d\tilde{\psi}^o}=\frac{d\underline{(\overline{\psi_\ga})}}{d\psi^o}\]
which gives the result. \end{proof}

%%%%lemW
\subsection{Lemma}
\label{lemW}
{\it Let $\gG=(N, M, \alpha, \beta, \Gamma, T, T', \nu)$ be a measured quantum groupoid, $W$ its pseudo-multiplicative unitary, $(e_i)_{i\in I}$ an orthogonal $(\alpha, \nu)$-basis of $H$; then, we have, for all $a\in\gN_{\widehat{\Phi}^c}\cap\gN_{\hat{T}^c}$, $\zeta\in D(_\alpha H, \nu)\cap D(H_{\hat{\beta}}, \nu^o)$ :
\[\sum_i\Lambda_{\widehat{\Phi}^c}(\omega_{ J_{\widehat{\Phi}}J_\Phi\zeta, J_{\widehat{\Phi}}J_\Phi e_i}\underset{N^o}{_{\hat{\alpha}}*_\beta}id)\widehat{\Gamma}^c(a)\underset{\nu}{_\beta\otimes_\alpha}e_i=
W^*(\Lambda_{\widehat{\Phi}^c}(a)\underset{\nu^o}{_\alpha\otimes_{\hat{\beta}}}\zeta)\]}

\begin{proof}
Let us first remark that $J_\Phi J_{\widehat{\Phi}}\zeta$ and $J_\Phi J_{\widehat{\Phi}}e_i$ belong to $D(_{\hat{\alpha}}H, \nu)$, and that $\Lambda_{\widehat{\Phi}^c}(a)$ belongs, thanks to (\cite{E5} 2.2), to $D(_\alpha H, \nu)$. Applying then the definition (\cite{E5} 3.6 (i)) of the pseudo-multiplicative unitary $W^{^c}$ of the measured quantum group $\widehat{\gG}^c$, we get that :
\[\Lambda_{\widehat{\Phi}^c}((\omega_{J_{\widehat{\Phi}}J_\Phi \zeta,  J_{\widehat{\Phi}}J_\Phi e_i}\underset{N^o}{_{\hat{\alpha}}*_\beta}id)\widehat{\Gamma}^c(a))=(\omega_{J_{\widehat{\Phi}}J_\Phi \zeta, J_{\widehat{\Phi}}J_\Phi e_i}*id)(\widehat{W}{^ c}^*)\Lambda_{\widehat{\Phi}^c}(a)\]
As $\widehat{W}{^ c}^*=(\widehat{W^o})^*=\sigma W^o\sigma$, we get :
\[(\omega_{ J_{\widehat{\Phi}}J_\Phi\zeta,  J_{\widehat{\Phi}}J_\Phi e_i}*id)(\widehat{W}{^ c}^*)=
(id*\omega_{J_{\widehat{\Phi}}J_\Phi \zeta, J_{\widehat{\Phi}}J_\Phi e_i})(W^o)\]
and, using \cite{E5} 3.12 (v) and 3.11(iii), we get :
\[(id*\omega_{ J_{\widehat{\Phi}}J_\Phi\zeta, J_{\widehat{\Phi}J_\Phi }e_i})(W^o)=
J_{\widehat{\Phi}}(id*\omega_{J_\Phi\zeta, J_\Phi e_i})(W)J_{\widehat{\Phi}}=(id*\omega_{\zeta, e_i})(W^*)\]
from which we get the result. 
\end{proof}

%%%%propW*sigma
\subsection{Proposition}
\label{propW*sigma}
{\it Let $\gG=(N, M, \alpha, \beta, \Gamma, T, T', \nu)$ be a measured quantum groupoid; let $A$ be a von Neumann algebra, $(b, \mathfrak a)$ a weighted action of $\gG$ on $A$, $(1\underset{N}{_b\otimes_\alpha}\hat{\beta}, \underline{\mathfrak a})$ be the action of $\gG$ on $A\underset{N}{_b*_\alpha}\mathcal L(H)$ introduced in \ref{crossed} and $\tilde{\Theta}$ the isomorphism introduced in \ref{tildeTheta} which sends 
$(A\underset{N}{_b*_\alpha}\mathcal L(H))\rtimes_{\underline{\mathfrak a}}\gG$ onto 
$(A\rtimes_\mathfrak a\gG)\underset{N}{_\beta*_\alpha}\mathcal L(H)$; let $\psi$ be a lifted weight on $A$, and $\underline{\psi}$ be the normal semi-finite faithful weight on $A\underset{N}{_b*_\alpha}\mathcal L(H)$ introduced in \ref{psibarre}, and $\widetilde{(\underline{\psi})}$ its dual weight on $(A\underset{N}{_b*_\alpha}\mathcal L(H))\rtimes_{\underline{\mathfrak a}}\gG$; let $\underline{(\tilde{\psi})}$ be the normal semi-finite faithful weight on $(A\rtimes_\mathfrak a\gG)\underset{N}{_\beta*_\alpha}\mathcal L(H)$ introduced by applying \ref{psibarre} to the weight $\tilde{\psi}$ on $A\rtimes_\mathfrak a\gG$. Then :
\newline
(i) for any $X\in \gN_{\widetilde{(\underline{\psi})}}$, $\tilde{\Theta}(X)$ belongs to $\gN_{\underline{(\tilde{\psi})}}$, and :
\[\Lambda_{\underline{(\tilde{\psi})}}(\tilde{\Theta}(X))=(1_H\underset{N}{_\beta\otimes_a}1_{H_\psi}\underset{N}{_b\otimes_\alpha}W^*\sigma_\nu)\Lambda_{\widetilde{(\underline{\psi})}}(X)\]
(ii) we have : $J_{\underline{(\tilde{\psi})}}(1_H\underset{N}{_\beta\otimes_a}1_{H_\psi}\underset{N}{_b\otimes_\alpha}W^*\sigma_\nu)=(1_H\underset{N}{_\beta\otimes_a}1_{H_\psi}\underset{N}{_b\otimes_\alpha}W^*\sigma_\nu)J_{\widetilde{(\underline{\psi})}}$. }

\begin{proof}
The fact that $\tilde{\Theta}(X)$ belongs to $\gN_{\underline{(\tilde{\psi})}}$ is a straightforward corollary of \ref{psitildetheta}(ii). Let us take $x$ in $\gN_\psi$, $\xi$ in $D(_\alpha H, \nu)\cap\mathcal D(\Delta_{\widehat{\Phi}}^{-1/2})$, such that $\Delta_{\widehat{\Phi}}^{-1/2}\zeta$ belongs to $D(_\alpha H, \nu)$, $\eta$ in $D(_\alpha H, \nu)$, and $a$ in $\gN_{\widehat{\Phi}^c}\cap\gN_{\hat{T}^c}$. Then, by \ref{psibarre3}(i), we get that $\rho_\eta^{b, \alpha}x(\rho^{b, \alpha}_\xi)^*$ belongs to $\gN_{\underline{\psi}}$, and, by (\cite{E5} 13.3), $(1\underset{N}{_b\otimes_\alpha}a)\underline{\ga}(\rho_\eta^{b, \alpha}x(\rho^{b, \alpha}_\xi)^*)$ belongs to $\gN_{\widetilde{(\underline{\psi})}}$. Moreover, we have, where $(e_i)_{i\in I}$ is an orthogonal $(\alpha, \nu)$-basis of $H$ :
\begin{align*}
\Lambda_{\underline{(\tilde{\psi})}}(\tilde{\Theta}((1\underset{N}{_b\otimes_\alpha}a)\underline{\ga}(\rho_\eta^{b, \alpha}x(\rho^{b, \alpha}_\xi)^*))
&=
\Lambda_{\underline{(\tilde{\psi})}}(\tilde{\Theta}(1\underset{N}{_b\otimes_\alpha}a)\tilde{\Theta}\underline{\ga}(\rho_\eta^{b, \alpha}x(\rho^{b, \alpha}_\xi)^*))\\
&=
\Lambda_{\underline{(\tilde{\psi})}}(\tilde{\Theta}(1\underset{N}{_b\otimes_\alpha}a)(\ga\underset{N}{_b*_\alpha}id)(\rho_\eta^{b, \alpha}x(\rho^{b, \alpha}_\xi)^*))\\
&=
\Lambda_{\underline{(\tilde{\psi})}}(\tilde{\Theta}(1\underset{N}{_b\otimes_\alpha}a)\rho_\eta^{\beta, \alpha}\ga(x)(\rho^{\beta, \alpha}_\xi)^*)\\
&=
\sum_i \Lambda_{\underline{(\tilde{\psi})}}(\rho_{e_i}^{\beta, \alpha}(\rho_{e_i}^{\beta, \alpha})^*\tilde{\Theta}(1\underset{N}{_b\otimes_\alpha}a)\rho_\eta^{\beta, \alpha}\ga(x)(\rho^{\beta, \alpha}_\xi)^*)
\end{align*}
Then, using \ref{tildeTheta}, we get that $\tilde{\Theta}(1\underset{N}{_b\otimes_\alpha}a)=1\underset{N}{_b\otimes_\alpha}(1\underset{N^o}{_\beta\otimes_\alpha}J_\Phi J_{\widehat{\Phi}})\widehat{\Gamma}^{oc}(a)(1\underset{N^o}{_\beta\otimes_\alpha}J_{\widehat{\Phi}}J_\Phi )$, and, therefore, that :
\[(\rho_{e_i}^{\beta, \alpha})^*\tilde{\Theta}(1\underset{N}{_b\otimes_\alpha}a)\rho_\eta^{\beta, \alpha}
=1\underset{N}{_b\otimes_\alpha}(id\underset{N^o}{_\beta*_{\hat{\alpha}}}\omega_{J_{\widehat{\Phi}}J_\Phi \eta, J_{\widehat{\Phi}}J_\Phi e_i})\widehat{\Gamma}^{oc}(a)\]
and, we get then, applying \ref{psibarre3}(i) to the weight $\tilde{\psi}$, that :
\[\Lambda_{\underline{(\tilde{\psi})}}(\tilde{\Theta}((1\underset{N}{_b\otimes_\alpha}a)\underline{\ga}(\rho_\eta^{b, \alpha}x(\rho^{b, \alpha}_\xi)^*))
=
\sum_i\Lambda_{\underline{(\tilde{\psi})}}(\rho_{e_i}^{\beta, \alpha}
(1\underset{N}{_b\otimes_\alpha}(id\underset{N^o}{_\beta*_{\hat{\alpha}}}\omega_{J_{\widehat{\Phi}}J_\Phi \eta, J_{\widehat{\Phi}}J_\Phi e_i})\widehat{\Gamma}^{oc}(a))
\ga(x)(\rho^{\beta, \alpha}_\xi)^*)\]
is equal to :
\[\sum_i J_{\widehat{\Phi}}\Delta_{\widehat{\Phi}}^{-1/2}\xi\underset{\nu}{_\beta\otimes_{\tilde{a}}}\Lambda_{\tilde{\psi}}((1\underset{N}{_b\otimes_\alpha}(id\underset{N^o}{_\beta*_{\hat{\alpha}}}\omega_{J_{\widehat{\Phi}}J_\Phi \eta, J_{\widehat{\Phi}}J_\Phi e_i})\widehat{\Gamma}^{oc}(a))\ga(x))\underset{\nu}{_\beta\otimes_\alpha}e_i\]
where, for all $n\in N$, we put $\tilde{a}(n)=J_{\tilde{\psi}}(1\underset{N}{_b\otimes_\alpha}\beta(n^*))J_{\tilde{\psi}}$. We then get, by \ref{crossed}, that $\tilde{a}(n)=U^\ga_\psi(1\underset{N^o}{_a\otimes_\beta}\alpha(n))(U^\ga_\psi)^*=a(n)\underset{N}{_b\otimes_\alpha}1$. And, therefore, using now (\cite{E5} 13.3), we get that $\Lambda_{\underline{(\tilde{\psi})}}(\tilde{\Theta}((1\underset{N}{_b\otimes_\alpha}a)\underline{\ga}(\rho_\eta^{b, \alpha}x(\rho^{b, \alpha}_\xi)^*))
$ is equal to :
\[\sum_i J_{\widehat{\Phi}}\Delta_{\widehat{\Phi}}^{-1/2}\xi\underset{\nu}{_\beta\otimes_{a}}\Lambda_\psi(x)\underset{\nu}{_b\otimes_\alpha}\Lambda_{\widehat{\Phi}^c}((id\underset{N^o}{_\beta*_{\hat{\alpha}}}\omega_{J_{\widehat{\Phi}}J_\Phi \eta, J_{\widehat{\Phi}}J_\Phi e_i})\widehat{\Gamma}^{oc}(a))\underset{\nu}{_\beta\otimes_\alpha}e_i\]
which, thanks to \ref{lemW} is equal to :
\[(1_H\underset{N}{_\beta\otimes_a}1_{H_\psi}\underset{N}{_b\otimes_\alpha}W^*\sigma_\nu)
(J_{\widehat{\Phi}}\Delta_{\widehat{\Phi}}^{-1/2}\xi\underset{\nu}{_b\otimes_a}\Lambda_\psi (x)\underset{\nu}{_b\otimes_\alpha}\zeta\underset{\nu}{_{\hat{\beta}}\otimes_\alpha}\Lambda_{\widehat{\Phi}^c}(a))\]
which, using \ref{psibarre3}(i) again, is equal to :
\[(1_H\underset{N}{_\beta\otimes_a}1_{H_\psi}\underset{N}{_b\otimes_\alpha}W^*\sigma_\nu)
(\Lambda_{\underline{\psi}}((\rho_\eta^{b, \alpha}x(\rho^{b, \alpha}_\xi)^*))\underset{\nu}{_{\hat{\beta}}\otimes_\alpha}\Lambda_{\widehat{\Phi}^c}(a))\]
and, by (\cite{E5} 13.3) again, to :
\[(1_H\underset{N}{_\beta\otimes_a}1_{H_\psi}\underset{N}{_b\otimes_\alpha}W^*\sigma_\nu)
\Lambda_{\widetilde{(\underline{\psi})}}((1\underset{N}{_b\otimes_\alpha}a)\underline{\ga}(\rho_\eta^{b, \alpha}x(\rho^{b, \alpha}_\xi)^*))\]
Using now \ref{psibarre3}(i), we get that, for any $a$ in $\gN_{\widehat{\Phi}^c}\cap\gN_{\hat{T}^c}$ and $Y$ in $\gN_{\underline{\psi}}$ :
\[\Lambda_{\underline{(\tilde{\psi})}}(\tilde{\Theta}((1\underset{N}{_b\otimes_\alpha}a)\underline{\ga}(Y))
=(1_H\underset{N}{_\beta\otimes_a}1_{H_\psi}\underset{N}{_b\otimes_\alpha}W^*\sigma_\nu)
\Lambda_{\widetilde{(\underline{\psi})}}((1\underset{N}{_b\otimes_\alpha}a)\underline{\ga}(Y))\]
and, using now (\cite{E5} 13.3), we finish the proof of (i). 
\newline
Let's suppose now that $X$ is analytic with respect to $\widetilde{(\underline{\psi})}$, such that $\sigma_{-i/2}^{\widetilde{(\underline{\psi})}}(X^*)$ belongs to $\gN_{\widetilde{(\underline{\psi})}}$. Then, using \ref{psitildetheta}(ii) and (i), we get that $\tilde{\Theta}(X)$ is analytic with respect to $\underline{(\tilde{\psi})}$, and that $\sigma_{-i/2}^{\underline{(\tilde{\psi})}}(\tilde{\Theta}(X^*))$ belongs to $\gN_{\underline{(\tilde{\psi})}}$. More precisely, we then get :
\begin{align*}
J_{\underline{(\tilde{\psi})}}(1_H\underset{N}{_\beta\otimes_a}1_{H_\psi}\underset{N}{_b\otimes_\alpha}W^*\sigma_\nu)\Lambda_{\widetilde{(\underline{\psi})}}(X)
&=
J_{\underline{(\tilde{\psi})}}\Lambda_{\underline{(\tilde{\psi})}}(\tilde{\Theta}(X))\\
&=
\Lambda_{\underline{(\tilde{\psi})}}(\sigma_{-i/2}^{\underline{(\tilde{\psi})}}(\tilde{\Theta}(X^*)))\\
&=
\Lambda_{\underline{(\tilde{\psi})}}(\tilde{\Theta}(\sigma_{-i/2}^{\widetilde{(\underline{\psi})}}(X^*)))\\
&=
(1_H\underset{N}{_\beta\otimes_a}1_{H_\psi}\underset{N}{_b\otimes_\alpha}W^*\sigma_\nu)\Lambda_{\widetilde{(\underline{\psi})}}(\sigma_{-i/2}^{\widetilde{(\underline{\psi})}}(X^*)))\\
&=
(1_H\underset{N}{_\beta\otimes_a}1_{H_\psi}\underset{N}{_b\otimes_\alpha}W^*\sigma_\nu)J_{\widetilde{(\underline{\psi})}}\Lambda_{\widetilde{(\underline{\psi})}}(X)
\end{align*}
which, by density, gives (ii). \end{proof}

%%%%%%%%%%Upsibarre
\subsection{Proposition}
\label{Upsibarre}
{\it Let $\gG=(N, M, \alpha, \beta, \Gamma, T, T', \nu)$ be a measured quantum groupoid; let $A$ be a von Neumann algebra, $(b, \mathfrak a)$ a weighted action of $\gG$ on $A$, and $(1\underset{N}{_b\otimes_\alpha}\hat{\beta}, \underline{\mathfrak a})$ be the action of $\gG$ on $A\underset{N}{_b*_\alpha}\mathcal L(H)$ introduced in \ref{crossed}; let $\psi$ be a lifted weight on $A$, and $\underline{\psi}$ be the normal semi-finite faithful weight on $A\underset{N}{_b*_\alpha}\mathcal L(H)$ introduced in \ref{psibarre}. Then, the unitary $U^{\underline{\ga}}_{\underline{\psi}}$ satisfies :
\[U^{\underline{\ga}}_{\underline{\psi}}=
(1_H\underset{N}{_\beta\otimes_a}1_{H_\psi}\underset{N}{_b\otimes_\alpha}\sigma W\sigma)(1_H\underset{N}{_\beta\otimes_\alpha}(id\underset{N^o}{_a*_\beta}\varsigma_N)(U_\psi^\ga\underset{N^o}{_a\otimes_\beta}1_H))\sigma_1^{\beta, \alpha}(W^{o*}\underset{N}{_\beta\otimes_a}1_{H_\psi}\underset{N}{_b\otimes_\alpha}1_H)(\sigma_1^{\beta, \hat{\alpha}})^*\]
where $\sigma_1^{\beta, \alpha}$ is the flip from $(H\underset{\nu^o}{_\alpha\otimes_\beta}H)\underset{\nu}{_\beta\otimes_a}H_\psi\underset{\nu}{_b\otimes_\alpha}H$ onto $H\underset{\nu}{_\beta\otimes_\alpha}((H_\psi\underset{\nu}{_b\otimes_\alpha}H)\underset{\nu^o}{_a\otimes_\beta}H)$, and 
$\sigma_1^{\beta, \hat{\alpha}}$ is the flip from $H\underset{\nu}{_\beta\otimes_{\hat{\alpha}}}H\underset{\nu}{_\beta\otimes_a}H_\psi\underset{\nu}{_b\otimes_\alpha}H$ onto $(H\underset{\nu}{_\beta\otimes_a}H_\psi\underset{\nu}{_b\otimes_\alpha}H)\underset{\nu^o}{_{\hat{\alpha}}\otimes_\beta}H$. }

\begin{proof}
Let us recall (\ref{crossed}) that $(1\underset{N}{_b\otimes_\alpha}\hat{\beta}, \underline{\ga})$ is an action of $\gG$ on $A\underset{N}{_b*_\alpha}\mathcal L(H)$. Let $\underline{a}$ be the representation 
of $N$ on $H_{\underline{\psi}}$ defined, for all $n\in N$ by :
\[\underline{a}(n)=J_{\underline{\psi}}\pi_{\underline{\psi}}(1\underset{N}{_b\otimes_\alpha}\hat{\beta}(n^*))J_{\underline{\psi}}\]
Using \ref{psibarre3} (iii) and (ii), we get that :
\[\underline{a}(n)=J_{\widehat{\Phi}}\hat{\beta}(n^*)J_{\widehat{\Phi}}\underset{N}{_\beta\otimes_a}1_{H_\psi}\underset{N}{_b\otimes_\alpha}1_{H}=\hat{\alpha}(n)\underset{N}{_\beta\otimes_a}1_{H_\psi}\underset{N}{_b\otimes_\alpha}1_{H}\]
and, therefore, $U^{\underline{\ga}}_{\underline{\psi}}$ is a unitary from $(H\underset{\nu}{_\beta\otimes_a}H_\psi\underset{\nu}{_b\otimes_\alpha}H)\underset{\nu^o}{_{\hat{\alpha}}\otimes_\beta}H$ onto $H\underset{\nu}{_\beta\otimes_a}H_\psi\underset{\nu}{_b\otimes_\alpha}H\underset{\nu}{_{\hat{\beta}}\otimes_\alpha}H$ given by the formula :
\[U^{\underline{\ga}}_{\underline{\psi}}=J_{\widetilde{(\underline{\psi})}}(J_{\underline{\psi}}\underset{\nu^o}{_{\hat{\alpha}}\otimes_\beta}J_{\widehat{\Phi}})\]
We have, using  \ref{propW*sigma} :
\[(1_H\underset{N}{_\beta\otimes_a}1_{H_\psi}\underset{N}{_b\otimes_\alpha}W^*\sigma_\nu)U^{\underline{\ga}}_{\underline{\psi}}\sigma_1^{\beta, \hat{\alpha}}=
J_{\underline{(\tilde{\psi})}}(1_H\underset{N}{_\beta\otimes_a}1_{H_\psi}\underset{N}{_b\otimes_\alpha}W^*\sigma_\nu)(J_{\underline{\psi}}\underset{\nu^o}{_{\hat{\alpha}}\otimes_\beta}J_{\widehat{\Phi}})\sigma_1^{\beta, \hat{\alpha}}\]
Let $\xi_1$, $\xi_2$ in $D(H_\beta, \nu^o)$, $\xi_3$ in $D(_\alpha H, \nu)$, $\eta\in H_\psi$; then $J_{\widehat{\Phi}}\xi_1$ belongs to $D(_\alpha H, \nu)$, and let us define $\zeta_i\in D(H_\beta, \nu^o)$ and $\zeta'_i\in D(_\alpha H, \nu)$ such that :
\[W^*(J_{\widehat{\Phi}}\xi_1\underset{\nu^o}{_\alpha\otimes_{\hat{\beta}}}J_{\widehat{\Phi}}\xi_2)=lim_J\sum_{i\in J}(\zeta_i\underset{\nu}{_\beta\otimes_\alpha}\zeta'_i)\]
the limit being taken on the filter of finite subsets $J\subset I$. 
Let us look at the image of the vector $\xi_1\underset{\nu}{_\beta\otimes_{\hat{\alpha}}}\xi_2\underset{\nu}{_\beta\otimes_a}\eta\underset{\nu}{_b\otimes_\alpha}\xi_3$ under the unitary :
\[J_{\underline{(\tilde{\psi})}}(1_H\underset{N}{_\beta\otimes_a}1_{H_\psi}\underset{N}{_b\otimes_\alpha}W^*\sigma_\nu)(J_{\underline{\psi}}\underset{\nu^o}{_{\hat{\alpha}}\otimes_\beta}J_{\widehat{\Phi}})\sigma_1^{\beta, \hat{\alpha}}\] 
This vector is first sent by $\sigma_1^{\beta, \hat{\alpha}}$ on $(\xi_2\underset{\nu}{_\beta\otimes_a}\eta\underset{\nu}{_b\otimes_\alpha}\xi_3)\underset{\nu^o}{_{\hat{\alpha}}\otimes_\beta}\xi_1$, then $(J_{\underline{\psi}}\underset{\nu^o}{_{\hat{\alpha}}\otimes_\beta}J_{\widehat{\Phi}})$ sends it on $J_{\widehat{\Phi}}\xi_3\underset{\nu}{_\beta\otimes_a}J_\psi\eta\underset{\nu}{_b\otimes_\alpha}J_{\widehat{\Phi}}\xi_2\underset{\nu}{_{\hat{\beta}}\otimes_\alpha}J_{\widehat{\Phi}}\xi_1$, then $(1_H\underset{N}{_\beta\otimes_a}1_{H_\psi}\underset{N}{_b\otimes_\alpha}W^*\sigma_\nu)$ sends it on 
\[J_{\widehat{\Phi}}\xi_3\underset{\nu}{_\beta\otimes_a}J_\psi\eta\underset{\nu}{_b\otimes_\alpha}W^*(J_{\widehat{\Phi}}\xi_1\underset{\nu^o}{_\alpha\otimes_{\hat{\beta}}}J_{\widehat{\Phi}}\xi_2)=
lim_J\sum_{i\in J}(J_{\widehat{\Phi}}\xi_3\underset{\nu}{_\beta\otimes_a}J_\psi\eta\underset{\nu}{_b\otimes_\alpha}
\zeta_i\underset{\nu}{_\beta\otimes_\alpha}\zeta'_i)\]
and $J_{\underline{(\tilde{\psi})}}$ sends it then on :
\[lim_J\sum_{i\in J}(J_{\widehat{\Phi}}\zeta'_i\underset{\nu}{_\beta\otimes_a}J_{\tilde{\psi}}(J_\psi\eta\underset{\nu}{_b\otimes_\alpha}\zeta_i)\underset{\nu}{_\beta\otimes_\alpha}\xi_3)
=
lim_J\sum_{i\in J}(J_{\widehat{\Phi}}\zeta'_i\underset{\nu}{_\beta\otimes_a}U^\ga_\psi(\eta\underset{\nu^o}{_a\otimes_\beta}J_{\widehat{\Phi}}\zeta_i)\underset{\nu}{_\beta\otimes_\alpha}\xi_3)\]
which is equal to :
\[(1_H\underset{N}{_\beta\otimes_a}U^\ga_\psi\underset{N}{_b\otimes_\alpha}1_H)(1_H\underset{N}{_\beta\otimes_\alpha}\sigma_\nu\underset{N}{_b\otimes_\alpha}1_H)(\sigma W^{o*}\underset{N}{_\beta\otimes_a}1_{H_\psi}\underset{N}{_b\otimes_\alpha}1_H)(\xi_1\underset{\nu}{_\beta\otimes_{\hat{\alpha}}}\xi_2\underset{\nu}{_\beta\otimes_a}\eta\underset{\nu}{_b\otimes_\alpha}\xi_3)\]
from which, using again the density of finite sums of elementary tensors in the relative Hilbert tensor product, we get that :
\begin{multline*}
(1_H\underset{N}{_\beta\otimes_a}1_{H_\psi}\underset{N}{_b\otimes_\alpha}W^*\sigma_\nu)U^{\underline{\ga}}_{\underline{\psi}}\sigma_1^{\beta, \hat{\alpha}}=\\
(1_H\underset{N}{_\beta\otimes_a}U^\ga_\psi\underset{N}{_b\otimes_\alpha}1_H)(1_H\underset{N}{_\beta\otimes_\alpha}\sigma_\nu\underset{N}{_b\otimes_\alpha}1_H)(\sigma W^{o*}\underset{N}{_\beta\otimes_a}1_{H_\psi}\underset{N}{_b\otimes_\alpha}1_H)=\\
(1_H\underset{N}{_\beta\otimes_a}(id\underset{N^o}{_a*_\beta}\varsigma_N)(U_\psi^\ga\underset{N}{_b\otimes_\alpha}1_H))\sigma_1^{\beta, \alpha}(W^{o*}\underset{N}{_\beta\otimes_a}1_{H_\psi}\underset{N}{_b\otimes_\alpha}1_H)
\end{multline*}
from which we get the result. \end{proof}

%%%%Upsicorep
\subsection{Proposition}
\label{Upsicorep}
{\it Let $\gG=(N, M, \alpha, \beta, \Gamma, T, T', \nu)$ be a measured quantum groupoid; let $A$ be a von Neumann algebra, and let $(b, \mathfrak a)$ be a weighted action of $\gG$ on $A$; let $\psi$ be a lifted weight on $A$; 
then, the unitary $U^\ga_\psi$ introduced in \ref{crossed} is a copresentation of $\gG$.}

\begin{proof}
With the notations of \ref{Upsibarre}, we get, using \ref{Upsibarre}, that : 
\[1_H\underset{N}{_\beta\otimes_a}(id\underset{N^o}{_a*_\beta}\varsigma_N)(U_\psi^\ga\underset{N}{_b\otimes_\alpha}1_H)=
(1_H\underset{N}{_\beta\otimes_a}1_{H_\psi}\underset{N}{_b\otimes_\alpha}\sigma W^*\sigma)U^{\underline{\ga}}_{\underline{\psi}}\sigma_1^{\beta, \hat{\alpha}}(W^{o}\underset{N}{_\beta\otimes_a}1_{H_\psi}\underset{N}{_b\otimes_\alpha}1_H)(\sigma_1^{\beta, \alpha})^*\]
which we shall write, for simplification, with the usual leg numbering notation :
\[(U^\ga_{\psi})_{2,4}=\widehat{W}_{3,4}U^{\underline{\ga}}_{\underline{\psi}}(W^o)_{4, 1}\]
But $\widehat{W}$ is a corepresentation of $\gG^o$ (\cite{E5}, 5.6), $U^{\underline{\ga}}_{\underline{\psi}}$ is a corepresentation of $\gG$ by \ref{cora}, and $\sigma W^o\sigma$ is a corepresentation of $\gG^o$ by (\cite{E5}, 5.6 and 5.3). So, we get :
\begin{align*}
(id*\Gamma)(U^\ga_\psi)_{2,4,5}
&=
\widehat{W}_{3,5}\widehat{W}_{3,4}(U^{\underline{\ga}}_{\underline{\psi}})_{1,2,3,4}(U^{\underline{\ga}}_{\underline{\psi}})_{1,2,3,5}W^o_{5,1}W^o_{4,1}\\
&=\widehat{W}_{3,5}(U^\ga_{\psi})_{2,4}W^{o*}_{4,1}(U^{\underline{\ga}}_{\underline{\psi}})_{1,2,3,5}W^o_{5,1}W^o_{4,1}\\
&=(U^\ga_{\psi})_{2,4}W^{o*}_{4,1}\widehat{W}_{3,5}(U^{\underline{\ga}}_{\underline{\psi}})_{1,2,3,5}
W^o_{5,1}W^o_{4,1}\\
&=(U^\ga_{\psi})_{2,4}W^{o*}_{4,1}(U^\ga_{\psi})_{2,5}W^o_{4,1}\\
&=(U^\ga_{\psi})_{2,4}(U^\ga_{\psi})_{2,5}
\end{align*}
which shows that $U^\ga_\psi$ is a corepresentation. A more complete proof is a painful exercise we leave to the reader. \end{proof}

%%%standard
\subsection{Theorem}
\label{standard}
{\it Let $\gG=(N, M, \alpha, \beta, \Gamma, T, T', \nu)$ be a measured quantum groupoid; let $A$ be a von Neumann algebra, and let $(b, \mathfrak a)$ be a weighted action of $\gG$ on $A$; then, for any normal semi-finite faithful weight $\psi$ on $A$, the unitary $U^\ga_\psi$ introduced in \ref{crossed} is a standard implementation of $\ga$, in the sense of \ref{action}. }

\begin{proof}
As the action is weighted, there exists a normal semi-finite faithful weight $\psi$ on $A$ which is lifted from $\nu^o$; we get that $U^\ga_\psi$ is a corepresentation by \ref{Upsicorep}, and is therefore a standard implementation. Using now \ref{propu}, we easily get that it remains true for any normal semi-finite faithful weight on $A$, which is the result. \end{proof}

%%%remstandard
\subsection{Remark}
\label{remstandard}
In \ref{standard}, we had obtained that $U^\ga_\psi$ is a standard implementation of $\ga$, if there exists a normal-semi-finite faithful operator-valued weight from $A$ onto $b(N)$; this is true in particular in the following cases :
\newline
(i) $\gG$ is a locally compact quantum group ($N=\mathbb{C}$); this result was obtained in (\cite{V1} 4.4);
\newline
(ii) if $N$ is abelian and $b(N)\subset Z(A)$; in particular, if $\gG$ is a measured groupoid; we shall discuss this particular case in \ref{exgd}. More general, if $\gG$ is a continuous field of locally compact quantum groups (\ref{ex} (iv)), or is De Commer's example (\ref{ex} (v)). 
\newline
(iii) $A$ is a type I factor; if we write $A=\mathcal L(\gH)$, starting from any normal semi-finite weight on $b(N)'$, we get a normal faithful semi-finite operator-valued weight from $A$ to $b(N)$. More generally, this remains true if $A$ is a sum of type I factors;
\newline
(iv) $N$ is a sum of type I factors (in particular, if $N$ is a finite dimensional algebra, which is the case, in particular if $\gG$ is a finite dimensional quantum groupoid);
\newline
(v) $N$ and $A$ are semi-finite. 
\newline
In \ref{Uinv}, the result was proved if $\ga$ is a dual action. 

%%%exgd
\subsection{Example}
\label{exgd}
Let $\mathcal G$ be a measured groupoid, with $\mathcal G^{(0)}$ as its set of units, $r$ and $s$ its range and source application, $(\lambda^u)_{u\in \mathcal G^{(0)}}$ its Haar system, and $\nu$ a quasi-invariant measure; let $\mu=\int_{\mathcal G^{(0)}}\lambda^ud\nu$ ; let us consider the von Neumann algebra $L^\infty(\mathcal G, \mu)$, which is a $L^\infty(\mathcal G^{(0)})$-bimodule, thanks to the two homomorphisms $r_{\mathcal G}$ and $s_{\mathcal G}$ defined, for $f$ in $L^\infty(\mathcal G^{(0)})$ by $r_{\mathcal G}(f)=f\circ r$ and $s_{\mathcal G}(f)=f\circ s$. We have shown in (\cite{E5}, 3.1, 3.4 and 3.17) how it is possible to put a measured quantum groupoid structure on this von Neumann bimodule. 
\newline
An action $(b, \ga)$ of this measured quantum groupoid on a von Neumann algebra $A$ verifies that $b(L^\infty(\mathcal G^{(0)}))\subset Z(A)$, and, therefore, $A$ can be decomposed as $A=\int_{\mathcal G^{(0)}}^\oplus A^xd\nu(x)$ (\cite{E5}, 6.1); moreover, let $\psi$ be a normal semi-finite faithful on $A=\int_{\mathcal G^{(0)}}^\oplus A^xd\nu(x)$. Then $\psi$ is a lifted weight; more precisely, there exists a measurable field $\psi^x$ of normal semi-finite faithful weights, such that $\psi=\int_{\mathcal G^{(0)}}^\oplus \psi^x d\nu(x)$ in the sense of (\cite{T} 4.6), and $H_\psi=\int _{\mathcal G^{(0)}}^\oplus H_{\psi^x}d\nu (x)$. 
\newline
On the other hand, the action $\ga$ is (\cite{E5}, 6.3) an action of $\mathcal G$ in the sense of (\cite{Y3}, 3.1), i.e. for all $g\in\mathcal G$, there exists a family of $*$-isomorphisms $\ga_g$ from $A^{s(g)}$ onto $A^{r(g)}$, such that, if $(g_1, g_2)\in\mathcal G^{(2)}$, we have $\ga_{g_1g_2}=\ga_{g_1}\ga_{g_2}$, and such that, for any normal positive functional $\omega=\int_{\mathcal G^{(0)}}^\oplus \omega^x d\nu(x)$, and any $y=\int_{\mathcal G^{(0)}}^\oplus y^x d\nu (x)$, the function $g\mapsto \omega^{r(g)}(\ga_g(y^{s(g)}))$ is $\mu$-measurable. These $*$-isomorphisms have standard implementations $u_g : H_{\psi^{s(g)}}\to H_{\psi^{r(g)}}$ such that $\ga_g(y^{s(g)})=u_gy^{s(g)}u_g^*$. if $(g_1, g_2)\in\mathcal G^{(2)}$, we have $u_{g_1g_2}=u_{g_1}u_{g_2}$.
\newline
More precisely, the Hilbert space $H_\psi\underset{\nu}{_b\otimes_{r_\mathcal G}}L^2(\mathcal G, \mu)$ can be identified with $\int_{\mathcal G}^\oplus H_{\psi^{r(g)}}d\mu(g)$. We then get :
\[\ga (\int_{\mathcal G^{(0)}}^\oplus y^x d\nu (x))=\int_{\mathcal G}^\oplus \ga_g(y^{s(g)})d\mu(g)\]
In \cite{Y1} and \cite{Y2} is given a construction of the crossed product of $A$ by $\mathcal G$; using (\cite{Y3} 2.14), we see (\cite{E5}, 9.2) that this crossed-product is isomorphic to the definition given in (\cite{E5}, 9.1). Moreover, we get the same notion of dual action (\cite{E5}, 9.6) and of dual weight (\cite{E5}, 13.1). 
\newline
 As $b$ is central, we have $a=b$, and the Hilbert space $H_\psi\underset{\nu^o}{_a\otimes_{s_\mathcal G}}L^2(\mathcal G, \mu)$ can be identified with $\int_{\mathcal G}^\oplus H_{\psi^{s(g)}}d\mu(g)$.  Using then \cite{Y3}, 2.6, we get that $U^\ga_\psi=\int_\mathcal G^\oplus u_gd\mu(g)$, which is a unitary from $\int_{\mathcal G}^\oplus H_{\psi^{s(g)}}d\mu(g)$ onto $\int_{\mathcal G}^\oplus H_{\psi^{r(g)}}d\mu(g)$. 
%%%%gamma
\section{The $(b,\gamma)$ property for weights}
\label{gamma}
If $\gG=(N, M, \alpha, \beta, \Gamma, T, T', \nu)$ be a measured quantum groupoid, and if $b$ is a normal faithful non degenerate anti-homomorphism from $N$ into a von Neumann algebra $A$, we define the $(b, \gamma)$ property for normal faithful semi-finite weights on $A$ (\ref{defgamma}).  We define then, for such a weight, a normal semi-finite faithful weight $\underline{\psi_\delta}$ on $A\underset{N}{_b*_\alpha}\mathcal L(H)$ (\ref{thgamma}). We obtain then several technical results (\ref{cor2}, \ref{deltaDelta}, \ref{psibarre4}) which will be used in chapter \ref{bidualw}. 

%%%defgamma
\subsection{Definition}
\label{defgamma}
Let $\gG=(N, M, \alpha, \beta, \Gamma, T, T', \nu)$ be a measured quantum groupoid, and let $b$ be a normal faithful non degenerate anti-homomorphism from $N$ into a von Neumann algebra $A$; we shall say that a normal faithful semi-finite weight $\psi$ on $A$ satisfies the $(b,\gamma)$ property if, for all $n\in N$ and $t\in\mathbb{R}$, we have $\sigma_t^\psi (b(n))=b(\gamma_t(n))$, where $\gamma_t$ is the one-parameter automorphism group of $N$ defined by $\sigma_t^T(\beta(n))=\beta(\gamma_t(n))$ (\cite{E5}, 3.8 (v)). 

%%%delta
\subsection{Example}
\label{delta}
Let $\gG=(N, M, \alpha, \beta, \Gamma, T, T', \nu)$ be a measured quantum groupoid, $A$ a von Neumann algebra, $(b,\ga)$ an action of $\gG$ on $A$, $\psi$ a $\delta$-invariant normal faithful semi-finite weight on $A$ bearing the density property, as defined in (\cite{E5}) and recalled in \ref{action}. Then, $\psi$ satisfies the $(b, \gamma)$ property. 
\newline
Namely, for any $x\in A$, $t\in\mathbb{R}$, we have :
\[\mathfrak a(\sigma_t^\psi(x))=(\Delta_\psi^{it}\underset{\nu}{_b\otimes_\alpha}\delta^{-it}\Delta_{\widehat{\Phi}}^{-it})\mathfrak a(x)(\Delta_\psi^{-it}\underset{\nu}{_b\otimes_\alpha}\delta^{it}\Delta_{\widehat{\Phi}}^{it})\]
and, therefore, for any $n\in N$, we get, using (\cite{E5}, 3.8(ii)) :
\begin{align*}
\ga(\sigma_t^\psi(b(n)))
&=
(\Delta_\psi^{it}\underset{\nu}{_b\otimes_\alpha}\delta^{-it}\Delta_{\widehat{\Phi}}^{-it})(1\underset{N}{_b\otimes_\alpha}\beta(n))(\Delta_\psi^{-it}\underset{\nu}{_b\otimes_\alpha}\delta^{it}\Delta_{\widehat{\Phi}}^{it})\\
&=1\underset{N}{_b\otimes_\alpha}\delta^{-it}\Delta_{\widehat{\Phi}}^{-it}\beta(n)\delta^{it}\Delta_{\widehat{\Phi}}^{it}\\
&=1\underset{N}{_b\otimes_\alpha}\sigma_{t}^\Phi\sigma_{-t}^{\Phi\circ R}\tau_{-t}(\beta(n))\\
&=1\underset{N}{_b\otimes_\alpha}\sigma_{t}^\Phi(\beta(n))\\
&=1\underset{N}{_b\otimes_\alpha}\beta(\gamma_t(n))\\
&=\ga(b(\gamma_t(n)))
\end{align*}
from which we get the property, by the injectivity of $\ga$. 
%%%%%betahat
\subsection{Example}
\label{betahat}
Let $\gG=(N, M, \alpha, \beta, \Gamma, T, T', \nu)$ be a measured quantum groupoid, and let $(N, b, A)$ be a faithful weighted right von Neumann right-module, in the sense of \ref{defw}; let $\psi$ be a normal faithful semi-finite weight on $A$, lifted from $\nu^o$, and let $\underline{\psi}$ be the normal faithful semi-finite weight on $A\underset{N}{_b*_\alpha}\mathcal L(H)$ defined in \ref{psibarre}. Then, $\underline{\psi}$ satisfies the $(1\underset{N}{_b\otimes_\alpha}\hat{\beta}, \gamma)$ property. 
\newline
Namely, using \ref{psibarre2} and (\cite{E5} 3.10 (vii)), we get :
\[\sigma_t^{\underline{\psi}}(1\underset{N}{_b\otimes_\alpha}\hat{\beta}(n))=1\underset{N}{_b\otimes_\alpha}\Delta_{\widehat{\Phi}}^{-it}\hat{\beta}(n)\Delta_{\widehat{\Phi}}^{it}=1\underset{N}{_b\otimes_\alpha}\sigma_{-t}^{\widehat{\Phi}}(\hat{\beta}(n))=
1\underset{N}{_b\otimes_\alpha}\hat{\beta}(\hat{\gamma}_{-t}(n))=1\underset{N}{_b\otimes_\alpha}\hat{\beta}(\gamma_t(n))\]

%%%%%thgamma
\subsection{Theorem}
\label{thgamma}
{\it Let $\gG=(N, M, \alpha, \beta, \Gamma, T, T', \nu)$ be a measured quantum groupoid, and let $b$ be a normal faithful non degenerate anti-homomorphism from $N$ into a von Neumann algebra $A$; let $\psi$ be a normal faithful semi-finite weight on $A$ satisfying the $(b,\gamma)$ property; then :
\newline
(i) it is possible to define a one-parameter group of unitaries $\Delta_\psi^{it}\underset{\nu}{_b\otimes_\alpha}(\delta\Delta_{\widehat{\Phi}})^{-it}$ on $H_\psi\underset{\nu}{_b\otimes_\alpha}H$, with natural values on elementary tensors. We shall denote $\Delta_\psi^{1/2}\underset{\nu}{_b\otimes_\alpha}(\delta\Delta_{\widehat{\Phi}})^{-1/2}$ its analytic generator. 
\newline
(ii) there exists a normal semi-finite faithful weight $\underline{\psi}_\delta$ on $A\underset{N}{_b*_\alpha}\mathcal L(H)$ such that :
\[\frac{d\underline{\psi}_\delta}{d\psi^o}=\Delta_\psi^{1/2}\underset{\nu}{_b\otimes_\alpha}(\delta\Delta_{\widehat{\Phi}})^{-1/2}\]
(iii) for any $a$ in $\gN_\psi\cap\gN_\psi^*$, and $\xi\in D(_\alpha H, \nu)\cap \mathcal D((\delta\Delta_{\widehat{\Phi}})^{-1/2})$, such that $(\delta\Delta_{\widehat{\Phi}})^{-1/2}\xi$ belongs to $D(_\alpha H, \nu)$, we have :}
\[\underline{\psi}_\delta(\rho_\xi^{b, \alpha}aa^*(\rho_\xi^{b, \alpha})^*)=\|\Delta_\psi^{1/2}\Lambda_\psi (a)\underset{\nu}{_b\otimes_\alpha}(\delta\Delta_{\widehat{\Phi}})^{-1/2}\xi\|^2\]

\begin{proof}
Let $\eta\in D(_\alpha H, \nu)$, $n\in\gN_\nu$; then, we get :
\begin{align*}
\alpha(n)(\delta\Delta_{\widehat{\Phi}})^{-it}\eta
&=(\delta\Delta_{\widehat{\Phi}})^{-it}\sigma_t^{\widehat{\Phi}}\sigma_t^{\Phi\circ R}\sigma_{-t}^\Phi(\alpha(n))\eta\\
&=(\delta\Delta_{\widehat{\Phi}})^{-it}\alpha(\sigma^\nu_t\gamma_{-t}\tau_{-t}(n))\eta\\
&=(\delta\Delta_{\widehat{\Phi}})^{-it}R^{\alpha, \nu}(\eta)\Lambda_\nu(\gamma_{-t}(n))
\end{align*}
There exists a positive self-adjoint non singular operator $h$ on $H_\nu$ such that :
\[\Lambda_\nu(\gamma_t(n))=h^{it}\Lambda_\nu(n)\]
We then get that :
\[\alpha(n)(\delta\Delta_{\widehat{\Phi}})^{-it}\eta=(\delta\Delta_{\widehat{\Phi}})^{-it}R^{\alpha, \nu}(\eta)h^{-it}\Lambda_\nu(n)\]
from which we get that $(\delta\Delta_{\widehat{\Phi}})^{-it}\eta$ belongs to $D(_\alpha H, \nu)$, and that :
\[R^{\alpha, \nu}((\delta\Delta_{\widehat{\Phi}})^{-it}\eta)=(\delta\Delta_{\widehat{\Phi}})^{-it}R^{\alpha, \nu}(\eta)h^{-it}\]
from which we get that :
\[<(\delta\Delta_{\widehat{\Phi}})^{-it}\eta, (\delta\Delta_{\widehat{\Phi}})^{-it}\eta>_{\alpha, \nu}=
h^{it}<\eta, \eta>_{\alpha, \nu}h^{-it}\]
As we have, for all $m\in N$, $\gamma_t(m)=h^{it}mh^{-it}$, we therefore get that :
\[<(\delta\Delta_{\widehat{\Phi}})^{-it}\eta, (\delta\Delta_{\widehat{\Phi}})^{-it}\eta>_{\alpha, \nu}^o=
\gamma_t(<\eta, \eta>_{\alpha, \nu}^o)\]
and, therefore, for all $\xi\in H_\psi$ :
\begin{align*}
\|\Delta_\psi^{it}\xi\underset{\nu}{_b\otimes_\alpha}(\delta\Delta_{\widehat{\Phi}})^{-it}\eta\|^2
&=
(b(\gamma_t(<\eta, \eta>_{\alpha, \nu}^o))\Delta_\psi^{it}\xi|\Delta_\psi^{it}\xi)\\
&=
(\sigma_t^\psi(b(<\eta, \eta>_{\alpha, \nu}^o))\Delta_\psi^{it}\xi|\Delta_\psi^{it}\xi)\\
&=
\|\xi\underset{\nu}{_b\otimes_\alpha}\eta\|^2
\end{align*}
which is (i). 
\newline
As $(\Delta_\psi^{it}\underset{\nu}{_b\otimes_\alpha}(\delta\Delta_{\widehat{\Phi}})^{-it})(J_\psi xJ_\psi\underset{N}{_b\otimes_\alpha}1)(\Delta_\psi^{it}\underset{\nu}{_b\otimes_\alpha}(\delta\Delta_{\widehat{\Phi}})^{-it})=J_\psi \sigma_t^{\psi}(x)J_\psi\underset{N}{_b\otimes_\alpha}1$, we get (ii). Result (iii) is just a corollary of (ii) and \ref{subalgebra}(iv). \end{proof}

%%%corthgamma
\subsection{Corollary}
\label{corthgamma}
{\it Let $\gG$ be a measured quantum groupoid, and $(b,\mathfrak a)$ an action of $\gG$ on a von Neumann algebra $A$; let $\psi$ be a $\delta$-invariant weight on $A$, bearing the density condition, as defined in \ref{action}, and $\overline{\psi_\ga}$ the weight constructed on $A\underset{N}{_b*_\alpha}\mathcal L(H)$ by transporting the bidual weight (\ref{crossed}) of $\psi$. Using \ref{delta}, we can use \ref{thgamma} and define the weight $\underline{\psi}_\delta$ on $A\underset{N}{_b*_\alpha}\mathcal L(H)$
Then, we have : $\overline{\psi_\ga}=\underline{\psi}_\delta$. }

\begin{proof}
We have, in general, $\frac{d\overline{\psi_\ga}}{d\psi^o}=\Delta_{\tilde{\psi}}^{1/2}$ (\ref{crossed}). 
So, using \ref{Uinv}(ii) and \ref{thgamma}, we get the result. \end{proof}

%%%cor2
\subsection{Corollary}
\label{cor2}
{\it Let $\gG=(N, M, \alpha, \beta, \Gamma, T, T', \nu)$ be a measured quantum groupoid, and let $b$ be a normal faithful non degenerate anti-homomorphism from $N$ into a von Neumann algebra $A$; let 
$\psi_1$ (resp. $\psi_2$) be a normal faithful semi-finite weight on $A$ satisfying the $(b,\gamma)$ property; then :
\newline
(i) the cocycle $(D\psi_1:D\psi_2)_t$ belongs to $A\cap b(N)'$; 
\newline
(ii) we have : $(D\underline{\psi_1}_\delta : D\underline{\psi_2}_\delta)_t=(D\psi_1:D\psi_2)_t\underset{N}{_b\otimes_\alpha}1$. }

\begin{proof}
For any $x\in A$, we have :
\[\sigma_t^{\psi_1}(x)=(D\psi_1:D\psi_2)_t\sigma_t^{\psi_2}(x)(D\psi_1:D\psi_2)_t^*\]
and, therefore :
\[\sigma_t^{\psi_1}\circ\sigma_{-t}^{\psi_2}(x)=(D\psi_1:D\psi_2)_t x(D\psi_1:D\psi_2)_t^*\]
In particular, we get, for any $n\in N$ :
\[b(n)=(D\psi_1:D\psi_2)_t b(n)(D\psi_1:D\psi_2)_t^*\]
from which we get (i). 
Let $(\gH, \pi, J, \mathcal P)$ be a standard representation of the von Neumann algebra $A$; then $A^o$ is represented on $\gH$ by $JAJ$; for any normal semi-finite faithful weight $\psi$ on $A$, we have $\frac{d\psi}{d\psi^o}=\Delta_\psi^{1/2}$; moreover, we have then :
\begin{align*}
(\frac{d\psi_1}{d\psi_1^o})^{it}(D\psi_1^o:D\psi_2^o)_t(\frac{d\psi_2^o}{d\psi_2})^{it}
&=
(\frac{d\psi_1}{d\psi_1^o})^{it}(\frac{d\psi_1^o}{d\psi_1})^{it}(\frac{d\psi_2^o}{d\psi_1})^{-it}(\frac{d\psi_2^o}{d\psi_2})^{it}\\
&=
(\frac{d\psi_1}{d\psi_2^o})^{it}(\frac{d\psi_2}{d\psi_2^o})^{-it}\\
&=
(D\psi_1: D\psi_2)_{t}
\end{align*}
and, therefore $(D\psi_1^o:D\psi_2^o)_t=\Delta_{\psi_1}^{-it}(D\psi_1: D\psi_2)_{t}\Delta_{\psi_2}^{it}$. 
By similar arguments, we have on $\gH\underset{\nu}{_b\otimes_\alpha}H$ :
\begin{align*}
(D\underline{\psi_1}_\delta : D\underline{\psi_2}_\delta)_t
&=
(\frac{d\underline{\psi_1}_\delta}{d\psi_1^o})^{it}(\frac{d\psi_1^o}{d\underline{\psi_2}_\delta})^{it}\\
&=(\frac{d\underline{\psi_1}_\delta}{d\psi_1^o})^{it}(D\psi_1^o:D\psi_2^o)_t(\frac{d\underline{\psi_2}_\delta}{d\psi_2^o})^{-it}
\end{align*}
As $(D\psi_1^o:D\psi_2^o)_t$ belongs to $JAJ\underset{N}{_b\otimes_\alpha}1_H$ and is therefore equal to :
\[\Delta_{\psi_1}^{-it}(D\psi_1:D\psi_2)_{t}\Delta_{\psi_2}^{it}\underset{N}{_b\otimes_\alpha}1_H\]
we obtain, using \ref{thgamma}(ii), that $(D\underline{\psi_1}_\delta : D\underline{\psi_2}_\delta)_t$ is equal to :
\[(\Delta_{\psi_1}^{it}\underset{\nu}{_b\otimes_\alpha}(\delta\Delta_{\widehat{\Phi}})^{-it})(\Delta_{\psi_1}^{-it}(D\psi_1:D\psi_2)_{t}\Delta_{\psi_2}^{it}\underset{N}{_b\otimes_\alpha}1_H)(\Delta_{\psi_2}^{-it}\underset{\nu}{_b\otimes_\alpha}(\delta\Delta_{\widehat{\Phi}})^{it})\]
from which we get the result. \end{proof}

%%%%deltaDelta
\subsection{Proposition}
\label{deltaDelta}
{\it Let $\gG=(N, M, \alpha, \beta, \Gamma, T, T', \nu)$ be a measured quantum groupoid; it is possible to define one parameter groups of unitaries $\Delta_{\widehat{\Phi}}^{it}\underset{\nu}{_\beta\otimes_\alpha}\Delta_{\widehat{\Phi}}^{it}$ and $(\delta\Delta_{\widehat{\Phi}})^{it}\underset{\nu^o}{_\alpha\otimes_{\hat{\beta}}}\Delta_{\widehat{\Phi}}^{it}$, with natural values on elementary tensors, and 
we have :}
\[W(\Delta_{\widehat{\Phi}}^{it}\underset{\nu}{_\beta\otimes_\alpha}\Delta_{\widehat{\Phi}}^{it})W^*=
(\delta\Delta_{\widehat{\Phi}})^{it}\underset{\nu^o}{_\alpha\otimes_{\hat{\beta}}}\Delta_{\widehat{\Phi}}^{it}\]

\begin{proof}
From (\cite{E5} 3.10 (vi)), we get that $\Delta_{\widehat{\Phi}}$ is the closure of $PJ_\Phi\delta^{-1}J_\Phi$, where $P$ is the managing operator of the pseudo-multiplicative unitary $W$, and $\delta$ the modulus of $\gG$; in (\cite{E5} 3.8 (vii)), we had got that it is possible to define one parameter groups of unitaries $P^{it}\underset{\nu}{_\beta\otimes_\alpha}P^{it}$ and $P^{it}\underset{\nu^o}{_\alpha\otimes_{\hat{\beta}}}P^{it}$, with natural values on elementary tensors, and that :
\[W(P^{it}\underset{\nu}{_\beta\otimes_\alpha}P^{it})=(P^{it}\underset{\nu^o}{_\alpha\otimes_{\hat{\beta}}}P^{it})W\]
On the other hand, it is possible (\cite{E5}, 3.8 (vi)) to define a one parameter group of unitaries $\delta^{it}\underset{\nu}{_\beta\otimes_\alpha}\delta^{it}$, with natural values on elementary tensors, and that :
\[\delta^{it}\underset{\nu}{_\beta\otimes_\alpha}\delta^{it}=\Gamma(\delta^{it})=W^*(1\underset{N^o}{_\alpha\otimes_{\hat{\beta}}}\delta^{it})W\]
Moreover, we know, from (\cite{E5}, 3.11 (iii)), that :
\[W(J_{\widehat{\Phi}}\underset{\nu^o}{_\alpha\otimes_{\hat{\beta}}}J_\Phi)=(J_{\widehat{\Phi}}\underset{\nu^o}{_\alpha\otimes_{\hat{\beta}}}J_\Phi)W^*\]
and from (\cite{E5} 3.8 (vi)) that $J_{\widehat{\Phi}}\delta^{-it}J_{\widehat{\Phi}}=R(\delta^{it})=\delta^{-it}$. 
\newline
With all these data, we get that it is possible to define $\Delta_{\widehat{\Phi}}^{it}\underset{\nu}{_\beta\otimes_\alpha}\Delta_{\widehat{\Phi}}^{it}$ as :
\[\Delta_{\widehat{\Phi}}^{it}\underset{\nu}{_\beta\otimes_\alpha}\Delta_{\widehat{\Phi}}^{it}=(P^{it}\underset{\nu}{_\beta\otimes_\alpha}P^{it})(J_\Phi\delta^{it}J_\Phi\underset{N}{_\beta\otimes_\alpha}J_\Phi\delta^{it}J_\Phi)\]
and $(\delta\Delta_{\widehat{\Phi}})^{it}\underset{\nu^o}{_\alpha\otimes_{\hat{\beta}}}\Delta_{\widehat{\Phi}}^{it}$ as :
\[(\delta\Delta_{\widehat{\Phi}})^{it}\underset{\nu^o}{_\alpha\otimes_{\hat{\beta}}}\Delta_{\widehat{\Phi}}^{it}=(P^{it}\underset{\nu^o}{_\alpha\otimes_{\hat{\beta}}}P^{it})(J_{\widehat{\Phi}}\underset{\nu}{_\beta\otimes_\alpha}J_\Phi)(\delta^{it}\underset{\nu}{_\beta\otimes_\alpha}\delta^{it})(J_{\widehat{\Phi}}\underset{\nu^o}{_\alpha\otimes_{\hat{\beta}}}J_\Phi)(J_\Phi\delta^{it}J_\Phi\underset{N^o}{_\alpha\otimes_{\hat{\beta}}}1)\]
and to verify that :
\begin{align*}
W(\Delta_{\widehat{\Phi}}^{it}\underset{\nu}{_\beta\otimes_\alpha}\Delta_{\widehat{\Phi}}^{it})W^*
&=
W(P^{it}\underset{\nu}{_\beta\otimes_\alpha}P^{it})(J_\Phi\delta^{it}J_\Phi\underset{N}{_\beta\otimes_\alpha}J_\Phi\delta^{it}J_\Phi)W^*\\
&=(P^{it}\underset{\nu^o}{_\alpha\otimes_{\hat{\beta}}}P^{it})W(J_\Phi\delta^{it}J_\Phi\underset{N}{_\beta\otimes_\alpha}J_\Phi\delta^{it}J_\Phi)W^*\\
&=(P^{it}\underset{\nu^o}{_\alpha\otimes_{\hat{\beta}}}P^{it})(J_\Phi\delta^{it}J_\Phi\underset{N^o}{_\alpha\otimes_{\hat{\beta}}}1)W(1\underset{N}{_\beta\otimes_\alpha}J_\Phi\delta^{it}J_\Phi)W^*
\end{align*}
which is equal to :
\[(P^{it}\underset{\nu^o}{_\alpha\otimes_{\hat{\beta}}}P^{it})(J_\Phi\delta^{it}J_\Phi\underset{N^o}{_\alpha\otimes_{\hat{\beta}}}1)(J_{\widehat{\Phi}}\underset{\nu}{_\beta\otimes_\alpha}J_\Phi)W^*(1\underset{N^o}{_\alpha\otimes_{\hat{\beta}}}\delta^{it})W(J_{\widehat{\Phi}}\underset{\nu^o}{_\alpha\otimes_{\hat{\beta}}}J_\Phi)\]
and, therefore, to :
\[(P^{it}\underset{\nu^o}{_\alpha\otimes_{\hat{\beta}}}P^{it})(J_\Phi\delta^{it}J_\Phi\underset{N^o}{_\alpha\otimes_{\hat{\beta}}}1)(J_{\widehat{\Phi}}\underset{\nu}{_\beta\otimes_\alpha}J_\Phi)(\delta^{it}\underset{\nu}{_\beta\otimes_\alpha}\delta^{it})(J_{\widehat{\Phi}}\underset{\nu^o}{_\alpha\otimes_{\hat{\beta}}}J_\Phi)\]
or to :
\[(P^{it}\underset{\nu^o}{_\alpha\otimes_{\hat{\beta}}}P^{it})(J_\Phi\delta^{it}J_\Phi\underset{N^o}{_\alpha\otimes_{\hat{\beta}}}1)(\delta^{it}\underset{\nu^o}{_\alpha\otimes_{\hat{\beta}}}J_{\widehat{\Phi}}\delta^{it}J_{\widehat{\Phi}})=(\delta\Delta_{\widehat{\Phi}})^{it}\underset{\nu^o}{_\alpha\otimes_{\hat{\beta}}}\Delta_{\widehat{\Phi}}^{it}\]
which finishes the proof. \end{proof}

%%%%%psibarre4
\subsection{Proposition}
\label{psibarre4}
{\it Let $\gG=(N, M, \alpha, \beta, \Gamma, T, T', \nu)$ be a measured quantum groupoid, $(b,\ga)$ a weighted action of $\gG$ on a von Neumann algebra $A$, and $\psi$ a normal semi-finite faithful weight on $A$, lifted from $\nu^o$; then the von Neumann algebra $A\underset{N}{_b*_\alpha}\mathcal L(H)$ is a faithful right $N$-module in two different ways, using $1\underset{N}{_b\otimes_\alpha}\beta$, and $1\underset{N}{_b\otimes_\alpha}\hat{\beta}$; moreover, the weight $\underline{\psi}$ constructed in \ref{psibarre} is a lifted weight from $\nu$, using $1\underset{N}{_b\otimes_\alpha}\beta$, and, on the other hand, satisfies the $(1\underset{N}{_b\otimes_\alpha}\hat{\beta}, \gamma)$ property ; therefore, we can define a normal semi-finite faithful weight $\underline{\underline{\psi}}$ on $A\underset{N}{_b*_\alpha}\mathcal L(H)\underset{N}{_\beta*_\alpha}\mathcal L(H)$, and another normal semi-finite faithful weight $\underline{(\underline{\psi})}_\delta$ on $A\underset{N}{_b*_\alpha}\mathcal L(H)\underset{N}{_{\hat{\beta}}*_\alpha}\mathcal L(H)$. As in \ref{tildeTheta}, let us write, for any $Y$ in $\mathcal L(\gH\underset{\nu}{_b\otimes_\alpha}H\underset{\nu}{_{\hat{\beta}}\otimes_\alpha}H)$, 
\[\tilde{\Theta}(Y)=(1\underset{N}{_b\otimes_\alpha}W)^*(id\underset{N}{_b*_\alpha}\varsigma_N)(Y)(1\underset{N}{_b\otimes_\alpha}W)\]
which belongs to $\mathcal L(\gH\underset{\nu}{_b\otimes_\alpha}H\underset{\nu}{_\beta\otimes_\alpha}H)$. Then, we have :}
\[\underline{\underline{\psi}}\circ\tilde{\Theta}=\underline{(\underline{\psi})}_\delta\]

\begin{proof}
By definition, the weight $\underline{\underline{\psi}}$ is defined on $A\underset{N}{_b*_\alpha}\mathcal L(H)\underset{N}{_\beta*_\alpha}\mathcal L(H)$ by considering on $H_{\underline{\psi}}\underset{\nu}{_\beta\otimes_\alpha}H$ the spatial derivative :
\[\frac{d\underline{\underline{\psi}}}{d(\underline{\psi})^o}=\Delta_{\underline{\psi}}\underset{\nu}{_\beta\otimes_\alpha}\Delta_{\widehat{\Phi}}^{-1}\]
and, using \ref{psibarre3}, we therefore get, on $H\underset{\nu}{_\beta\otimes_a}H_\psi\underset{\nu}{_b\otimes_\alpha}H\underset{\nu}{_\beta\otimes_\alpha}H$, that :
\[\frac{d\underline{\underline{\psi}}}{d(\underline{\psi})^o}=\Delta_{\widehat{\Phi}}^{-1}\underset{\nu}{_\beta\otimes_a}\Delta_\psi\underset{\nu}{_b\otimes_\alpha}\Delta_{\widehat{\Phi}}^{-1}\underset{\nu}{_\beta\otimes_\alpha}\Delta_{\widehat{\Phi}}^{-1}\]
On the other hand, the weight $\underline{(\underline{\psi})}_\delta$ is defined on $A\underset{N}{_b*_\alpha}\mathcal L(H)\underset{N}{_{\hat{\beta}}*_\alpha}\mathcal L(H)$ by considering on $H_{\underline{\psi}}\underset{\nu}{_{\hat{\beta}}\otimes_\alpha}H=H\underset{\nu}{_\beta\otimes_a}H_\psi\underset{\nu}{_b\otimes_\alpha}H\underset{\nu}{_{\hat{\beta}}\otimes_\alpha}H$ the spatial derivative :
\[\frac{d\underline{(\underline{\psi})}_\delta}{d(\underline{\psi})^o}=\Delta_{\underline{\psi}}\underset{\nu}{_{\hat{\beta}}\otimes_\alpha}(\delta\Delta_{\widehat{\Phi}})^{-1}=\Delta_{\widehat{\Phi}}^{-1}\underset{\nu}{_\beta\otimes_a}\Delta_\psi\underset{\nu}{_b\otimes_\alpha}\Delta_{\widehat{\Phi}}^{-1}\underset{\nu}{_{\hat{\beta}}\otimes_\alpha}(\delta\Delta_{\widehat{\Phi}})^{-1}\]
from which we get, using \ref{deltaDelta} and the definition of $\tilde{\Theta}$ that :
\[\frac{d\underline{\underline{\psi}}}{d(\underline{\psi})^o}=(id\underset{N}{_\beta*_a}\tilde{\Theta})(\frac{d\underline{(\underline{\psi})}_\delta}{d(\underline{\psi})^o})\]
The weight $(\underline{\psi})^o$ is defined on $J_{\underline{\psi}}\pi_{\underline{\psi}}(A\underset{N}{_b*_\alpha}\mathcal L(H))J_{\underline{\psi}}$, which, using again \ref{psibarre3}, is equal to $\mathcal L(H)\underset{N}{_\beta*_a}A'\underset{N}{_b\otimes_\alpha}1_H$; we see, therefore, for $X\in \mathcal L(H)\underset{N}{_\beta*_a}A'$, that $(id\underset{N}{_\beta*_a}\tilde{\Theta})$ sends $X\underset{N}{_b\otimes_\alpha}1_H\underset{N}{_{\hat{\beta}}\otimes_\alpha}1_H$ on $X\underset{N}{_b\otimes_\alpha}1_H\underset{N}{_\beta\otimes_\alpha}1_H$, and leaves $(\underline{\psi})^o$ invariant. From which we deduce that :
\[\frac{d\underline{\underline{\psi}}\circ\tilde{\Theta}}{d(\underline{\psi})^o}=\frac{d\underline{(\underline{\psi})}_\delta}{d(\underline{\psi})^o}\]
from which we get the result. \end{proof}

%%%%%bidulaw
\section{Biduality of weights}
\label{bidualw}
In that chapter, following what had been done for locally compact quantum groups in \cite{Y4}, \cite{Y5}, and \cite{BV}, starting from an action $\ga$ of a measured quantum groupoid on a von Neumann algebra $A$, we define the Radon-Nikodym derivative of a lifted weight on $A$ with respect to this action (\ref{defder}); this operator is an $\ga$-cocycle (\ref{cocycle}), which measures, in a certain sense, how the weight $\psi$ behaves towards the action. In particular, we prove that this cocycle is equal to $1$ if and only if the weight is invariant by the action (\ref{thinv}, \ref{thinv2}). 
%%%%bidualw
\subsection{Theorem}
\label{thbidualw}
{\it Let $\gG=(N, M, \alpha, \beta, \Gamma, T, T', \nu)$ be a measured quantum groupoid, $(b,\ga)$ a weighted action of $\gG$ on a von Neumann algebra $A$, $\psi$ a normal semi-finite faithful weight on $A$ lifted from $\nu^o$; let $\tilde{\psi}$ be the dual weight on the crossed-product $A\rtimes_\ga\gG$, and let $\overline{\psi_\ga}$ be the normal semi-finite faithful weight on $A\underset{N}{_b*_\alpha}\mathcal L(H)$ obtained from the bidual weight $\tilde{\tilde{\psi}}$ and the isomorphism between $A\underset{N}{_b*_\alpha}\mathcal L(H)$ and the double crossed-product; let $\underline{\psi}$ be normal semi-finite faithful weight on $A\underset{N}{_b*_\alpha}\mathcal L(H)$ constructed in \ref{psibarre}. We have then :
\[(D\overline{\psi_\ga}:D\underline{\psi})_t=\Delta_{\tilde{\psi}}^{it}(\Delta_\psi^{-it}\underset{\nu}{_b\otimes_\alpha}\Delta_{\widehat{\Phi}}^{it})\]
Moreover, the unitaries $\Delta_{\tilde{\psi}}^{it}(\Delta_\psi^{-it}\underset{\nu}{_b\otimes_\alpha}\Delta_{\widehat{\Phi}}^{it})$ belong to $A\underset{N}{_b*_\alpha}(M\cap\beta(N)')$. }
\begin{proof}
We have $(D\overline{\psi_\ga}:D\underline{\psi})_t=
(\frac{d\overline{\psi_\ga}}{d\psi^o})^{it}(\frac{d\underline{\psi}}{d\psi^o})^{-it}$, from which we get the first result, by \ref{crossed} and \ref{psibarre2}. So, we get that the unitaries $\Delta_{\tilde{\psi}}^{it}(\Delta_\psi^{-it}\underset{\nu}{_b\otimes_\alpha}\Delta_{\widehat{\Phi}}^{it})$ belong to $A\underset{N}{_b*_\alpha}\mathcal L(H)$; let's take $x\in M'$; using \ref{corsigma3}, we have $\sigma_t^{\overline{\psi}_\ga}(1\underset{N}{_b\otimes_\alpha}x)=1\underset{N}{_b\otimes_\alpha}\Delta_{\widehat{\Phi}}^{-it}x\Delta_{\widehat{\Phi}}^{it}$, and, using \ref{psibarre2}, we get that $\sigma_t^{\underline{\psi}}(1\underset{N}{_b\otimes_\alpha}x)=1\underset{N}{_b\otimes_\alpha}\Delta_{\widehat{\Phi}}^{-it}x\Delta_{\widehat{\Phi}}^{it}$; therefore, we get that $(\frac{d\overline{\psi_\ga}}{d\psi^o})^{it}(\frac{d\underline{\psi}}{d\psi^o})^{-it}$ commutes with $1\underset{N}{_b\otimes_\alpha}x$, and, therefore, belongs to $A\underset{N}{_b*_\alpha}M$. 
\newline
Let $n\in N$; we have :
\[\sigma_t^{\underline{\psi}}(1\underset{N}{_b\otimes_\alpha}\beta(n))=1\underset{N}{_b\otimes_\alpha}\Delta_{\widehat{\Phi}}^{-it}\beta(n)\Delta_{\widehat{\Phi}}^{it}=1\underset{N}{_b\otimes_\alpha}\tau_{-t}(\beta(n))=1\underset{N}{_b\otimes_\alpha}\beta(\sigma_{-t}^\nu(n))\]
and, on the other hand :
\[\sigma_t^{\overline{\psi}_\ga}(1\underset{N}{_b\otimes_\alpha}\beta(n))=\sigma_t^{\overline{\psi}_\ga}(\ga(b(n)))=\ga(\sigma_t^\psi(b(n)))=\ga(b(\sigma_{-t}^\nu(n))=1\underset{N}{_b\otimes_\alpha}\beta(\sigma_{-t}^\nu(n))\]
which proves that both $\underline{\psi}$ and $\overline{\psi}_\ga$ are lifted weights from the weight $\nu^o$, and, therefore, that $(D\overline{\psi_\ga}:D\underline{\psi})_t$ belongs to $A\underset{N}{_b*_\alpha}\beta(N)'$, which finishes the proof. 
\end{proof}

%%%defder
\subsection{Definition}
\label{defder}
Let $\gG=(N, M, \alpha, \beta, \Gamma, T, T', \nu)$ be a measured quantum groupoid, $(b,\ga)$ a weighted action of $\gG$ on a von Neumann algebra $A$, $\psi$ a normal semi-finite faithful weight on $A$ lifted from $\nu^o$; we shall call the unitaries $(D\overline{\psi_\ga}:D\underline{\psi})_t\in A\underset{N}{_b*_\alpha}(M\cap\beta(N)')$ the Radon-Nikodym derivative of the weight $\psi$ with respect to the action $(b, \ga)$, and denote it, for simplification, $(D\psi\circ\ga:D\psi)_t$, following the notations of (\cite{BV}, 10.2). 

%%%%cocycle
\subsection{Theorem}
\label{cocycle}
{\it Let $\gG=(N, M, \alpha, \beta, \Gamma, T, T', \nu)$ be a measured quantum groupoid, $(b,\ga)$ a weighted action of $\gG$ on a von Neumann algebra $A$, $\psi$ a normal semi-finite faithful weight on $A$ lifted from $\nu^o$; the Radon-Nikodym derivative  $(D\psi\circ\ga:D\psi)_t$ introduced in \ref{defder} is a $\ga$-cocycle, i.e., we have :}
\[(id\underset{N}{_b*_\alpha}\Gamma)((D\psi\circ\ga:D\psi)_t)=(\ga\underset{N}{_b*_\alpha}id)((D\psi\circ\ga:D\psi)_t)((D\psi\circ\ga:D\psi)_t)\underset{N}{_\beta\otimes_\alpha}1)\]

\begin{proof}
For all $t\in\mathbb{R}$, $(\ga\underset{N}{_b*_\alpha}id)((D\psi\circ\ga:D\psi)_t)$ belongs to $A\underset{N}{_b*_\alpha}M\underset{N}{_\beta*_\alpha}M$, and the operator $\underline{\ga}((D\psi\circ\ga:D\psi)_t)=\tilde{\Theta}^{-1}(\ga\underset{N}{_b*_\alpha}id)((D\psi\circ\ga:D\psi)_t)$ belongs to $A\underset{N}{_b*_\alpha} M\underset{N}{_{\hat{\beta}}*_\alpha}M$ (where $\tilde{\Theta}$ had been defined in \ref{psibarre4}) . 
\newline
We have, using successively \ref{crossed}, \ref{corthgamma} and \ref{tildeTheta}(iii) :

\[\underline{\ga}((D\psi\circ\ga:D\psi)_t)
=
\underline{\ga}((D\overline{\psi_\ga}:D\underline{\psi})_t)
=
(D\overline{(\overline{\psi_\ga})_{\underline{\ga}}}:D\overline{(\underline{\psi})_{\underline{\ga}}})_t
=
(D\underline{(\overline{\psi_\ga})}_\delta :D\underline{(\overline{\psi_\ga})}\circ\tilde{\Theta})_t\]
On the other hand, using successively \ref{corunderline}(ii) and \ref{psibarre4} :
\begin{align*}
\tilde{\Theta}^{-1}((D\psi\circ\ga:D\psi)_t)\underset{N}{_\beta\otimes_\alpha}1)
&=
\tilde{\Theta}^{-1}((D\overline{\psi_\ga}:D\underline{\psi})_t\underset{N}{_\beta\otimes_\alpha}1)\\
&=
\tilde{\Theta}^{-1}(D\underline{\overline{\psi_\ga}}:D\underline{\underline{\psi}})_t)\\
&=
(D\underline{\overline{\psi_\ga}}\circ\tilde{\Theta}:D
\underline{\underline{\psi}}\circ\tilde{\Theta})_t\\
&=(D\underline{\overline{\psi_\ga}}\circ\tilde{\Theta}:D\underline{(\underline{\psi})}_\delta)_t
\end{align*}
and, therefore, we get that :
\begin{multline*}
\tilde{\Theta}^{-1}[(\ga\underset{N}{_b*_\alpha}id)((D\psi\circ\ga:D\psi)_t)((D\psi\circ\ga:D\psi)_t)\underset{N}{_\beta\otimes_\alpha}1)]\\
=
\underline{\ga}((D\psi\circ\ga:D\psi)_t)\tilde{\Theta}^{-1}((D\psi\circ\ga:D\psi)_t)\underset{N}{_\beta\otimes_\alpha}1)
\end{multline*}
is equal, using \ref{cor2}(ii), to :
\begin{align*}
(D\underline{(\overline{\psi_\ga})}_\delta :D\underline{(\overline{\psi_\ga})}\circ\tilde{\Theta})_t(D\underline{(\overline{\psi_\ga})}\circ\tilde{\Theta}:D\underline{(\underline{\psi})}_\delta)_t
&=
(D\underline{(\overline{\psi_\ga})}_\delta :D\underline{(\underline{\psi})}_\delta)_t\\
&=
(D\overline{\psi_\ga}:D\underline{\psi})_t\underset{N}{_{\hat{\beta}}\otimes_\alpha}1\\
&=(D\psi\circ\ga:D\psi)_t\underset{N}{_{\hat{\beta}}\otimes_\alpha}1
\end{align*}
from which we get that :
\begin{align*}
(\ga\underset{N}{_b*_\alpha}id)((D\psi\circ\ga:D\psi)_t)((D\psi\circ\ga:D\psi)_t)\underset{N}{_\beta\otimes_\alpha}1)
&=
\tilde{\Theta}((D\psi\circ\ga:D\psi)_t\underset{N}{_{\hat{\beta}}\otimes_\alpha}1)\\
&=
(id\underset{N}{_b*_\alpha}\Gamma)((D\psi\circ\ga:D\psi)_t)
\end{align*}
which is the result. \end{proof}

%%%%exlcg
\subsection{Example}
\label{exlcg}
Let $\bf{G}$ be a locally compact quantum group, and $\ga$ an action of $\bf{G}$ on a von Neumann algebra $A$; then this result had been obtained in (\cite{Y4}, 4.8 and \cite{Y5}, 3.7 and \cite{BV}, 10.3).

%%%exgd2
\subsection{Example}
\label{exgd2}
Let $\mathcal G$ be a measured groupoid; let us use all the notations introduced in \ref{exgd}. Let $(\ga)_{g\in\mathcal G}$ be an action of $\mathcal G$ on a von Neumann algebra $A=\int_{\mathcal G^{(0)}}^\oplus A^x d\nu(x)$, and $\psi=\int_{\mathcal G^{(0)}}^\oplus \psi^x d\nu(x)$ a normal semi-finite faithful weight on $A$. Then, the Radon-Nikodym derivative of $\psi$ with respect to the action $\ga$, is, using (\cite{Y3}, 2.6), given by :
\[(D\psi\circ\ga:D\psi)_t=\int_\mathcal G^\oplus (D\psi^{r(g)}:D\psi^{s(g)}\circ \ga_{g^{-1}})_t d\nu(g)\]
which is acting on $\int_\mathcal G^\oplus H_{\psi^{r(g)}} d\mu(g)=H_\psi\underset{\nu}{_b\otimes_{r_\mathcal G}}L^2(\mathcal G, \mu)$. 

%%%%inv
\subsection{Definition}
\label{inv}
Let $(b, \ga)$ an action of a measured quantum groupoid $\gG$ on a von Neumann algebra $A$. A normal semi-finite faithful weight $\psi$ on $A$ will be said invariant by $\ga$ if, for all $\eta\in D(_\alpha H, \nu)\cap D(H_\beta, \nu^o)$ and $x\in\gN_\psi$, we have :
\[\psi[(id\underset{N}{_b*_\alpha}\omega_\eta)\ga(x^*x)]=\|\Lambda_\psi (x)\underset{\nu^o}{_a\otimes_\beta}\eta\|^2\]
We shall always suppose that such weights bear the density property, defined in \ref{action}, as for $\delta$-invariant weights. 

%%%%thinv
\subsection{Theorem}
\label{thinv}
{\it Let $(b, \ga)$ an action of a measured groupoid $\gG$ on a von Neumann algebra $A$, $\psi$ a normal semi-finite faithful weight on $A$, invariant by $\ga$ in the sense of \ref{inv}, and bearing the density property, as defined in \ref{action}. Then, let $(e_i)_{i\in I}$ be an $(\alpha, \nu)$-orthogonal basis of $H$, $x\in\gN_\psi$, $\eta\in D(_\alpha H, \nu)\cap D(H_\beta, \nu^o)$ :
\newline
(i) for any $\xi\in D(_\alpha H, \nu)$, $(id\underset{N}{_b*_\alpha}\omega_{\eta, \xi})\ga(x)$ belongs to $\gN_\psi$; 
\newline
(ii) the sum $\sum_i\Lambda_\psi((id\underset{N}{_b*_\alpha}\omega_{\eta, e_i})\ga(x)\underset{\nu}{_b\otimes_\alpha}e_i$ is strongly converging; its limit does not depend upon the choice of the $(\alpha, \nu)$-othogonal basis of $H$, and allow us to define an isometry $V'_\psi$ from $H_\psi\underset{\nu^o}{_a\otimes_\beta}H$ to $H_\psi\underset{\nu}{_b\otimes_\alpha}H$ such that :
\[V'_\psi(\Lambda_\psi(x)\underset{\nu^o}{_a\otimes_\beta}\eta)=\sum_i\Lambda_\psi((id\underset{N}{_b*_\alpha}\omega_{\eta, e_i})\ga(x)\underset{\nu}{_b\otimes_\alpha}e_i\]
(iii) we have :
\[\Lambda_\psi((id\underset{N}{_b*_\alpha}\omega_{\eta, \xi})\ga(x))=(id*\omega_{\eta, \xi})(V'_\psi)\lambda_\psi (x)\]
(iv) for any $y\in A$, $z\in M'$, $n\in N$, we have :
\[\ga(y)V'_\psi=V'_\psi (y\underset{N^o}{_a\otimes_\beta}1)\]
\[(1\underset{N}{_b\otimes_\alpha}z)V'_\psi=V'_\psi(1\underset{N^o}{_a\otimes_\beta}z)\]
\[(a(n)\underset{N}{_b\otimes_\alpha}1)V'_\psi=V'_\psi(1\underset{N^o}{_a\otimes_\beta}\alpha(n))\]
\[(1\underset{N}{_b\otimes_\alpha}\beta(n))V'_\psi=V'_\psi(b(n)\underset{N^o}{_a\otimes_\beta}1)\]
\[(1\underset{N}{_b\otimes_\alpha}\hat{\beta}(n))V'_\psi=V'_\psi(1\underset{N^o}{_a\otimes_\beta}\hat{\beta}(n))\]
(v) the operator $V'_\psi$ is a unitary; moreover, it is a copresentation of $\gG$ on $_a(H_\psi)_b$ which implements $\ga$; 
\newline
(vi) we have :
\[V'_\psi(\Delta_\psi^{it}\underset{N^o}{_a\otimes_\beta}\Delta_{\widehat{\Phi}}^{-it})=(\Delta_\psi^{it}\underset{N}{_b\otimes_\alpha}\Delta_{\widehat{\Phi}}^{-it})V'_\psi\]
Moreover, the weight $\psi$ is lifted from $\nu^o$; more precisely, there exists a normal faithful semi-finite operator-valued weight $\gT$ from $A$ onto $b(N)$ such that $\psi=\nu^o\circ b^{-1}\circ \gT$, and, for all $x\in\gN_{\gT}\cap\gN_\psi$, we have : 
\[(\gT\underset{N}{_b*_\alpha}id)\ga(x^*x)=1\underset{N}{_b\otimes_\alpha}\beta\circ b^{-1}\gT(x^*x)=\ga(\gT(x^*x))\]
\[(\psi\underset{\nu}{_b*_\alpha}id)\ga(x^*x)=\beta\circ b^{-1}\gT(x^*x)\]
(vii) we have :
\[\ga(\sigma_t^\psi(y))=(\sigma_t^\psi\underset{N}{_b*_\alpha}\tau_t)\ga(y)\]
(viii) the standard implementation $U^{\ga}_\psi$ is equal to $V'_\psi$; 
\newline
(ix) the dual weight satisfies $\Delta_{\tilde{\psi}}^{it}=\Delta_\psi^{it}\underset{N}{_b\otimes_\alpha}\Delta_{\widehat{\Phi}}^{-it}$; 
\newline
(x) the Radon-Nikodym derivative $(D\psi\circ\ga:D\psi)_t$ is equal to $1$. }
\begin{proof}
Result (i) is identical to (\cite{E5}, 8.3(i)), and (ii) is similar to (\cite{E5}, 8.3(ii) and 8.4(i)); the proof of (iii) is similar to the proof of (\cite{E5}, 8.4(ii) and (iii)), and the proof of (iv) is similar (and somehow simpler) to the proof of (\cite{E5}, 8.4(iv) and (v)). Now result (v) is obtained in a similar way to (\cite{E5}, 8.5 and 8.6); by similar calculations to (\cite{E5}, 8.7 and 8.8(i)), we obtain that, for all $t\in\mathbb{R}$, we have $\sigma_t^\psi(b(n))=b(\sigma_{-t}^\nu(n))$, which gives the existence of a normal faithful semi-finite operator-valued weight $\gT$ from $A$ onto $b(N)$ such that $\psi=\nu^o\circ b^{-1}\circ \gT$. For any $x\in\gN_\psi\cap\gN_{\gT}$, the vector $\Lambda_\psi(x)$ belongs to $\mathcal D(_\alpha H, \nu)$, and we have, for any $\eta\in H$ :
\[\|\Lambda_\psi (x)\underset{\nu^o}{_a\otimes_\beta}\eta\|^2=(\beta\circ b^{-1}\gT(x^*x)\eta|\eta)\]
So, using the density property and \ref{inv}, we get, for all $x\in \gN_\psi\cap\gN_{\gT}$, that :
\[(\psi\underset{\nu}{_b*_\alpha}id)\ga(x^*x)=\beta\circ b^{-1}\gT(x^*x)\]
and, therefore, that :
\[(\gT\underset{N}{_b*_\alpha}id)\ga(x^*x)=1\underset{N}{_b\otimes_\alpha}\beta\circ b^{-1}\gT(x^*x)=\ga(\gT(x^*x)\]
we finish the proof of (vi) in a similar way to (\cite{E5}, 8.8(ii)). Then (vii) is a straightforward corollary of (vi) and (v), and (viii) and (ix) are obtained in a similar way to \ref{Uinv}(i) and (ii). As $\Delta_{\tilde{\psi}}=\frac{d\overline{\psi_\ga}}{d\psi^o}$ (\cite{E5} 13.6) and $\Delta_\psi\underset{N}{_b\otimes_\alpha}\Delta_{\widehat{\Phi}}^{-1}=\frac{d\underline{\psi}}{d\psi^o}$ by \ref{psibarre}(ii), we infer from (ix) that $\overline{\psi_\ga}=\underline{\psi}$, which, by \ref{defder}, finishes the proof. \end{proof}

%%%%corinv
\subsection{Corollary}
\label{corinv}
{\it  Let $(b, \ga)$ be an action of a measured quantum groupoid $\gG$ on a von Neumann algebra $A$; let $\psi_1$, $\psi_2$ be two invariant normal faithful semi-finite weights on $A$, as defined in \ref{inv}, and let us suppose that both $\psi_1$ and $\psi_2$ bear the density property, as defined in \ref{action}. Then, for all $t\in\mathbb{R}$, $(D\psi_1:D\psi_2)_t$ belongs to $A^\ga$. }
\begin{proof}
The proof is similar to (\cite{E5}, 8.11). \end{proof}

%%%thinv2
\subsection{Theorem}
\label{thinv2}
{\it Let $(b, \ga)$ be a weighted action of a measured quantum groupoid $\gG$ on a von Neumann algebra $A$, and $\psi$ a normal semi-finite faithful weight on $A$, lifted from $\nu^o$. If the Radon-Nikodym derivative $(D\psi\circ\ga:D\psi)_t$ is equal to $1$, then the weight $\psi$ is invariant by $\ga$ in the sense of \ref{inv}. }

\begin{proof}
Let $\xi\in D(_\alpha H, \nu)\cap D(H_\beta, \nu^o)\cap\mathcal D(\Delta_{\widehat{\Phi}}^{-1/2})$ such that $\Delta_{\widehat{\Phi}}^{-1/2}\xi$ belongs to $D(_\alpha H, \nu)$; let us remark first that if $y$ belongs to $\gN_{\widehat{\Phi}}\cap\gN_{\widehat{\Phi}}^*\cap\gN_{\hat{T}}\cap\gN_{\hat{T}}^*$, and is analytic with respect to $\sigma_t^{\widehat{\Phi}}$, and such that $\sigma_z(x)$ belongs to 
$\gN_{\widehat{\Phi}}\cap\gN_{\widehat{\Phi}}^*\cap \gN_{\hat{T}}\gN_{\hat{T}}^*$, for all $z\in\mathbb{C}$, then $\Lambda_{\widehat{\Phi}}(z)$ satisfies all those conditions, and this gives that the set of such elements $\xi$ is dense in $H$. 
\newline
Let $\eta$ be in $D(_\alpha H, \nu)\cap \mathcal D(\Delta_{\widehat{\Phi}}^{-1/2})$ such that $\Delta_{\widehat{\Phi}}^{-1/2}\eta$ belongs to $D(_\alpha H, \nu)$, and $x\in \gN_\psi$, analytic with respect to $\psi$, such that $\sigma_{-i/2}(x^*)$ belongs to $\gN_\psi$. Then, we have, using \ref{psibarre3}(i) applied to $\underline{\nu^o}$ :
\begin{multline*}
((\psi\underset{\nu}{_b*_\alpha}id)\ga (x^*x)\xi\underset{N^o}{_\alpha\otimes_\beta}J_{\widehat{\phi}}\Delta_{\widehat{\Phi}}^{-1/2}\eta|\xi\underset{N^o}{_\alpha\otimes_\beta}J_{\widehat{\phi}}\Delta_{\widehat{\Phi}}^{-1/2}\eta)=\\
((\psi\underset{\nu}{_b*_\alpha}id)\ga (x^*x)\Lambda_{\underline{\nu^o}}(\theta^{\alpha, \nu}(\xi, \eta))|\Lambda_{\underline{\nu^o}}(\theta^{\alpha, \nu}(\xi, \eta))=\\
\underline{\nu^o}(\theta^{\alpha, \nu}(\xi, \eta)^*(\psi\underset{\nu}{_b*_\alpha}id)\ga (x^*x)\theta^{\alpha, \nu}(\xi, \eta)
\end{multline*}
which is equal, using \ref{psibarre2} and \ref{ex2}, to :
\[\underline{\psi}(1\underset{N}{_b\otimes_\alpha}\theta^{\alpha, \nu}(\xi, \eta))^*\ga(x^*x)(1\underset{N}{_b\otimes_\alpha}\theta^{\alpha, \nu}(\xi, \eta)))\]
By hypothesis, as $\overline{\psi_\ga}=\underline{\psi}$ by \ref{defder}, we get, using \ref{action} that $\sigma_t^{\underline{\psi}}(\ga(x))=\sigma_t^{\overline{\psi_\ga}}(\ga(x))=\ga(\sigma_t^\psi(x))$. Moreover, we can write, thanks to the hypothesis and to \ref{psibarre3} applied to $\underline{\nu^o}$ :
\[J_{\underline{\nu^o}}\Lambda_{\underline{\nu^o}}(\theta^{\alpha, \nu}(\xi, \eta))=J_{\widehat{\Phi}}\xi\underset{\nu}{_\beta\otimes_\alpha}\Delta_{\widehat{\Phi}}^{-1/2}\eta=\Lambda_{\underline{\nu^o}}(\theta^{\alpha, \nu}(\Delta_{\widehat{\Phi}}^{-1/2}\eta, \Delta_{\widehat{\Phi}}^{1/2}\xi))\]
from which we get that $[\ga(x)(1\underset{N}{_b\otimes_\alpha}\theta^{\alpha, \nu}(\xi, \eta))]^*$ belongs to $\mathcal D(\sigma_{-i/2}^{\underline{\psi}})$, and, therefore, that:
\[((\psi\underset{\nu}{_b*_\alpha}id)\ga (x^*x)\xi\underset{\nu^o}{_\alpha\otimes_\beta}J_{\widehat{\phi}}\Delta_{\widehat{\Phi}}^{-1/2}\eta|\xi\underset{\nu^o}{_\alpha\otimes_\beta}J_{\widehat{\phi}}\Delta_{\widehat{\Phi}}^{-1/2}\eta)\]
is equal to :
\[\|\Lambda_{\underline{\psi}}(\sigma_{-i/2}^{\underline{\psi}}([\ga(x)(1\underset{N}{_b\otimes_\alpha}\theta^{\alpha, \nu}(\xi, \eta))]^*)\|^2=
\|\Lambda_{\underline{\psi}}((1\underset{N}{_b\otimes_\alpha}\theta^{\alpha, \nu}(\Delta_{\widehat{\Phi}}^{-1/2}\eta, \Delta_{\widehat{\Phi}}^{1/2}\xi))\ga(\sigma_{-i/2}^\psi(x^*)))\|^2\]
which, thanks again to the hypothesis and to \ref{psia}, is equal to :
\begin{align*}
\|\Lambda_{\psi}(\sigma_{-i/2}^\psi(x^*))\underset{\nu}{_b\otimes_\alpha}J_{\widehat{\Phi}}\xi\underset{\nu}{_\beta\otimes_\alpha}\Delta_{\widehat{\Phi}}^{-1/2}\eta\|^2
&=
\|J_\psi\Lambda_\psi(x)\underset{\nu}{_b\otimes_\alpha}J_{\widehat{\Phi}}\xi\underset{\nu}{_\beta\otimes_\alpha}\Delta_{\widehat{\Phi}}^{-1/2}\eta\|^2\\
&=
\|\Lambda_\psi(x)\underset{\nu^o}{_a\otimes_\beta}\xi\underset{\nu^o}{_\alpha\otimes_\beta}J_{\widehat{\Phi}}\Delta_{\widehat{\Phi}}^{-1/2}\eta\|^2
\end{align*}
So, finally, we get the equality :
\[((\psi\underset{\nu}{_b*_\alpha}id)\ga (x^*x)\xi\underset{\nu^o}{_\alpha\otimes_\beta}J_{\widehat{\phi}}\Delta_{\widehat{\Phi}}^{-1/2}\eta|\xi\underset{\nu^o}{_\alpha\otimes_\beta}J_{\widehat{\phi}}\Delta_{\widehat{\Phi}}^{-1/2}\eta)=
\|\Lambda_\psi(x)\underset{\nu^o}{_a\otimes_\beta}\xi\underset{\nu^o}{_\alpha\otimes_\beta}J_{\widehat{\Phi}}\Delta_{\widehat{\Phi}}^{-1/2}\eta\|^2\]
which, by continuity, remains true for any $x\in\gN_\psi$ and $\xi\in D(_\alpha H, \nu)\cap D(H_\beta, \nu^o)$; from which we infer that :
\begin{multline*}
((\psi\underset{\nu}{_b*_\alpha}id)\ga (x^*x)\alpha(<J_{\widehat{\Phi}}\Delta_{\widehat{\Phi}}^{-1/2}\eta, J_{\widehat{\Phi}}\Delta_{\widehat{\Phi}}^{-1/2}\eta>_{\beta, \nu^o})\xi|\xi)=\\
(\Lambda_\psi(x)\underset{\nu^o}{_a\otimes_\beta}\alpha(<J_{\widehat{\Phi}}\Delta_{\widehat{\Phi}}^{-1/2}\eta, J_{\widehat{\Phi}}\Delta_{\widehat{\Phi}}^{-1/2}\eta>_{\beta, \nu^o})\xi|\Lambda_\psi(x)\underset{\nu^o}{_a\otimes_\beta}\xi)
\end{multline*}
from which, by density of the elements of the form $<J_{\widehat{\Phi}}\Delta_{\widehat{\Phi}}^{-1/2}\eta, J_{\widehat{\Phi}}\Delta_{\widehat{\Phi}}^{-1/2}\eta>_{\beta, \nu^o}$ in $N^+$, we get, for any $n\in N^+$ :
\[((\psi\underset{\nu}{_b*_\alpha}id)\ga (x^*x)\alpha(n)\xi|\xi)=(\Lambda_\psi(x)\underset{\nu^o}{_a\otimes_\beta}\alpha(n)\xi|\Lambda_\psi(x)\underset{\nu^o}{_a\otimes_\beta}\xi)\]
from which we get the result, by density of $D(_\alpha H, \nu)\cap D(H_\beta, \nu^o)$. \end{proof}

%%%propcocycle
\subsection{Proposition}
\label{propcocycle}
{\it Let $\gG$ be a measured quantum groupoid, $(b,\ga)$ a weighted action of $\gG$ on a von Neumann algebra $A$, $\psi_1$ and $\psi_2$ two normal semi-finite faithful weights on $A$, lifted from $\nu^o$, and $(D\psi_1\circ\ga :D\psi_1)_t$, $(D\psi_2\circ\ga: D\psi_2)_t$ their Radon-Nikodym derivatives with respect to the action $(b, \ga)$, as defined in \ref{defder}. Then, the Radon-Nikodym derivative $(D\psi_1:D\psi_2)_t$ belongs to $A\cap b(N)'$, and we have, for all $t\in\mathbb{R}$ :}
\[(D\psi_2\circ\ga: D\psi_2)_t=\ga((D\psi_2:D\psi_1)_t)(D\psi_1\circ\ga :D\psi_1)_t((D\psi_2:D\psi_1)_t^*\underset{N}{_b\otimes_\alpha}1)\]

\begin{proof}
As $\psi_1$ and $\psi_2$ are lifted weights from $\nu$, we get that $(D\psi_1:D\psi_2)_t$ belongs to $A\cap b(N)'$ by (\cite{T}, 4.22(iii)); moreover, we have :
\[(D\overline{\psi_{2\ga}}:D\underline{\psi_2})_t=(D\overline{\psi_{2\ga}}:D\overline{\psi_{1\ga}})_t(D\overline{\psi_{1\ga}}:D\underline{\psi_1})_t(D\underline{\psi_1}:D\underline{\psi_2})_t\]
from which we get the result, using \ref{action}, \ref{defder} and \ref{corunderline}(ii). \end{proof}

%%%%corcocycle
\subsection{Corollary}
\label{corcocycle}
{\it Let $\gG$ be a measured quantum groupoid, $(b,\ga)$ a weighted action of $\gG$ on a von Neumann algebra $A$; then, are equivalent :
\newline
(i) there exists a normal semi-finite faithful weight on $A$, which is invariant and bears the density condition; 
\newline
(ii) there exists a normal semi-finite faithful weight $\psi$ on $A$, lifted from $\nu^o$, and a $\sigma_t^\psi$-cocycle $u_t$ on $A\cap b(N)'$ such that $(D\psi\circ\ga :D\psi)_t=\ga(u_t^*)(u_t\underset{N}{_b\otimes_\alpha}1)$; 
\newline
(iii) for any normal semi-finite faithful weight $\psi$ on $A$, lifted from $\nu^o$, there exists a $\sigma_t^\psi$-cocycle $u_t$ on $A\cap b(N)'$ such that $(D\psi\circ\ga :D\psi)_t=\ga(u_t^*)(u_t\underset{N}{_b\otimes_\alpha}1)$.}

\begin{proof}
Let suppose (i), and let $\varphi$ be an invariant weight on $A$, bearing the density condition; then, by \ref{thinv}(vi), the weight is lifted, and, if $\psi$ is any another lifted weight on $A$, $u_t=(D\varphi :D\psi)_t$ is a $\sigma_t^\psi$-cocycle in $A\cap b(N)'$ by (\cite{T}, 4.22(iii)); moreover, using \ref{propcocycle}, we get (iii). 
\newline
Conversely, if we suppose (ii), there exists a normal semi-finite faithful weight $\varphi$ on $A$ such that $u_t=(D\varphi:D\psi)_t$; as $\psi$ is lifted, and $u_t$ belongs to $A\cap b(N)'$, we know, using (\cite{T}, 4.22(iii)), that $\varphi$ is lifted, too. Using now \ref{propcocycle}, we get that $(D\varphi\circ\ga:D\varphi)_t=1$, which, thanks to \ref{thinv2}, gives the result. \end{proof}

%%%%%%bibli

\end{document}